\documentclass[11pt]{article}
\usepackage{amsmath,amsbsy,amsfonts}
\usepackage[frenchb]{babel}
\newcommand{\nl}[0]{\vskip 0.2cm{}\noindent}
\begin{document}

\title{Paquets d'Arthur pour les groupes classiques; point de vue combinatoire}
\author{C. M{\oe}glin\\ Institut de Math\'ematiques de Jussieu, CNRS Paris\\ moeglin@math.jussieu.fr}

\date{}
\maketitle
Le but de cet article est de terminer la description combinatoire des paquets d'Arthur pour des groupes classiques p-adiques.  Les paquets sont associ\'es \`a des morphismes, $\psi$, du groupe $W_{F}\times SL(2,{\mathbb C}) \times SL(2,{\mathbb C})$ dans le groupe dual (au sens de Langlands) de $G$, continus et born\'es sur $W_{F}$ et alg\'ebriques sur les 2 copies de $SL(2,{\mathbb C})$. Quand $G$ est un groupe symplectique ou orthogonal, le groupe dual a une repr\'esentation naturelle dans un groupe lin\'eaire de la forme $GL(n^*, {\mathbb C})$; en composant $\psi$ avec cette repr\'esentation naturelle, on trouve un morphisme \`a valeurs dans $GL(n^*,{\mathbb C})$. Gr\^ace \`a la correspondance locale de Langlands, on sait associer \`a $\psi$ une repr\'esentation irr\'eductible unitaire de $GL(n^*,F)$ qui est temp\'er\'ee exactement quand $\psi$ est trivial sur la 2e copie de $SL(2,{\mathbb C})$. Notons $\pi_{GL}(\psi)$ cette repr\'esentation.

Cet article fait suite \`a \cite{paquetsdiscrets} o\`u l'on a donn\'e une d\'efinition des repr\'esentations dans un paquet d'Arthur en supposant que le centralisateur du morphisme dans le $L$-groupe est un groupe fini et  pos\'e les d\'efinitions g\'en\'erales. Le point est ici d'\'etudier la r\'eductibilit\'e de certaines induites. On peut expliquer  ce que nous avons en vue de la fa\c{c}on suivante; prenons $\pi_{0}$ une composante locale d'une forme automorphe de carr\'e int\'egrable pour un groupe de m\^eme type que $G$; fixons $\rho$ une repr\'esentation cuspidale irr\'eductible d'un groupe lin\'eaire. On suppose que $\rho\simeq \rho^*$ et soient $a,b$ des entiers. On d\'efinit $Sp(b,St(a,\rho))$ comme l'unique sous-module irr\'eductible de l'induite:
$$
St(\rho,a)\vert \vert^{-(b-1)/2} \times \cdots \times St(\rho,a)\vert \,\vert^{(b-1)/2},$$o\`u $St(\rho,a)$ est la repr\'esentation de Steinberg associ\'ee \`a $\rho$ et $a$ et o\`u $\vert\,\vert$ est (comme dans tout ce travail) la valeur absolue du d\'eterminant. Ici $Sp$ est un diminutif de Speh analogue \`a $St$ pour tenir compte du fait que l'importance d'un point de vue global de ces repr\'esentations avaient \'et\'e mise en \'evidence par Birgit Speh. Ces repr\'esentations sont exactement les repr\'esentations des groupes lin\'eaires associ\'ee \`a des repr\'esentations irr\'eductibles (born\'ees, alg\'ebriques) de $W_{F}\times SL(2,{\mathbb C})\times SL(2,{\mathbb C})$ via la correspondance de Langlands maintenant d\'emontr\'ee (cf. \cite{harris}, \cite{henniart}, \cite{zelevinsky}). On note $\rho\otimes [a]\otimes [b]$ la repr\'esentation irr\'eductible de $W_{F}\times SL(2,{\mathbb C})\times SL(2,{\mathbb C})$ correspondant \`a $Sp(b,St(\rho,a))$.

On veut d\'ecomposer l'induite $Sp(b,St(\rho,a))\times \pi_{0}$; par unitarit\'e on sait que cette induite est semi-simple et on veut montrer que sa longueur est inf\'erieure ou \'egale \`a $inf(a,b)+1$ cette borne pouvant \^etre atteinte. Prenant ensuite un facteur direct, on veut aussi montrer que si on r\'einduit avec plusieurs copies de $Sp(b,St(\rho,a))$ on garde alors l'irr\'eductibilit\'e. 

Ne sachant rien sur $\pi_{0}$, on ne pourra pas faire grand chose, donc ce que l'on fait ici est de d\'emontrer ces r\'esultats pour $\pi_{0}$ dans un paquet d'Arthur tel que construit dans \cite{paquetsdiscrets}. On s'attend \`a avoir l'irr\'eductibilit\'e de cette induite si $\rho\otimes [a]\otimes [b]$ ne se factorise pas par un groupe de m\^eme type que $^LG$ et ceci est d\'emontr\'e en \ref{inductionfin}. Dans le cas oppos\'e, le r\'esultat est d\'emontr\'e en \ref{inductionsuite}. Signalons qu'une des difficult\'es annexes que l'on rencontre,  est que l'on ne sait pas que les repr\'esentations que nous consid\'erons sont unitaires et il faut donc aussi d\'emontrer que les induites sont semi-simples. On note $\psi_{0}$ la repr\'esentation de $W_{F}\times SL(2,{\mathbb C})\times SL(2,{\mathbb C})$ associ\'ee \`a $\pi_{0}$. La d\'emonstration peut se r\'esumer ainsi: en \cite{paquetsdiscrets}, on a d\'ecrit les repr\'esentations qui doivent \^etre associ\'ees au morphisme $\rho\otimes [a]\otimes [b]\oplus \rho \otimes [a]\otimes [b]\oplus \psi_{0}$. Comme on a fix\'e $\pi_{0}$ dans le paquet associ\'e \`a $\psi_{0}$, nous avons une limitation sur les param\`etres pouvant intervenir et il faut d\'emontrer que la somme des repr\'esentations associ\'ees \`a ces param\`etres d\'ecompose effectivement l'induite. Ce n'est en fait pas difficile, si l'on admet que $\pi_{0}$ est unitaire, c'est m\^eme assez facile. C'est donc un r\'esultat assez pr\'ecis qui est d\'emontr\'e ici, le seul ennui est que l'on ne sait toujours pas pour quels param\`etres la repr\'esentation associ\'ee n'est pas la fausse repr\'esentation $0$.

\

L'article commence par rappeler les d\'efinitions g\'en\'erales de \cite{paquetsdiscrets} et par donner quelques propri\'et\'es g\'en\'erales des repr\'esentations dans les paquets d'Arthur: une param\'etrisation de leur d\'ecomposition en composantes irr\'eductibles et des propri\'et\'es de leurs modules de Jacquet. Puis on d\'emontre la d\'ecomposition des induites. A la fin en \ref{conclusion}, on donne le cas le plus g\'en\'eral que je connais o\`u la d\'ecomposition ne fait pas intervenir cette incertitude sur la nullit\'e de certaines composantes. Je n'ai pas donn\'e d'exemles o\`u la longueur de l'induite est strictement comprise entre $1$ et $inf(a,b)+1$;  il y en a s\^urement mais ils sont aussi s\^urement techniquement p\'enibles \`a construire.

 \section{Paquets g\'en\'eraux\label{paquetsgeneraux}}
 \subsection{Blocs de Jordan\label{blocsdejordan}}
 Les paquets de repr\'esentations tels que d\'efinis par Arthur sont associ\'es \`a des morphismes, $\psi$, de $W_{F}\times SL(2,{\mathbb C}) \times SL(2,{\mathbb C})$ dans le $L$-groupe de $G$, continues, born\'ees sur $W_{F}$ et alg\'ebrique sur les \'e copies de $SL(2,{\mathbb C})$. Pour classifier \`a conjugaison pr\`es, on prolonge $\psi$ en une repr\'esentation gr\^ace \`a l'inclusion naturelle du $L$-groupe dans le groupe lin\'eaire correspondant \`a la repr\'esentation standard et on d\'ecompose cette repr\'esentation en composantes irr\'eductibles. 
 
 Une repr\'esentation irr\'eductible de $W_{F}\times SL(2,{\mathbb C})\times SL(2,{\mathbb C})$ est la donn\'ee d'une repr\'esentation irr\'eductible, $\rho$ de $W_{F}$ dans un groupe lin\'eaire dont on notera $d_{\rho}$ le rang (si n\'ecessaire) et de 2 entiers $a,b$ qui donnent la dimension de 2 repr\'esentations irr\'eductibles de $SL(2,{\mathbb C})$, les repr\'esentations irr\'eductibles de $SL(2,{\mathbb C})$ \'etant uniquement  caract\'eris\'ees, \`a isomorphisme pr\`es, par leur dimension. On dit que la repr\'esentation irr\'eductible $(\rho,a,b)$ a la condition de parit\'e si elle se factorise par un groupe de m\^eme type que $^LG$; pour cela, il faut d'abord que $\rho$ soit autoduale, c'est-\`a-dire isomorphe \`a $\rho^*$, ce qui entra\^{\i}ne que $\rho$ se factorise par un groupe orthogonal (on pose alors $\eta_{\rho}=+$) ou par un groupe symplectique (on pose alors $\eta_{\rho}=-$). Ensuite, il faut que $\eta_{\rho}(-1)^{a+b}=+1$ si $^LG$ est orthogonal et $+1$ si $^LG$ est symplectique.

 A un tel triplet $(\rho,a,b)$ on associe le quadruplet $(\rho,A,B,\zeta)$ o\`u $A=(a+b)/2-1$, $B=\vert (a-b)/2\vert$ et $\zeta$ est le signe de $a-b$ quand ce nombre est non nul; quand $B$ est nul, on prend $\zeta=+$; ce choix influe sur la d\'efinition, il ne change pas le paquet dans son ensemble, ce paquet \'etant d\'efini par une \'egalit\'e de transfert qui ''voit'' $\zeta$ uniquement sous la forme d'un signe; ceci est expliqu\'e dans \cite{selecta} par. 4 et 5 et en particulier 5.6. Ce choix change donc uniquement la param\'etrisation \`a l'int\'er\'erieur des paquets et seulement dans certains cas mais qui sont les cas consid\'er\'es ici. On garde la terminologie, $(\rho,A,B,\zeta)$ a la condition de parit\'e si le triplet $(\rho,a,b)$ qui lui correspond l'a; ici cela se traduit par le fait que $A,B$ sont entiers ou demi-entiers non entiers.

On a une classification de ces morphismes $\psi$, vus \`a conjugaison pr\`es,  en donnant chaque composante isotypique avec sa multiplicit\'e. On pose $Jord(\psi):=\{(\rho,A,B,\zeta)\}$ o\`u l'on fait intervenir chaque $(\rho,A,B,\zeta)$ avec la multiplicit\'e de la repr\'esentation associ\'ee dans la d\'ecomposition de $\psi$. On remarque que si $(\rho,A,B,\zeta)$ classifie une des composantes isotypiques de $\psi$ mais n'a pas la condition de parit\'e, deux cas sont possibles. Le premier cas est celui o\`u $\rho\simeq \rho^*$; dans ce cas, la repr\'esentation irr\'eductible se factorise par un groupe classique de type oppos\'e \`a $^LG$ et la multiplicit\'e avec laquelle cette repr\'esentation appara\^{\i}t est n\'ecessairement paire. L'autre cas est celui o\`u $\rho\not\simeq \rho^*$ et dans ce cas la repr\'esentation irr\'eductible correspondant \`a $(\rho^*,A,B,\zeta)$ intervient avec la m\^eme multiplicit\'e que $(\rho,A,B,\zeta)$. La repr\'esentation associ\'ee \`a $\psi$ se d\'ecompose donc en produit de repr\'esentations index\'ees par les $(\rho,A,B,\zeta)\cup (\rho^*,A,B,\zeta)$ o\`u la repr\'esentation correspondante est la composante isotypique correspondant \`a $(\rho,A,B,\zeta)$ si $\rho\simeq \rho^*$ et la somme de cette composante isotypique avec celle correspondant \`a $(\rho^*,A,B,\zeta)$ sinon. Chacune de ces repr\'esentations est \`a valeurs dans un groupe de m\^eme type que $^LG$ et le centralisateur de $\psi$ dans $^LG$ est un produit des centralisateurs de chacune de ces repr\'esentations. Pour calculer ce centralisateur, on distingue donc:

soit $(\rho,A,B,\zeta)$ ayant la condition de parit\'e, alors le centralisateur correspondant est un groupe orthogonal de la forme $O(m,{\mathbb C})$ o\`u $m$ est la longueur de la repr\'esentation c'est-\`a-dire la multiplicit\'e de $(\rho,A,B,\zeta)$ dans $Jord(\psi)$;

soit $(\rho,A,B,\zeta)$ avec $\rho\simeq \rho^*$ mais n'ayant pas la condition de parit\'e; alors le centralisateur correspondant est un groupe symplectique $Sp(2m,{\mathbb C})$ o\`u $2m $ est la multiplicit\'e de $(\rho,A,B,\zeta)$ dans $Jord(\psi)$;

soit $(\rho,A,B,\zeta)$ avec $\rho\not\simeq \rho^*$, alors le centralisateur correspondant est un groupe lin\'eaire $GL(m,{\mathbb C})$ o\`u $m$ est la multiplicit\'e de $(\rho,A,B,\zeta)$ dans $Jord(\psi)$.
 
 \
 
 Le caract\`ere $\epsilon$, dans tout ce papier, est un morphisme du centralisateur de $\psi$ dans le $L$-groupe \`a valeurs dans $\{\pm 1\}$. Avec la description ci-dessus du centralisateur, on voit que $\epsilon$ est triviale sur les facteurs correspondant \`a $(\rho,A,B,\zeta)$ n'ayant pas la condition de parit\'e. Et sur les facteurs correspondant \`a $(\rho,A,B,\zeta)$ ayant la condition de parit\'e, $\epsilon$ a 2 possibilit\'es, soit \^etre trivial soit \^etre le caract\`ere signe; il faut tenir compte du fait que $^LG$ peut avoir un centre; la restriction de $\epsilon$ au centre de $^LG$ est d\'etermin\'e par le groupe. Si $G$ est un groupe symplectique cette restriction doit \^etre triviale, au sens strict, $^LG$ est alors le groupe sp\'ecial ortogonal de dimension impair et c'est donc bien ce qu'il faut. Si $G$ est un groupe orthogonal, on fixe la forme orthogonale, elle est a donc un invariant de Hasse qui vaut $\pm 1$ et on demande que la restriction de $\epsilon$ au centre de $^LG$ soit trivial exactement si l'invariant de Hasse vaut $+1$; cela est fix\'e au d\'ebut et n'intervient quasiment pas; on pose $ \eta_{G}=+1$ si $G$ est un groupe symplectique ou si c'est un groupe orthogonal pour une forme orthogonale fix\'ee d'invariant de Hasse $+1$ et $-1$ sinon; il est clair que dans le cas des groupes orthogonaux sur un espace de dimension paire, on a fait un choix qui n'est pas dict\'e par le groupe, ou encore que les param\'etrisations donn\'ees d\'ependent de la forme orthogonale et pas seulement du groupe; c'est assez naturel. On identifiera donc $\epsilon$ \`a une application de $Jord(\psi)$ dans $\{ \pm 1\}$ qui vaut $+1$ sur tout $(\rho,A,B,\zeta)$ pour lequel $\epsilon$ restreint au facteur du centralisateur de $\psi$ correspondant \`a $(\rho,A,B,\zeta)$ est trivial et $-1$ sinon et telle que 
 $$
 \prod_{(\rho,A,B,\zeta)\in Jord(\psi)}\epsilon(\rho,A,B,\zeta)=\eta_{G}
 $$
 o\`u dans le produit $Jord(\psi)$ est vu avec multiplicit\'e.
 
\

Pour simplifier, on dira que $\psi$ v\'erifie la condition de parit\'e si tout \'el\'ement $(\rho,A,B,\zeta)$ de $Jord(\psi)$ v\'erifie la condition de parit\'e. 
\subsection{Notations\label{notations}}
Soit $(\rho,A,B,\zeta)$ comme ci-dessus, on d\'efinit la repr\'esentation $S(\rho, A,B,\zeta)$, de la fa\c{c}on suivante.  On consid\`ere le tableau:
$$
\begin{matrix}
&\zeta B &\cdots &-\zeta A\\
&\vdots &\vdots &\vdots\\
&\zeta A &\cdots &-\zeta B
\end{matrix}
$$
Si $\zeta=+$ les lignes sont des segments d\'ecroissants et les colonnes des segments croissants tandis que si $\zeta=-$ c'est l'inverse. Dans un groupe lin\'eaire convenable \`a chaque ligne correspond une induite par exemple pour la $\ell$-i\`eme ligne c'est l'induite:
$$
\rho\vert\,\vert^{\zeta (B+\ell-1)}\times \cdots \times \rho\vert\,\vert^{-\zeta (A-\ell+1)}
$$
qui a un unique sous-module irr\'eductible not\'e $<\rho\vert\,\vert^{\zeta (B+\ell-1)}, \cdots, \rho\vert\,\vert^{-\zeta (A-\ell+1)}>$. On consid\`ere alors l'induite de tous ces sous-modules; cette induite a elle aussi un unique sous-module irr\'eductible d'apr\`es les r\'esultats maintenant classiques de Zelevinsky et c'est ce sous-module irr\'eductible que l'on note $S(\rho,A,B,\zeta)$. En d'autres termes dans la classification de Zelevinsky, $S(\rho,A,B,\zeta)$ est la repr\'esentation irr\'eductible qui correspond aux multi-segments form\'es par les lignes du tableau. C'est aussi la m\^eme repr\'esentation que nous aurions obtenue en consid\'erant les multi-segments form\'es par les colonnes du tableau.

\

Donnons tout de suite une d\'efinition voisine. Ici on fixe $T$ un entier et on d\'efinit $S(\zeta,A,B,T)$ une repr\'esentation irr\'eductible d'un groupe lin\'eaire convenable. On a fait dispara\^{\i}tre $\rho$ de la notation car cette repr\'esentation sera fix\'ee par le contexte. Et on consid\`ere le tableau:
$${\cal C}(\zeta,A,B,T):=
\begin{matrix}
&\zeta (B+T) &\cdots &\zeta (A+T)\\
&\vdots &\vdots &\vdots\\
&\zeta (B+1) &\cdots &\zeta (A+1)
\end{matrix}
$$
Et $S(\zeta,A,B,T)$ est la repr\'esentation irr\'eductible bas\'ee sur $\rho$ correspondant aux multi-segments form\'es par les lignes du tableau. On peut l\`a aussi remplacer lignes par colonnes sans changer la d\'efinition.

\

On utilisera aussi constamment la notation $Jac_{x}\pi$; elle signifie ceci. Soit $\pi$ une repr\'esentation tr\`es g\'en\'eralement irr\'eductible mais pas toujours du groupe $G$ et soit $\rho$ une repr\'esentation irr\'eductible fix\'ee d'un groupe lin\'eaire $GL(d_{\rho},F)$, ce qui d\'efinit $d_{\rho}$. On suppose que l'indice de Witt de la forme dont $G$ est le groupe d'automorphisme est sup\'erieur \`a $d_{\rho}$; ainsi $G$ a un sous-groupe parabolique dont le Levi est isomorphe \`a $GL(d_{\rho},F)\times G'$ o\`u $G'$ est un groupe de m\^eme type que $G$ mais de rang plus petit. On consid\`ere le module de Jacquet de $\pi$ le long du radical unipotent de ce parabolique dans le groupe de Grothendieck associ\'e aux repr\'esentations lisses de longueur finie du groupe $GL(d_{\rho},F)\times G'$ et on d\'ecompose en fonction du support cuspidal pour l'action de $GL(d_{\rho},F)$. On ne garde (c'est-\`a-dire on projette) que les constituants irr\'eductibles sur lesquels $GL(d_{\rho},F)$ agit par $\rho\vert\,\vert^x$, ils sont donc de la forme $\rho\vert\,\vert^x\otimes \sigma'$ et $Jac_{x}\pi$ est par d\'efinition l'\'el\'ement du groupe de Grothendieck associ\'e aux repr\'esentations lisses de longueur finie de $G'$, qui est la somme (avec multiplicit\'e) de ces repr\'esentations $\sigma'$.
On pose $Jac_{x}\pi=0$ si l'indice de Witt de la forme est strictement inf\'erieur \`a $d_{\rho}$ par une convention \'evidente.

\

On peut donner la m\^eme d\'efinition pour $G$ remplacer par un groupe lin\'eaire g\'en\'eral. Pour les groupes lin\'eaires, il n'y a pas de raison de privil\'egier la gauche \`a la droite et on peut donc aussi d\'efinir $Jac^d_{x}\pi$ quand on utilise le L\'evi $G'\times GL(d_{\rho},F)$; on peut ainsi d\'efinir pour $x,y$ fix\'e $Jac_{x}Jac^d_{y}\pi$ en composant les 2 applications et on a clairement $Jac_{x}Jac^d_{y}=Jac^d_{y}Jac_{x}$.  On utilisera peu les d\'efinitions pour les groupes lin\'eaires et uniquement dans le cas ci-dessous; ici on suppose (comme en fait dans tout le travail que $\rho\simeq \rho^*$).
Soit $\pi$ et $\delta$ une repr\'esentation d'un groupe lin\'eaire. On consid\`ere l'induite $\delta\times \pi$ comme repr\'esentation d'un groupe de m\^eme type que $G$, alors:
$$
Jac_{x}\bigl(\delta\times \pi\bigr)= Jac_{x}\delta \times \pi \oplus Jac^d_{-x}\delta \times \pi \oplus \delta\times Jac_{x}\pi\eqno(1)
$$
qui r\'esulte des formules de Bernstein-Zelevinsky.

\

On est amen\'e \`a g\'en\'eraliser ces d\'efinitions \`a un ensemble ${\cal E}:=\{x,\cdots, y\}$ totalement ordonn\'e o\`u on pose par d\'efinition
$$
Jac_{x,\cdots, y}=Jac_{y}\circ \cdots \circ Jac_{x}
$$
La formule (1) ci-dessus se g\'en\'eralise \`a un ensemble totalement ordonn\'e ${\cal E}$ de fa\c{c}on un peu plus compliqu\'ee: on suppose que $\rho\simeq \rho^*$ pour simplifier les notations. Alors $Jac_{x\in {\cal E}}\bigl(\delta\times \pi\bigr)$ est la somme sur tous les d\'ecoupages de ${\cal E}$ en trois sous-ensembles (dont certains peuvent \^etre vides) ${\cal E}=\cup_{i\in [1,3]}{\cal E}_{i}$ dont chacun est totalement ordonn\'e par restriction de l'ordre sur ${\cal E}$; et le terme correspondant \`a un tel d\'ecoupage est:
$$
Jac_{x\in {\cal E}_{1}}Jac^d_{x\in -^t{\cal E}_{2}}(\delta)\times Jac_{x\in {\cal E}_{3}}\sigma,
$$
o\`u $-^t{\cal E}_{2}$ est l'ensemble des \'el\'ements de $-{\cal E}_{2}$ mais avec l'ordre renvers\'e. C'est une application directe de (1).

\subsection{D\'efinition g\'en\'erale\label{definitiongenerale}}
 On met ici un ordre total sur $Jord(\psi)$.
Pour tout $\rho$, on ordonne totalement les quadruplets $(\rho,A,B,\zeta)$ ayant la condition de parit\'e en consid\'erant que $(\rho,A,B,\zeta)>(\rho,A',B',\zeta')$ au moins si

 soit $B>B'$, soit $B=B'$ et $A>A'$, soit $B=B',A=A'$ et $\zeta=+$. 
 
 M\^eme sans tenir compte du fait qu'il peut y avoir plusieurs $\rho$, (ce qui n'a aucune importance), si $Jord(\psi)$ a de la multiplicit\'e cela ne suffit pas pour d\'efinir un ordre total, il faut prolonger arbitrairement. Soit $\psi$, $\tilde{\psi}$ des morphismes pour des groupes de rang \'eventuellement diff\'erent; on ordonne totalement $Jord(\psi)$ et $Jord(\tilde{\psi})$. On dit que $\tilde{\psi}$ domine $\psi$ s'il existe une bijection, $b$, de $Jord(\tilde{\psi})$ sur $Jord(\psi)$ respectant l'ordre et v\'erifiant pour tout $(\rho,A,B,\zeta)\in Jord(\psi)$, il existe $T_{(\rho,A,B,\zeta)}$ tel que $b^{-1}(\rho,A,B,\zeta)=(\rho,A+T_{\rho,A,B,\zeta},B+T_{\rho,A,B,\zeta},\zeta)$. Il n'est pas tr\`es agr\'eable de fixer l'ordre sur $Jord(\tilde{\psi})$, il vaut mieux fixer $b$ la bijection avec la condition pr\'ec\'edente et  tel qu'il existe un ordre sur $Jord(\tilde{\psi})$ respect\'e par $b$; l'ordre est clairement uniquement d\'etermin\'e par $b$. On finira par montrer que les d\'efinitions ne d\'ependent pas du choix de $b$.
 
 On peut signaler que c'est dans le choix de ''dominant'' que la convention $\zeta=+$ si $B=0$ intervient car $b^{-1}(\rho,A,B,\zeta)$ aura aussi $\zeta=+$ met $B+T_{\rho,A,B,\zeta}\neq 0$.
  
 \
 
 On remarque que si $\psi$ v\'erifie la condition de parit\'e, il existe $\tilde{\psi}$ dominant $\psi$ et tel que $\tilde{\psi}$ est de restriction discr\`ete \`a la diagonale, c'est-\`a-dire que la restriction de $\tilde{\psi}$ \`a $W_{F}\times SL(2,{\mathbb C})$, o\`u $SL(2,{\mathbb C})$ est plong\'e diagonalement dans $SL(2,{\mathbb C})\times SL(2,{\mathbb C})$ est sans multiplicit\'e en tant que repr\'esentation de ce groupe.
 
\subsection{Hypoth\`ese\label{hypothese}}
On suppose que $\psi$ v\'erifie la condition de parit\'e.
On fixe $\tilde{\psi}$ dominant $\psi$ et pour pouvoir construire les repr\'esentations associ\'ees \`a $\psi$ on a besoin de savoir construire celles associ\'ees \`a $\tilde{\psi}$; l'hypoth\`ese que nous ferons dans tout le papier sera donc que l'on sait construire ces repr\'esentations, c'est-\`a-dire que l'on consid\`ere le morphisme $\tilde{\psi}\circ \Delta$ de restriction de $\tilde{\psi}$ \`a $W_{F}$ fois la diagonale de $SL(2,{\mathbb C})$. C'est un morphisme de $W_{F}\times SL(2,{\mathbb C})$ dans un groupe dual analogue \`a celui de $G$ mais de rang \'eventuellement plus grand. En composant avec la repr\'esentation naturelle de ce $L$-groupe, on voit $\tilde{\psi}\circ \Delta$ comme une repr\'esentation de $W_{F}\times SL(2,{\mathbb C})$ n\'ecessairement semi-simple; toute sous-repr\'esentation irr\'eductible est associ\'e \`a une morphisme $\psi'$ de $W_{F}\times SL(2,{\mathbb C})$ n\'ecessairement dans le groupe dual d'un groupe de m\^eme type que $G$; on suppose alors que l'on sait construire toutes les s\'eries discr\`etes associ\'ees \`a $\psi'$ pour tout choix de $\psi'$, c'est-\`a-dire qu'\`a tout caract\`ere, $\epsilon'$,  du centralisateur de $\psi'$, on sait associer une s\'erie discr\`ete $\pi(\psi',\epsilon')$ irr\'eductible avec les propri\'et\'es suivantes:
on d\'ecompose $\psi'$ en somme de repr\'esentations irr\'eductibles.  On rappelle qu'une composante irr\'eductible est  de la forme $\rho\otimes \sigma_{a}$ o\`u $\rho$ est une repr\'esentation irr\'eductible de $W_{F}$ autoduale et $\sigma_{a}$ la repr\'esentation irr\'eductible de $SL(2,{\mathbb C})$ de dimension a. Alors on note encore $\rho$ la repr\'esentation cuspidale du groupe lin\'eaire correspondant \`a $\rho$ par la correspondance de Langlands locale et on note $St(\rho,a)$ la repr\'esentation de Steinberg associ\'ee \`a $\rho$ et l'entier $a$. Pour $\rho$ comme ci-dessus, on pose $\eta_{\rho}=+$ si le morphisme associ\'e se factorise par un groupe orthogonal et $\eta_{\rho}=-1$ s'il se factorise par un groupe symplectique. On a d\'ej\`a d\'efini $\eta_{G}$ en \ref{blocsdejordan}.

\

Les hypoth\`ese que nous faisons sur l'ensemble des repr\'esentations $\pi(\psi',\epsilon')$ dont nous venons de supposer l'existence comme s\'eries discr\`etes irr\'eductibles est que

\

\sl pour tout $\epsilon'$ et pour tout $\rho\simeq \rho^*$ une repr\'esentation cuspidale irr\'eductible autoduale d'un groupe lin\'eaire et tout entier $a$, l'induite
$St(\rho,a)\times \pi(\psi',\epsilon')$ est irr\'eductible si est seulement si soit $\eta_{\rho}(-1)^a\neq \eta_{G}$ soit $\rho\otimes \sigma_{a}$ est une composante irr\'eductible de la repr\'esentation $\psi'$.\rm

\

C'est comme cela que l'on a introduit les blocs de Jordan en \cite{europe} et on a d\'emontr\'e cette hypoth\`ese pour les groupes unitaires (en suivant les id\'ees d'Arthur) en \cite{discretunitaire}; pour les groupes orthogonaux on renvoie \`a \cite{discretorthogonaux} et plus g\'en\'eralement
on a discut\'e dans l'introduction de \cite{paquetsdiscrets} l'\'etat de ces hypoth\`eses. 

On a choisi dans ce travail de travailler avec les groupes orthogonaux et symplectiques pour avoir explicitement un $L$-groupe tr\`es commode. Les m\'ethodes se g\'en\'eralisent au cas des groupes unitaires sans changement sauf dans les notations (cf \cite{discretunitaire}).
 \subsection{Construction\label{construction}}
Ici on continue avec $\psi$ ayant la propri\'et\'e de parit\'e (le cas g\'en\'eral s'obtient  \`a la fin cf \ref{inductionfin}). On fixe $\tilde{\psi}$ de restriction discr\`ete \`a la diagonale, dominant $\psi$, \`a l'aide d'une bijection $b$. Pour tout caract\`ere, $\epsilon$, du centralisateur de $\psi$ on d\'eduit un caract\`ere du centralisateur de $\tilde{\psi}$, not\'e $\tilde{\epsilon}$ tout simplement en posant $$\tilde{\epsilon}(\tilde{\rho},\tilde{A},\tilde{B},\tilde{\zeta})=\epsilon(b(\tilde{\rho},\tilde{A},\tilde{B},\tilde{\zeta})).
$$

%%%%%%%
 On rappelle qu'avec les notations ci-dessus on a d\'efini ${\cal C}(\zeta,A,B,T_{\rho,A,B,\zeta})$ et on ordonne totalement cet ensemble en prenant d'abord les lignes puis les colonnes. Et si $\psi$ est discret, on a pos\'e:
 $$
 \pi(\psi,\epsilon):=\circ_{(\rho,A,B,\zeta)\in Jord(\psi)}Jac_{{\cal C}(\zeta,A,B,T_{\rho,A,B,\zeta})}\pi(\tilde{\psi},\epsilon),$$
 o\`u $\tilde{\psi}$ est un morphisme de restriction discr\`ete \`a la diagonale dominant $\psi$ et o\`u $(\rho,A,B,\zeta)$ sont pris dans l'ordre croissant (on commence d'abord par calculer le $Jac$ pour le plus petit) et o\`u l'on a transport\'e $\epsilon$ par la bijection fix\'ee entre $Jord(\psi)$ et $Jord(\tilde{\psi})$. 
  
 On a montr\'e qu'avec cette d\'efinition $\pi(\psi,\epsilon)$ est bien d\'efini, c'est-\`a-dire ne d\'epend pas du choix de $\tilde{\psi}$.  Dans cette d\'efinition, il n'y a aucune raison de se limiter aux caract\`eres du centralisateur de $\psi$; on peut directement partir de $\tilde{\epsilon}$ un caract\`ere du centralisateur de $\tilde{\psi}$, on d\'efinit alors 
 $$
 \pi_{b}(\psi,\tilde{\epsilon})= \circ_{(\rho,A,B,\zeta)\in Jord(\psi)}Jac_{{\cal C}(\zeta,A,B,T_{\rho,A,B,\zeta})}\pi(\tilde{\psi},\tilde{\epsilon}),
 $$
 cette d\'efinition n'a de sens qu'une fois $b$ fix\'e. L'ind\'ependance d\'emontr\'ee en \cite{paquetsdiscrets} se formule alors ainsi: soit  $\tilde{\psi'}$, de restriction discr\`ete \`a la diagonale, dominant aussi $\psi$ avec une bijection $b'$, on transporte $\tilde{\epsilon}$ en un caract\`ere, $\tilde{\epsilon}' $ du centralisateur de $\tilde{\psi}'$ \`a l'aide de $b^{-1}b'$ et 
 $$
 \pi_{b}(\psi,\tilde{\epsilon})=\pi_{b'}(\psi, \tilde{\epsilon}').
 $$
 D'ici peu, nous allons montrer que cette g\'en\'eralisation est fictive puisque les repr\'esentations apparemment suppl\'ementaires sont en fait nulles.

 \

 On rappelle  que l'on a une description pr\'ecise de $\pi(\tilde{\psi},\tilde{\epsilon})$ puisque $\tilde{\psi}$ est de restriction discr\`ete \`a la diagonale. Avec les notations ci-dessus:

 \sl $\pi(\tilde{\psi},\tilde{\epsilon})$ est une repr\'esentation semi-simple sans multiplicit\'e dont les composantes irr\'eductibles sont en bijection avec la donn\'ee pour tout $(\rho,A,B,\zeta,T_{\rho,A,B,\zeta})$ d'un entier $\underline{\ell}(\rho,A,B,\zeta,T_{\rho,A,B,\zeta})$ positif ou nul, inf\'erieur ou \'egal \`a $(A-B+1)/2$ et d'un signe $\underline{\eta}(\rho,A,B,\zeta,T_{\rho,A,B,\zeta})$ v\'erifiant:
 $$
 \prod_{C\in [B+\underline{\ell}(\rho,A,B,\zeta,T_{\rho,A,B,\zeta}),A-\underline{\ell}(\rho,A,B,\zeta,T_{\rho,A,B,\zeta})]} \biggl(\underline{\eta}(\rho,A,B,\zeta,T_{\rho,A,B,\zeta})(-1)^C\biggr)=\tilde{\epsilon}(\rho,A,B,\zeta,T_{\rho,A,B,\zeta}),
 $$o\`u un produit vide vaut $+1$. \rm
 
 La condition sur le signe, fixe la parit\'e de $\underline{\ell}(\rho,A,B,\zeta,T_{\rho,A,B,\zeta})$ si $A-B+1$ est pair et fixe $\underline{\eta}(\rho,A,B,\zeta,T_{\rho,A,B,\zeta})$ si $A-B$ est pair.
\nl
On suppose que $\psi$ a la condition de parit\'e; on fixe  $\tilde{\psi}$ et la bijection $b$ de $Jord(\tilde{\psi})$ sur $Jord(\psi)$ et un carac\-t\`ere $\tilde{\epsilon}$ du centralisateur de $\tilde{\psi}$. On fixe des fonctions $\tilde{\underline{\ell}}$ et $\tilde{\underline{\eta}}$ de fa\c{c}on \`a fixer une composante irr\'eductible de $\pi(\tilde{\psi},\tilde{\epsilon})$; et on pose, $$\pi_{b}(\psi,\tilde{\epsilon},\tilde{\underline{\ell}},\tilde{\underline{\eta}}):=$$
$$\circ_{(\rho,A,B,\zeta)\in Jord(\psi)}Jac_{{\cal C}(\zeta,A,B,T_{\rho,A,B,\zeta})}\biggl(\pi(\tilde{\psi},\tilde{\epsilon},\tilde{\underline{\ell}},\tilde{\underline{\eta}})\biggr);$$  en anticipant les r\'esultats ult\'erieurs, on montrera que ces modules de Jacquets sont nuls sauf si $\tilde{\underline{\ell}}$ et $\tilde{\underline{\eta}}$ proviennent de fonctions sur $Jord(\psi)$, au sens le plus naturel pour $\tilde{\underline{\ell}}$ et avec une petite torsion pour $\tilde{\underline{\eta}}$; mais cela imposera que $\tilde{\epsilon}$ est l'image r\'eciproque d'un caract\`ere du centralisateur de $\psi$. C'est ce qui permettra de d\'efinir ind\'ependamment de $b$ les repr\'esentations $\pi(\psi,\epsilon,\underline{\ell},\underline{\eta})$. Pour le moment on sait que $\pi(\psi,\epsilon)$ est la somme des repr\'esentations $\pi_{b}(\psi,\epsilon,\tilde{\underline{\ell}},\tilde{\underline{\eta}})$ pour un choix de $b$ et les fonctions $\tilde{\underline{\ell}}$ et $\tilde{\underline{\eta}}$ 
d\'efinies sur $Jord(\tilde{\psi})$ variant.

  %%%%%%%%%

\section{Propri\'et\'es g\'en\'erales\label{proprietesgenerales}}

\subsection{Irr\'eductibilit\'e, \'enonc\'e\label{irreductibilite}}

\bf Proposition: \sl pour $\tilde{\underline{\ell}}$ et $\tilde{\underline{\eta}}$ fix\'es la repr\'esentation $\pi_{b}(\psi,\tilde{\epsilon},\tilde{\underline{t}},\tilde{\underline{\eta}})$ est soit nulle soit irr\'eductible. De plus l'ensemble des repr\'esentations non nulles obtenues ainsi (pour $\psi$ et $b$ fix\'e) est sans multiplicit\'e.\rm

\

Une cons\'equence imm\'ediate de la proposition est que $\pi(\psi,\epsilon)$ est semi-simple sans multiplicit\'e.
On d\'emontre cette proposition en m\^eme temps que le lemme ci-dessous dont on a besoin dans la d\'emonstration et qui a son int\'er\^et en soi.

\subsection{Une propri\'et\'e simple des modules de Jacquet \label{jacquetsimple}}

\bf Lemme: \sl soit $x\in {\mathbb R}$ et soit $m$ le cardinal de l'ensemble $\{(\rho,A,B,\zeta)\in Jord(\psi); \zeta B=x\}$ alors $Jac_{x, \cdots, x} \pi_{b}(\psi,\tilde{\epsilon})=0$ d\`es que le nombre de $x$ apparaissant est strictement sup\'erieur \`a $m$. \rm

Le lemme et la proposition sont  vraies si $\psi$ est de restriction discr\`ete \`a la diagonale (cf. \cite{paquetsdiscrets}, introduction). On reprend les notations $\tilde{\psi}, \tilde{\epsilon}$ et $b$ qui servent \`a d\'efinir $\pi(\psi,\tilde{\epsilon})$, pour simplifier la notation on remplace $\tilde{\underline{\ell}}$ et $\tilde{\underline{\eta}}$ par $\ell$ et $\eta$ \'ecrits en indice; ces fonctions n'interviennent pas dans la d\'emonstration, elles permettent juste de travailler avec des repr\'esentations qui vont \^etre irr\'eductibles et la r\'eciprocit\'e de Frobenius est plus simple. On suppose comme on en a le droit que si $(\rho,A_{i},B_{i},\zeta_{i})\in Jord(\tilde{\psi})$ pour $i=1,2$ distincts alors soit $A_{1}<< B_{2}$ soit $B_{1}>> A_{2}$.

On a donc le lemme et la proposition pour $\pi(\tilde{\psi},\tilde{\epsilon})_{\ell,\eta}$. On les d\'emontre de proche en proche pour chaque module de Jacquet interm\'ediaire qui servent \`a la d\'efinition. On peut donc simplement consid\'erer la situation suivante:  $\psi', \psi''$ sont des morphismes (non n\'ecessairement de restriction discr\`ete \`a la diagonale) qui dominent $\psi$ et tels que $\psi'$ domine $\psi''$ \`a l'aide d'une bijection not\'ee $b'$; deplus il existe $(\rho,A',B',\zeta')\in Jord(\psi')$ y intervenant avec multiplicit\'e 1 et  tel que $B'\geq 1$,   $b'$ est l'identit\'e sur tous les \'el\'ements de $Jord(\psi')$ autre que $(\rho,A',B',\zeta')$ et v\'erifie $b'(\rho,A',B',\zeta')=(\rho,A'-1,B'-1,\zeta')$ et tous les \'el\'ements de $Jord(\psi')$ strictement plus grands que $(\rho,A',B',\zeta')$ sont de la forme $(\tilde{\rho},\tilde{A},\tilde{B},\tilde{\zeta})$ avec $\tilde{B}>>B'$ et les $\tilde{B}$ sont tr\`es diff\'erents les uns des autres. On admet le lemme et la proposition pour $\pi(\psi',\tilde{\epsilon})_{\ell,\eta}$ et on les montre pour $$\pi(\psi'',\tilde{\epsilon})_{\ell,\eta}=Jac_{\zeta' B', \cdots, \zeta' A'}\pi(\psi',\tilde{\epsilon})_{\ell,\eta}.$$ 
Le r\'esultat est trivialement vrai si $\pi(\psi'',\tilde{\epsilon})_{\ell,\eta}=0$; on suppose donc que l'on n'est pas dans ce cas de nullit\'e. On fixe $x$ comme dans l'\'enonc\'e du lemme et on note $m(\psi',x)$ le plus grand entier tel que $Jac_{x, \cdots, x}\pi(\psi',\tilde{\epsilon})\neq 0$ o\`u il y a $m(\psi',x)$-copies de $x$; ce nombre d\'epend \'evidemment aussi de $\tilde{\epsilon}$ mais cela n'a pas d'importance ici, $\tilde{\epsilon}$ \'etant fix\'e et n'intervenant pas pour la borne de l'\'enonc\'e. On d\'efinit de m\^eme $m(\psi'',x)$ et on montre que:
$$
m(\psi'',x)\begin{cases}\leq m(\psi',x) \hbox{ si }x\neq \zeta'(B'-1)\\
\leq m(\psi',x)-1 \hbox{ si } x=\zeta ' B'\\
\leq m(\psi',x)+1 \hbox { si } x= \zeta' (B'-1).
\end{cases}.
$$Cela suffira pour prouver le lemme.

En suivant les d\'efinitions, on a $Jac_{x, \cdots, x}Jac_{\zeta' (B'), \cdots, \zeta' A'}(\pi(\psi',\tilde{\epsilon})\neq 0$ pour $m(\psi'',x)$ copies de $x$. Par r\'eciprocit\'e de Frobenius, il existe une inclusion de la forme
$$
\pi(\psi',\tilde{\epsilon})_{\ell,\eta}\hookrightarrow \rho\vert\,\vert^{\zeta B}\times \cdots \times \rho\vert\,\vert^{\zeta A}\times \rho\vert\,\vert^x \times \cdots \times \rho\vert\,\vert^x\times \tau
$$
o\`u $\tau$ est une repr\'esentation irr\'eductible non nulle. Comme $Jac_{y}\pi(\psi',\tilde{\epsilon})=0$ pour tout $y\in [\zeta B,\zeta A]$ ce morphisme se factorise en 
$$
\pi(\psi',\tilde{\epsilon})_{\ell,\eta}\hookrightarrow <\rho\vert\,\vert^{\zeta B}, \cdots , \rho\vert\,\vert^{\zeta A}>\times \rho\vert\,\vert^x \times \cdots \times \rho\vert\,\vert^x\times \tau \eqno(*)
$$
o\`u $<\rho\vert\,\vert^{\zeta B}, \cdots , \rho\vert\,\vert^{\zeta A}>$ est l'unique sous-module de l'induite \'evidente car cette repr\'esentation est le seul sous-quotient de cette induite \`a v\'erifier $Jac_{y}=0$ pour tout $y$ comme ci-dessus. Supposons maintenant que $x\neq \zeta (B'-1) $ et $x\neq \zeta (A'+1)$.  Sous cette hypoth\`ese $$ <\rho\vert\,\vert^{\zeta B}, \cdots , \rho\vert\,\vert^{\zeta A}>\times \rho\vert\,\vert^x$$ est une induite irr\'eductible d'apr\`es les r\'esultats de Zelevinsky et on peut donc \'echanger les facteurs. On trouve donc que $m(\psi',x)\geq m(\psi'',x)$ si $x\neq \zeta B'$ et $m(\psi',\zeta' B') \geq m(\psi'',\zeta' B')+1$. Cela donne le lemme pour tout $x\neq \zeta' (B'-1)$ et $x\neq \zeta'(A'+1)$. Supposons que $x=\zeta (B'-1)$;  les r\'esultats de Zelevinsky montre que l'induite $$<\rho\vert\,\vert^{\zeta B}, \cdots , \rho\vert\,\vert^{\zeta A}>\times \rho\vert\,\vert^x$$ est de longueur $2$. Soit $\sigma$ l'un de ces sous-quotients; Zelevinsky montre aussi que $\sigma\times \vert\,\vert^x$ est irr\'eductible. Ainsi, pour un bon choix de $\sigma$, (*) se factorise en un morphisme non nul
$$
\pi(\psi',\tilde{\epsilon})_{\ell,\eta}\hookrightarrow \rho\vert\,\vert^x \times \cdots \rho\vert\,\vert^x \times \sigma \times \tau,
$$
o\`u il y a $m(\psi'',x)-1$ copies de $x$. Le lemme est donc aussi d\'emontr\'e pour $x=\zeta'(B'-1)$. Il reste $x=\zeta' (A'+1)$. Ici on veut d\'emontrer que $m(\psi'',\zeta'(A'+1))=m(\psi',\zeta'(A'+1))=0$. Il faut donc d\'emontrer que $Jac_{\zeta' (B'), \cdots, \zeta' A',\zeta' (A'+1) }\pi(\psi',\tilde{\epsilon})=0$. On suppose par l'absurde que ce module de Jacquet est non nul d'o\`u une inclusion pour $\tau'$ une repr\'esentation convenable (r\'eciprocit\'e de Frobenius)
$$
\pi(\psi',\tilde{\epsilon})_{\ell,\eta}\hookrightarrow \rho\vert\,\vert^{\zeta' B'} \times \cdots \times  \rho\vert\,\vert^{\zeta' (A'+1)}\times \tau';$$ 
on sait d\'ej\`a que $Jac_{y}\pi(\psi',\tilde{\epsilon})=0$ pour tout $y\in ]\zeta' B',\zeta'(A'+1)]$ et le morphisme ci-dessus se factorise donc par le sous-module, $<\rho\vert\,\vert^{\zeta' B'}, \cdots, \rho\vert\,\vert^{\zeta' (A'+1)}>\times \tau'$ de l'induite $\rho\vert\,\vert^{\zeta' B'} \times \cdots \rho\vert\,\vert^{\zeta' (A'+1)}\times \tau'$.  On se rappelle que l'on commence par un morphisme $\tilde{\psi}$ de restriction discr\`ete \`a la diagonale et $\tilde{\psi}$ domine $\psi'$ avec un bijection $b'$; pour passer de $\tilde{\psi}$ \`a $\psi'$ on a fait diminuer tous les blocs de Jordan plus petits que $(b')^{-1}(\rho,A',B',\zeta')$; on distingue 2 cas.  Le premier est celui o\`u $(b')^{-1}(\rho,A',B',\zeta')=(\rho,A',B',\zeta')$. Comme les blocs de Jordan de $\tilde{\psi}$ sont tr\`es diff\'erents les uns des autres par hypoth\`ese de d\'epart, il est \`a peu pr\`es imm\'ediat que $Jac_{\zeta' B', \cdots, \zeta' (A'+1) }\pi(\psi',\tilde{\epsilon})_{\ell,\eta}\neq 0$ n\'ecessite aussi que $Jac_{\zeta' B', \cdots, \zeta' (A'+1)}\pi(\tilde{\psi},\tilde{\epsilon})_{\ell,\eta}\neq 0$. Et cela on sait que ce n'est pas possible car $\tilde{\psi}$ est de restriction discr\`ete \`a la diagonale et il n'existe pas $(\rho,A''',B''',\zeta''')$ avec $A'+1 \in [B''',A''']$ (cf. \cite{paquetsdiscrets}, 3.4). Le 2e cas est donc celui o\`u il existe $\psi'''$ tel que $Jord(\psi''')$ se d\'eduit de $Jord(\psi')$ en rempla\c{c}ant $(\rho,A',B',\zeta')$ par $(\rho,A'+1,B'+1,\zeta')$ et par d\'efinition on a
$$
\pi(\psi',\tilde{\epsilon})_{\ell,\eta}=Jac_{\zeta' (B'+1), \cdots, \zeta' (A'+1)}\pi(\psi''',\tilde{\epsilon})_{\ell,\eta}.
$$
On sait par r\'ecurrence que $\pi(\psi',\tilde{\epsilon})_{\ell,\eta}$ et $\pi(\psi''',\tilde{\epsilon})_{\ell,\eta}$ sont irr\'eductibles d'o\`u une inclusion par r\'eciprocit\'e de Frobenius
$$
\pi(\psi''',\tilde{\epsilon})_{\ell,\eta}\hookrightarrow \rho\vert\,\vert^{\zeta' (B'+1)}\times \cdots \times \rho\vert\,\vert^{\zeta' (A'+1)}\times \pi(\psi',\tilde{\epsilon})_{\ell,\eta}.
$$En utilisant la propri\'et\'e du module de Jacquet de $\pi(\psi',\tilde{\epsilon})_{\ell,\eta}$, on obtient un prolongement de cette inclusion en une inclusion
$$
\pi(\psi''',\tilde{\epsilon})_{\ell,\eta}\hookrightarrow \rho\vert\,\vert^{\zeta' (B'+1)}\times\cdots \times  \rho\vert\,\vert^{\zeta' (A'+1)}\times <\rho\vert\,\vert^{\zeta' B'}, \cdots, \rho\vert\,\vert^{\zeta' (A'+1)}>\times \tau'.
$$
Pour tout $x\in [\zeta'(B'+1),\zeta'(A'+1)]$, l'induite $\rho\vert\,\vert^x\times <\rho\vert\,\vert^{\zeta' B'}, \cdots, \rho\vert\,\vert^{\zeta' (A'+1)}>$ est irr\'eductible et on peut donc mettre $<\rho\vert\,\vert^{\zeta' B'}, \cdots, \rho\vert\,\vert^{\zeta' (A'+1)}>$ en premi\`ere place dans l'induction. D'o\`u $$Jac_{\zeta' B'}\pi(\psi''',\tilde{\epsilon})_{\ell,\eta}\neq 0;$$ or pout tout $(\rho,A''',B''',\zeta''')\in Jord(\psi''')$, $\zeta'B'\neq \zeta'''B'''$ et ceci est donc impossible. On a donc prouv\'e le lemme pour $\pi(\psi',\tilde{\epsilon})_{\ell,\eta}$.

\nl
Montrons maintenant la proposition. On suppose que $\pi(\psi',\tilde{\epsilon})_{\ell,\eta}$ est non nul et on va montrer que cette repr\'esentation est irr\'eductible et qu'elle d\'etermine uniquement la repr\'esentation $\pi(\tilde{\psi},\tilde{\epsilon})_{\ell,\eta}$, ceci prouvera la proposition de proche en proche.

Par r\'eciprocit\'e de Frobenius pour un quotient irr\'eductible $X$ bien choisi de $\pi(\psi',\tilde{\epsilon})_{\ell,\eta}$, il existe une inclusion
$$
\pi(\psi'',\tilde{\epsilon})_{\ell,\eta} \hookrightarrow \rho\vert\,\vert^{\zeta' B'}\times \cdots \times \rho\vert\,\vert^{\zeta' A'}\times X.
$$
On calcule $Jac_{\zeta' B', \cdots , \zeta' A'}$ du membre de droite et on trouve $X$ car $B'>0$ et pour tout $x\in [\zeta' B', \zeta' A']$, $Jac_{x}X=0$ puisque $Jac_{x}\pi(\psi',\tilde{\epsilon})_{\ell,\eta}=0$ d'apr\`es le lemme que l'on vient de d\'emontrer. Cela prouve que $X=\pi(\psi',\tilde{\epsilon})_{\ell,\eta}$ et  l'irr\'eductibilit\'e cherch\'ee. Cette preuve  montre aussi que si $\pi(\psi',\tilde{\epsilon})_{\ell,\eta}$ est non nul, il d\'etermine uniquement $\pi(\psi'',\tilde{\epsilon})_{\ell,\eta}$ puisque c'est l'unique sous-module irr\'eductible de l'induite \'ecrite. Et en remontant, on obtient le corollaire ci-dessous.
\subsection{Une cons\'equence technique\label{consequence}}
\bf Corollaire: \sl on suppose ici que $\pi(\psi,\tilde{\epsilon},\underline{\ell},\underline{\eta})\neq 0$ et qu'il existe $(\rho,A,B,\zeta)\in Jord(\psi)$ tel que d'une part
$ Jac_{\zeta B, \cdots, \zeta A}\pi(\psi,\tilde{\epsilon},\underline{\ell},\underline{\eta})\neq 0
$ et d'autre part
$\forall x\in ]B,A], \not\exists (\rho,A',B',\zeta)\in Jord(\psi), x=\zeta B'.$

Alors il existe une repr\'esentation irr\'eductible $\tau$ tel que $\pi(\psi,\tilde{\epsilon},\underline{\ell},\underline{\eta})$ soit un sous-module irr\'eductible de l'induite
$
<\rho\vert\,\vert^{\zeta B},\cdots, \rho\vert\,\vert^{\zeta A}>\times \tau.$\rm

\

La r\'eciprocit\'e de Frobenius donne l'existence de $\tau$ avec une inclusion
$$
\pi(\psi,\tilde{\epsilon},\underline{\ell},\underline{\eta})
\hookrightarrow \rho\vert\,\vert^{\zeta B}\times \cdots \times \rho\vert\,\vert^{\zeta A}\times \tau.
$$
Dans cette inclusion on peut remplacer l'induite $\rho\vert\,\vert^{\zeta B}\times \cdots \times \rho\vert\,\vert^{\zeta A}$ par un sous-quotient irr\'eductible $\tau'$. Mais gr\^ace \`a \ref{jacquetsimple} et \`a l'hypoth\`ese sur $Jord(\psi)$ on sait que $Jac_{x}\tau'=0$ pour tout $x\in ]\zeta B,\zeta A]$; n\'ec\'essairement $\tau'$ est la repr\'esentation $<\rho\vert\,\vert^{\zeta B},\cdots, \rho\vert\,\vert^{\zeta A}>$.

\section{D\'ecomposition des repr\'esentations d\'efinies\label{definition}}
\subsection{Un crit\`ere de nullit\'e\label{nullite}}
On veut arriver \`a une description des composantes irr\'eductibles (ou malheureusement 0) de  $\pi(\psi,\epsilon)$ en utilisant des fonctions $\underline{\ell}$ et $\underline{\eta}$ d\'efinies directement sur $Jord(\psi)$ et ne d\'ependant pas du choix d'une bijection $b$. On a en vue le r\'esultat suivant; fixons $\tilde{\psi}$ dominant $\psi$ \`a l'aide d'une bijection $b$; soit  $(\rho,A,B,\zeta)\in Jord(\psi)$ y intervenant avec multiplicit\'e au moins 2. Pour $i=1,2$, soit $T_{i}\in {\mathbb N}$ tel que $(\rho,A+T_{i},B+T_{i},\zeta)$ soit des \'el\'ements de $Jord(\tilde{\psi})$ d'image par $b$ une des copies de $(\rho,A,B,\zeta)\in Jord(\psi)$ et cons\'ecutifs pour l'ordre dans $Jord(\tilde{\psi})$; on fixe des fonctions $\tilde{\underline{\ell}}$ et $\tilde{\underline{\eta}}$ sur $Jord(\tilde{\psi})$ de fa\c{c}on \`a pouvoir d\'efinir $\pi_{b}(\psi,\epsilon,\tilde{\underline{\ell}},\tilde{\underline{\eta}})$; ces fonctions d\'efinissent un caract\`ere $\tilde{\epsilon}$ du centralisateur de $\tilde{\psi}$ (si le lecteur pr\'ef\`ere, il fixe $\tilde{\epsilon}$ et impose des conditions aux fonctions). On pose, pour $i=1,2$, $\ell_{i}:=\underline{\ell}(\rho,A+T_{i},B+T_{i},\zeta)$ et $\eta_{i}:=\underline{\eta}(\rho,A+T_{i},B+T_{i},\zeta)$.

\nl
\bf Lemme: \sl la repr\'esentation $\pi_{b}(\psi,\tilde{\epsilon},\tilde{\underline{\ell}},\tilde{\underline{\eta}})$
est nulle sauf \'eventuellement si $\ell_{1}=\ell_{2}$ et $\eta_{1}\eta_{2}=(-1)^{A-B}$.\rm
\nl
Avant de d\'emontrer le lemme remarquons ses cons\'equences sur la description des repr\'esentations $\pi(\psi,\epsilon)$.
\subsection{Param\'etrisation des composantes irr\'eductibles\label{decomposition}}
La premi\`ere cons\'equence du lemme  est que si  la repr\'esentation $\pi(\psi,\tilde{\epsilon}, \tilde{\underline{\ell}}, \tilde{\underline{\eta}})$ est non nulle
$$
\tilde{\epsilon}(\rho,A+T_{1},B+T_{1},\zeta)=\eta_{1}^{A-B+1}(-1)^{[(A-B+1)/2]+\ell_{1}}$$
$$=\eta_{2}^{A-B+1} (-1)^{(A-B)(A-B+1)}(-1)^{[(A-B+1)/2]+\ell_{2}}=\tilde{\epsilon}(\rho,A+T_{2},B+T_{2},\zeta).
$$
Donc en particulier les repr\'esentations $\pi_{b}(\psi,\tilde{\epsilon})$ sont nulles sauf \'eventuellement si $\tilde{\epsilon}$ provient d'un caract\`ere $\epsilon$ du centralisateur de $\psi$. 

 Montrons aussi comment ce lemme permet de d\'efinir  $\pi(\psi,\epsilon,\underline{\ell},\underline{\eta})$ pour des fonctions $\underline{\ell}$ et $\underline{\eta}$ d\'efinies sur $Jord(\psi)$ et v\'erifiant, pour tout $(\rho,A,B,\zeta)\in Jord(\psi)$
$$
\underline{\ell}(\rho,A,B,\zeta)\in [0, [(A-B+1)/2]; \qquad \underline{\eta}(\rho,A,B,\zeta)^{A-B+1}(-1)^{[(A-B+1)/2]+\underline{\ell}(\rho,A,B,\zeta)}=\epsilon (\rho,A,B,\zeta).
$$
On d\'efinira d'abord une fonction $\tilde{\underline{\ell}}$  sur $Jord(\tilde{\psi})$ en posant pour tout $(\rho',A',B',\zeta') \in Jord(\tilde{\psi})$, $\tilde{\underline{\ell}}(\rho',A',B',\zeta')=\underline{\ell}(b(\rho',A',B',\zeta'))$. Pour d\'efinir $\tilde{\underline{\eta}}(\rho',A',B',\zeta')$, il y a un choix \`a faire si $A-B+1$ est pair;  on peut envoyer $Jord(\psi)$ dans le m\^eme ensemble vu sans multiplicit\'e et en composant $b$ avec cette application on obtient une application non injective de $Jord(\psi')$ dans cet ensemble, not\'e $b'$. On voudrait poser $\tilde{\underline{\eta}}(\rho',A',B',\zeta')=\underline{\eta}(b'(\rho',A',B',\zeta'))$ mais la difficult\'e est que cette fonction ne v\'erifie pas l'alternance du lemme sur les fibres de $b'$ si $A'-B'$ est impair; il faut donc r\'etablir cette alternance en faisant un choix. Les fibres de $b'$ sont totalement ordonn\'ees par l'ordre mis sur $Jord(\tilde{\psi})$ et on demande l'\'egalit\'e $\tilde{\underline{\eta}}(\rho',A',B',\zeta')=\underline{\eta}(b(\rho',A',B',\zeta'))$ pour tous les \'el\'ements en position impaire dans la fibre (en commen\c{c}ant par le plus petit) et l'\'egalit\'e $\tilde{\underline{\eta}}(\rho',A',B',\zeta')=(-1)^{A'-B'}\underline{\eta}(b(\rho',A',B',\zeta'))$ pour tous les \'el\'ements en position paire. Et on pose:
$$
\pi(\psi,\epsilon,\underline{\ell},\underline{\eta})=\pi_{b}(\tilde{\psi},\epsilon,\tilde{\underline{\ell}},\tilde{\underline{\eta}}).
$$
Le lemme montre que $\pi(\psi,\epsilon)$ est la somme sans multiplicit\'e des repr\'esentations (irr\'eductibles ou nulles) $\pi(\psi,\epsilon,\underline{\ell},\underline{\eta})$ que nous venons de d\'efinir.
%%%%%%%
\subsection{Une r\'eduction pour \ref{nullite} \label{reduction}}
On garde les notations de \ref{decomposition}. On suppose que $T_{1}<T_{2}$ et on consid\`ere l'ensemble des nombres:
$${\cal E}:=\begin{matrix}
&\zeta (B+T_{2}), &\cdots, &\zeta( A+T_{2})\\
&\vdots, &\vdots, &\vdots\\
&\zeta(B+T_{1}+1),&\cdots, &\zeta (A+T_{1}+1). 
\end{matrix}
$$Cet ensemble a d\'ej\`a \'et\'e introduit avec la notation plus compliqu\'ee ${\cal C}(\zeta,A+T_{1},B+T_{1},T_{2}-T_{1})$.
On ordonne ${\cal E}$ en prenant de gauche \`a droite puis de haut en bas et le but de cette sous-section est de montrer que le lemme de \ref{definition} r\'esulte du r\'esultat plus pr\'ecis:
\nl
\bf Lemme: \sl $Jac_{x\in {\cal E}}\pi_{b}(\tilde{\psi},\tilde{\epsilon},\tilde{\underline{\ell}},\tilde{\underline{\eta}})$ est non nul exactement quand $\ell_{1}=\ell_{2}$ et $\eta_{1}=\eta_{2}(-1)^{A-B}$.
\nl
\rm
On admet donc ce lemme et on d\'emontre \ref{nullite}. On note $\psi_{<}$ le morphisme  d\'efini par:
$$
Jord(\psi_{<})=\{ b(\rho',A',B',\zeta'); (\rho',A',B',\zeta')\in Jord(\tilde{\psi}), (\rho',A',B',\zeta')<(\rho,A+T_{1},B +T_{1},\zeta)\}$$
$$\cup \{(\rho',A',B',\zeta')\in Jord(\tilde{\psi}), (\rho',A',B',\zeta')\geq (\rho,A+T_{1},B +T_{1},\zeta)\}.
$$
On d\'efinit $\psi_{\leq}$ en posant
$$
Jord(\psi_{\leq})=\{ b(\rho',A',B',\zeta'); (\rho',A',B',\zeta')\in Jord(\tilde{\psi}), (\rho',A',B',\zeta')\leq(\rho,A+T_{2},B +T_{2},\zeta)\}$$
$$\cup \{(\rho',A',B',\zeta')\in Jord(\tilde{\psi}), (\rho',A',B',\zeta')> (\rho,A+T_{2},B +T_{2},\zeta)\}.
$$
On a d\'efini les ensembles ${\cal C}(\zeta,A,B,T_{i})$ pour $i=1,2$ et on a par d\'efinition
$$
\pi_{b}(\psi_{\leq },\tilde{\epsilon},\tilde{\underline{\ell}},\tilde{\underline{\eta}})=Jac_{x\in {\cal C}(\zeta,A,B,T_{2})}Jac_{x\in {\cal C}(\zeta,A,B,T_{1})}\pi_{b}(\psi_{< },\tilde{\epsilon},\tilde{\underline{\ell}},\tilde{\underline{\eta}}).
$$
Par d\'efinition $\pi_{b}(\psi,\tilde{\epsilon},\tilde{\underline{\ell}},\tilde{\underline{\eta}})$ est un module de Jacquet convenable de $\pi_{b}(\psi_{\leq },\tilde{\epsilon},\tilde{\underline{\ell}},\tilde{\underline{\eta}})$ et il suffit donc de d\'emontrer que $\pi_{b}(\psi_{\leq },\tilde{\epsilon},\tilde{\underline{\ell}},\tilde{\underline{\eta}})=0$ si les conditions du lemme ne sont pas remplies.
On utilise les repr\'esentations $S(\zeta,A,B,T_{i})$ de \ref{notations} et on a (cf. \ref{consequence}) une inclusion
$$
\pi_{b}(\psi_{<},\tilde{\epsilon},\tilde{\underline{\ell}},\tilde{\underline{\eta}})\hookrightarrow
S(\zeta,A,B,T_{1})\times S(\zeta,A,B,T_{2}) \times \pi_{b}(\psi_{\leq },\tilde{\epsilon},\tilde{\underline{\ell}},\tilde{\underline{\eta}}).
$$
Le point ici  est que l'induite $S(\zeta,A,B,T_{2})\times S(\zeta,A,B,T_{1})$ est irr\'eductible; ceci sera d\'emontr\'e en \ref{unresultatdirreductibilite}. On peut donc \'echanger les 2 termes de cette induite. On utilise ensuite le fait que
$$
S(\zeta,A,B,T_{2})\hookrightarrow S(\zeta,A+T_{1},B+T_{1},T_{2}-T_{1})\times S(\zeta,A,B,T_{1}).$$
Cela montre l'existence d'une inclusion:
$$
\pi_{b}(\psi_{<},\tilde{\epsilon},\tilde{\underline{\ell}},\tilde{\underline{\eta}})\hookrightarrow
S(\zeta,A+T_{1},B+T_{1},T_{2}-T_{1})\times S(\zeta,A,B,T_{1})\times S(\zeta,A,B,T_{1})\times 
\pi_{b}(\psi_{\leq },\tilde{\epsilon},\tilde{\underline{\ell}},\tilde{\underline{\eta}}).
$$
En particulier pour que $\pi_{b}(\psi_{\leq },\tilde{\epsilon},\tilde{\underline{\ell}},\tilde{\underline{\eta}})$ soit non nul, il faut que 
$$
Jac_{x\in {\cal C}(\zeta,A+T_{1},B+T_{1},T_{2}-T_{1})}\pi_{b}(\psi_{<},\tilde{\epsilon},\tilde{\underline{\ell}},\tilde{\underline{\eta}})\eqno(1)
$$
soit non nul. On sait que $\pi_{b}(\psi_{<},\tilde{\epsilon},\tilde{\underline{\ell}},\tilde{\underline{\eta}})$ est un module de Jacquet convenable de $\pi(\tilde{\psi},\tilde{\epsilon},\tilde{\underline{\ell}},\tilde{\underline{\eta}})$. Mais comme $T_{1}$ est suppos\'e grand, il n'y a pas de liaison entre l'ensemble 
$ {\cal C}(\zeta,A+T_{1},B+T_{1},T_{2}-T_{1})$ et l'ensemble qui sert \`a passer de $\pi(\tilde{\psi},\tilde{\epsilon},\tilde{\underline{\ell}},\tilde{\underline{\eta}})$ \`a $\pi_{b}(\psi_{<},\tilde{\epsilon},\tilde{\underline{\ell}},\tilde{\underline{\eta}})$. On peut donc commuter les op\'erations et la non nullit\'e de 
(1) n\'ecessite que 
$
Jac_{x\in {\cal C}(\zeta,A+T_{1},B+T_{1},T_{2}-T_{1})}\pi_{b}(\tilde{\psi},\tilde{\epsilon},\tilde{\underline{\ell}},\tilde{\underline{\eta}})\neq 0
.
$
Maintenant, nous sommes dans la situation du lemme ci-dessus qui donne exactement la condition n\'ecessaire de non nullit\'e. Cela termine la r\'eduction.
%%%%%%
\subsection{Preuve du crit\`ere de \ref{reduction}\label{nullitenonnullite}}
On simplifie les notations; on fixe $\psi$ un morphisme de restriction discr\`ete \`a la diagonale et $(\rho,A,B,\zeta)\in Jord(\psi)$. Et on suppose qu'il existe un entier $T$ tel que l'\'el\'ement $(\rho,A+T,B+T,\zeta)$ soit aussi dans $Jord(\psi)$ est soit exactement l'\'el\'ement suivant $(\rho,A,B,\zeta)$. On fixe $\underline{\ell}$ et $\underline{\eta}$ de fa\c{c}on \`a d\'efinir la repr\'esentation irr\'eductible
$\pi(\psi,\epsilon,\underline{\ell},\underline{\eta})$. On pose encore, $\ell_{1}=\underline{\ell}(\rho,A,B,\zeta)$ et $\ell_{2}=\underline{\ell}(\rho,A+T,B+T,\zeta)$ et on d\'efinit de m\^eme $\eta_{i}$ pour $i=1,2$ en rempla\c{c}ant $\ell$ par $\eta$. Et on d\'emontre
\nl
\bf Lemme: \sl $Jac_{x\in {\cal C}(\zeta,A,B,T)}\pi(\psi,\epsilon,\underline{\ell},\underline{\eta})\neq 0$ exactement quand $\ell_{1}=\ell_{2}$ et $\eta_{1}=\eta_{2}(-1)^{A-B}$.\rm
\nl
La d\'emonstration se fait par r\'ecurrence sur $A-B$. Le cas o\`u $A=B$ est facile: ici n\'ecessairement $\ell_{1}=\ell_{2}=0$ et $\eta_{1}=\epsilon(\rho,A,B,\zeta)$ tandis que $\eta_{2}=\epsilon(\rho,A+T,B+T,\zeta)$. Le lemme s'exprime alors ainsi
$$
Jac_{\zeta (B+T), \cdots, \zeta (B+1)}\pi(\psi,\epsilon,\underline{\ell},\underline{\eta})\neq 0
$$
exactement quand $\epsilon(\rho,A,B,\zeta)$ $=$ $\epsilon(\rho,A+T,B+T,\zeta)$. Si $\psi$ est \'el\'ementaire c'est exactement la  d\'efinition; on est aussi dans la situation d'un morphisme \'el\'ementaire si pour tout $(\rho',A',B',\zeta')\in Jord(\psi)$, $\underline{\ell}(\rho',A',B',\zeta')=0$. Supposons qu'il n'en soit pas ainsi et fixons $(\rho',A',B',\zeta')\in Jord(\psi)$ tel que $\underline{\ell}(\rho',A',B',\zeta')\neq 0$. Alors on sait que $\pi(\psi,\epsilon,\underline{\ell},\underline{\eta})$ est l'unique sous-module irr\'eductible de l'induite:
$$
<\rho'\vert\,\vert^{\zeta' B'}, \cdots, \rho'\vert\,\vert^{-\zeta' A'}>\times \pi(\psi',\epsilon,\underline{\ell}',\underline{\eta}),
$$
o\`u $\psi'$ se d\'eduit de $\psi$ en rempla\c{c}ant $(\rho',A',B',\zeta')$ par $(\rho',A'-1,B'+1,\zeta')$ et $\underline{\ell}'$ est $\underline{\ell}$ sauf  $\underline{\ell}'(\rho',A'-1,B'+1,\zeta')=\underline{\ell}(\rho',A',B',\zeta')-1$; les autres fonctions se d\'eduisent naturellement. Comme $\psi$ est suppos\'e de restriction discr\`ete \`a la diagonale, on a certainement soit $A'<B$ soit $B'>A+T=B+T$. Avec cela, on v\'erifie l'isomorphisme
$$
<\rho'\vert\,\vert^{\zeta' B'}, \cdots, \rho'\vert\,\vert^{-\zeta' A'}>\times <\rho\vert\,\vert^{\zeta (B+T)}, \cdots, \rho\vert\,\vert^{\zeta(B+1)}>\simeq$$
$$ <\rho\vert\,\vert^{\zeta (B+T)}, \cdots, \rho\vert\,\vert^{\zeta(B+1)}>\times <\rho'\vert\,\vert^{\zeta' B'}, \cdots, \rho'\vert\,\vert^{-\zeta' A'}>
$$
 dans le $GL$ convenable: si $\zeta'=\zeta$, les segments associ\'es ont m\^eme propri\'et\'e de croissance et comme $B' \notin [B+T,B]$, ils ne sont pas li\'es et si $\zeta'=-\zeta$, tout \'el\'ement de $[\zeta (B+T),\zeta (B+1)]$  vu comme segment r\'eduit \`a un \'el\'ement, n'est pas li\'e au segment $[\zeta' B',-\zeta' A']$. D'o\`u le fait que l'on peut \'echanger les 2 repr\'esentations dans l'induite. Ceci entra\^{\i}ne l'\'equivalence:
 $$
 Jac_{\zeta (B+T), \cdots, \zeta(B+1)}\pi(\psi,\epsilon,\underline{\ell},\underline{\eta})\neq 0 \Leftrightarrow
 Jac_{\zeta (B+T), \cdots, \zeta(B+1)} \pi(\psi',\epsilon,\underline{\ell}',\underline{\eta})\neq 0.
 $$
 On se ram\`ene ainsi au cas o\`u $\psi$ est \'el\'ementaire d\'ej\`a trait\'e.
 
 On suppose donc maintenant que $A>B$. On note $\psi'$ le morphisme qui se d\'eduit de $\psi$ en enlevant les 2 blocs $(\rho,A,B,\zeta)$ et $(\rho,A+T,B+T,\zeta)$.

On suppose d'abord que $\ell_{1}>0$.  On sait que $\pi(\psi,\epsilon,\underline{\ell},\underline{\eta})$ est l'unique sous-module irr\'eductible de l'induite:
$$
<\rho\vert\,\vert^{\zeta B}, \cdots, \rho\vert\,\vert^{-\zeta A}> \times \pi(\psi',\epsilon',\underline{\ell}',\underline{\eta}', (\rho,A-1,B+1,\zeta,\ell_{1}-1,\eta_{1}), (\rho,A+T,B+T,\zeta,\ell_{2},\eta_{2})). \eqno(1)
$$ 
On  calcule $Jac_{x\in {\cal C}(\zeta,A,B,T)}$ du terme de droite de (1). On le fait ligne par ligne c'est-\`a-dire que l'on consid\`ere pour $j \in [0,T]$ d\'ecroissant, l'ensemble ${\cal C}(\zeta, A+j ,B+j,T-j)$ qui n'est autre que les $T-j$ premi\`eres lignes du tableau repr\'esentant ${\cal C}(\zeta,A,B,T)$. Pour $j=T$, l'ensemble que l'on vient de d\'efinir est vide. Et progressivement sur $j$ d\'ecroissant de $T$ \`a $0$, on va montrer que $Jac_{x\in {\cal C}(\zeta, A+j,B+j,T-j)}$ du terme de droite de (2) vaut:
$$
<\rho\vert\,\vert^{\zeta B}, \cdots, \rho\vert\,\vert^{-\zeta A}> \times 
$$
$$ Jac_{x\in {\cal C}(\zeta, A+j,B+j,T-j)}\pi(\psi',\epsilon', (\rho,A-1,B+1,\zeta,t_{1}-1,\eta_{1}),(\rho,A+T,B+T,\zeta, t_{2},\eta_{2})). \eqno (2)_{j}
$$Pour $j=T$, il n'y a rien \`a d\'emontrer.
On l'admet pour $j \geq 1$ et on le d\'emontre pour $j-1$; comme $j\geq 1$, par d\'efinition
 $$Jac_{x\in {\cal C}(\zeta, A+j,B+j,T-j)}\pi(\psi',\epsilon', (\rho,A-1,B+1,\zeta,\ell_{1}-1,\eta_{1}),(\rho,A+T,B+T,\zeta, \ell_{2},\eta_{2}))=$$
 $$\pi(\psi',\epsilon',\underline{t}',\underline{\eta}', (\rho,A-1,B+1,\zeta, \ell_{1}-1,\eta_{1})(\rho,A+j,B+j,\zeta,\ell_{2},\eta_{2})).$$
On passe de $j$ \`a $j-1$ en rajoutant la ligne $\zeta (B+j), \cdots, \zeta (A+j)$.  Il faut donc calculer $Jac_{x\in [\zeta (B+j), \zeta (A+j)]}$ du r\'esultat pour $j$ pour obtenir  $Jac_{x\in {\cal C}(\zeta,A+j-1,B+j-1,T-j+1)}$. 

On d\'ecompose donc l'intervalle $[\zeta (B+j),\zeta(A+j)]$ en trois sous-ensembles ${\cal E}_{i}$ pour $i\in [1,3]$ tel que:
$
Jac_{x\in {\cal E}_{1}}Jac^d_{x\in -^t{\cal E}_{2}}<\rho\vert\,\vert^{\zeta B}, \cdots, \rho\vert\,\vert^{-\zeta A}>\neq 0
$ et 
$$
Jac_{x\in {\cal E}_{3}}\pi(\psi',\epsilon',\underline{\ell}',\underline{\eta}', (\rho,A-1,B+1,\zeta, \ell_{1}-1,\eta_{1})(\rho,A+j,B+j,\zeta,\ell_{2},\eta_{2}))\neq 0.$$
L'\'el\'ement $\zeta (A+j)$ ne peut \^etre que dans ${\cal E}_{3}$ pour des raisons de support cuspidal. On note $\zeta B_{0}$ l' \'el\'ement de ${\cal E}_{3}$ qui satisfait $[\zeta B_{0},\zeta (A+j)]\in {\cal E}_{3}$ et $\zeta (B_{0}-1)\notin {\cal E}_{3}$; un tel \'el\'ement existe n\'ecessairement. Comme $\zeta (B_{0}-1)\notin {\cal E}_{3}$, la non nullit\'e relative \`a ${\cal E}_{3}$ n\'ecessite que
$$
Jac_{x\in [\zeta (B_{0}), \zeta (A+j)]}\pi(\psi',\epsilon',\underline{\ell}',\underline{\eta}', (\rho,A-1,B+1,\zeta, \ell_{1}-1,\eta_{1})(\rho,A+j,B+j,\zeta,\ell_{2},\eta_{2}))\neq 0.
$$
D'apr\`es \ref{jacquetsimple}, il exite donc $(\rho,A'',B'',\zeta'')\in Jord(\psi')\cup \{(\rho,A-1,B+1,\zeta),(\rho,A+j,B+j,\zeta)\}$ tel que $\zeta B_{0}=\zeta'' B''$. Comme $B_{0}\in [B,A+T]$ et que l'on sait que $\psi$ est de restriction discr\`ete \`a la diagonale, $(\rho,A'',B'',\zeta'')$ ne peut \^etre que $(\rho,B+1,A-1,\zeta)$ ou $(\rho,B+j,A+j,\zeta)$. Si $j=1$, les 2 possibilit\'es sont les m\^emes et si $j>1$, seul la deuxi\`eme est possible. Donc dans tous les cas, $B_{0}=B+j$ et ${\cal E}_{3}={\cal E}$. On a donc bien prouv\'e (2)$_{j}$ pour tout $j$.

On suppose maintenant que l'on a aussi $\ell_{2}\neq 0$; on continue les calculs pr\'ec\'edents en utilisant l'inclusion: $$\pi(\psi',\epsilon',\underline{\ell}',\underline{\eta}',(\rho,A-1,B+1,\zeta,\ell_{1}-1,\eta_{1})(\rho, A+T,B+T,\ell_{2},\eta_{2})) \hookrightarrow
<\rho\vert\,\vert^{\zeta (B+T)}, \cdots, \rho\vert\,\vert^{-\zeta (A+T)}> $$
$$\times \pi(\psi',\epsilon',\underline{\ell}',\underline{\eta}', (\rho,A-1,B+1,\zeta,\ell_{1}-1,\eta_{1}),(\rho,A+T-1,B+T+1,\zeta, \ell_{2}-1,\eta_{2})).\eqno(3)
$$
On garde les notations pr\'ec\'edentes et ici on d\'emontre de proche en proche sur $j$ d\'ecroissant que 
$$
Jac_{x\in {\cal C}(\zeta,A+j,B+j,T-j)} (3)= <\rho\vert\,\vert^{\zeta (B+j)}, \cdots, \rho\vert\,\vert^{-\zeta (A+j)}> \times$$
$$ \pi(\psi',\epsilon',\underline{\ell}',\underline{\eta}',(\rho,A-1,B+1,\zeta,\ell_{1}-1,\eta_{1}),
(\rho,A-1+j,B+1+j,\zeta,\ell_{2}-1,\eta_{2})).\eqno(4)_{j}$$
On suppose cette expression connue pour $j \geq 1$ et on la d\'emontre pour $j-1$, c'est \`a dire qu'il faut calculer $Jac_{x \in [\zeta (B+j),\zeta (A+j)]}$ du terme de droite (4)$_{j}$. On d\'ecompose encore $[\zeta (B+j),\zeta (A+j)]$ en 3 sous-ensembles ${\cal E}_{i}$ pour $i\in [1,3]$ tel que:
$$
Jac_{x\in {\cal E}_{1}}Jac^d_{x\in -^t{\cal E}_{2}} <\rho\vert\,\vert^{\zeta (B+j)}, \cdots, \rho\vert\,\vert^{-\zeta (A+j)}> \neq 0
$$
$$
Jac_{x\in {\cal E}_{3}}\pi(\psi',\epsilon',\underline{\ell}',\underline{\eta}',(\rho,A-1,B+1,\zeta,\ell_{1}-1,\eta_{1}),
(\rho,A-1+j,B+1+j,\zeta,\ell_{2}-1,\eta_{2}))\neq 0.
$$
L'\'el\'ement $\zeta (B+j)$ ne peut \^etre dans ${\cal E}_{3}$ par les propri\'et\'es des blocs de Jordan de $\psi'$ (cf. \ref{proprietesgenerales}); ainsi il est dans ${\cal E}_{1}$ car il ne peut \^etre le premier \'el\'ement de ${\cal E}_{2}$; si ${\cal E}_{1}$ n'est pas r\'eduit \`a $\zeta(B+j)$ il contient $\zeta (B+j-1)$; c'est impossible puisque cet \'el\'ement n'est pas dans l'intervalle. Ainsi ${\cal E}_{1}=\{\zeta (B+j)\}$. L'ensemble ${\cal E}_{2}$ s'il n'est pas vide a comme premier \'el\'ement $\zeta (A+j)$ et il sera donc alors n\'ecessairement r\'eduit \`a cet \'el\'ement; donc ${\cal E}_{3}$ contient au moins $[\zeta( B+j+1), \zeta (A+j-1)]$. Ainsi le module de Jacquet cherch\'e est
$$
Jac_{\zeta (A+j)}\biggl(<\rho\vert\,\vert^{\zeta (B+j-1)}, \cdots, \rho\vert\,\vert^{-\zeta (A+j)}>\times $$
$$Jac_{\zeta(B+j+1), \cdots, \zeta (A+j-1)} \pi(\psi',\epsilon',\underline{\ell}',\underline{\eta}',(\rho,A-1,B+1,\zeta,\ell_{1}-1,\eta_{1}),
(\rho,A-1+j,B+1+j,\zeta,\ell_{2}-1,\eta_{2}))\biggr)
$$
et la deuxi\`eme ligne n'est autre par d\'efinition que $$ \pi(\psi',\epsilon',\underline{\ell}',\underline{\eta}',(\rho,A-1,B+1,\zeta,\ell_{1}-1,\eta_{1}),
(\rho,A-2+j,B+j,\zeta,\ell_{2}-1,\eta_{2})).$$ D'apr\`es \ref{jacquetsimple}, le $Jac_{\zeta(A+j)}$ de cette repr\'esentation est nulle et c'est donc ${\cal E}_{2}$ qui contient $\zeta (A+j)$. On obtient alors (2)$_{j-1}$, comme cherch\'e.

Avec ces 2 \'etapes, on a montr\'e que si $\ell_{1}\ell_{2}\neq 0$, $Jac_{x\in {\cal C}(\zeta,A,B,T)}\pi(\psi,\epsilon,\underline{\ell},\underline{\eta})$ est soit nul soit est un sous-module de:
$$
<\rho\vert\,\vert^{\zeta B}, \cdots, \rho\vert\,\vert^{-\zeta A}>\times <\rho\vert\,\vert^{\zeta B}, \cdots, \rho\vert\,\vert^{-\zeta A}>
$$
$$
\times Jac_{x\in{\cal C}(\zeta,A-1,B+1,T)}\pi(\psi',\epsilon',\underline{t}',\underline{\eta}', (\rho,A-1,B+1,\zeta,t_{1}-1,\eta_{1}),(\rho,A+T-1,B+T+1,\zeta,t_{2}-1,\eta_{2})).\eqno(4)
$$
On peut appliquer la proposition par r\'ecurrence pour savoir que (4) est non nul pr\'ecis\'ement quand $\ell_{1}-1=\ell_{2}-1$ et $\eta_{1}\eta_{2}=(-1)^{A-B+2}=(-1)^{A-B}$. Ce sont exactement les conditions n\'ecessaires de l'\'enonc\'e pour la non nullit\'e.

Reste \`a d\'emontrer qu'elles sont suffisantes; cela vient \`a remonter les constructions pr\'ec\'edentes. On suppose donc que $\ell_{1}=\ell_{2}$ et que $\eta_{1}\eta_{2}=(-1)^{A-B}$. Par r\'ecurrence, on sait alors que 
$$
\pi(\psi',\epsilon',\underline{\ell}',\underline{\eta}', (\rho,A-1,B+1,\zeta,\ell_{1}-1,\eta_{1}),(\rho,A+T-1,B+T+1,\zeta,\ell_{2}-1,\eta_{2})) \hookrightarrow  S(\zeta,A-1,B+1,T) \times \tau, $$
o\`u $\tau$ est une repr\'esentation convenable, irr\'eductible non nulle.

D'o\`u par d\'efinition, l'inclusion de $\pi(\psi,\epsilon,\underline{\ell},\underline{\eta})$ dans
$$
<\rho\vert\,\vert^{\zeta B}, \cdots, \rho\vert\,\vert^{-\zeta A}> \times <\rho\vert\,\vert^{\zeta (B+T)}, \cdots, \rho\vert\,\vert^{-\zeta (A+T)}>\times S(\zeta,A-1,B+1,T)\times \tau. \eqno(5)
$$
Il faut rappeler que la repr\'esentation $S(\zeta,A-1,B+1,T)$ est associ\'ee au tableau:
$$
\begin{matrix}
&\zeta (B+T+1) &\cdots &\zeta (A+T-1)\\
&\vdots &\vdots&\vdots\\
&\zeta (B+2) &\cdots &\zeta A
\end{matrix}
$$
Ce sont donc des formules standard qui prouvent l'irr\'eductibilit\'e de l'induite dans un GL convenable de
$$
<\rho\vert\,\vert^{\zeta B}, \cdots, \rho\vert\,\vert^{-\zeta A}>\times S(\rho,A-1,B+1,T)$$et de $<\rho\vert\,\vert^{-\zeta (A+1)}, \cdots, \rho\vert\,\vert^{-\zeta (A+T)}>\times S(\rho,A-1,B+1,T)$. 
On utilise l'inclusion
$$<\rho\vert\,\vert^{\zeta B}, \cdots, \rho\vert\,\vert^{-\zeta A}>\times
<\rho\vert\,\vert^{\zeta (B+T)}, \cdots, \rho\vert\,\vert^{-\zeta (A+T)}>\hookrightarrow$$
$$<\rho\vert\,\vert^{\zeta B}, \cdots, \rho\vert\,\vert^{-\zeta A}>\times <\rho\vert\,\vert^{\zeta (B+T)}, \cdots, \rho\vert\,\vert^{-\zeta A}> \times  <\rho\vert\,\vert^{-\zeta (A+1)}, \cdots, \rho\vert\,\vert^{-\zeta (A+T)}>$$
$$\simeq <\rho\vert\,\vert^{\zeta (B+T)}, \cdots, \rho\vert\,\vert^{-\zeta A}> \times <\rho\vert\,\vert^{\zeta B}, \cdots, \rho\vert\,\vert^{-\zeta A}>\times <\rho\vert\,\vert^{-\zeta (A+1)}, \cdots, \rho\vert\,\vert^{-\zeta (A+T)}>$$
$$\hookrightarrow 
<\rho\vert\,\vert^{\zeta (B+T)}, \cdots, \rho\vert\,\vert^{\zeta (B+1)}>\times <\rho\vert\,\vert^{\zeta B}, \cdots, \rho\vert\,\vert^{-\zeta A}>\times <\rho\vert\,\vert^{\zeta B}, \cdots, \rho\vert\,\vert^{-\zeta A}>\times $$
$$<\rho\vert\,\vert^{-\zeta (A+1)}, \cdots, \rho\vert\,\vert^{-\zeta (A+T)}>;
$$
et l'inclusion (5) se prolonge en une inclusion  dans l'induite
$$<\rho\vert\,\vert^{\zeta (B+T)}, \cdots, \rho\vert\,\vert^{\zeta (B+1)}>\times S(\rho,A-1,B+1,T)
$$
$$
\times <\rho\vert\,\vert^{\zeta B}, \cdots, \rho\vert\,\vert^{-\zeta A}> \times <\rho\vert\,\vert^{\zeta B}, \cdots, \rho\vert\,\vert^{-\zeta A}> \times <\rho\vert\,\vert^{-\zeta (A+1)}, \cdots, \rho\vert\,\vert^{-\zeta (A+T)}> \times \tau.
$$
On a besoin de savoir que
$$<\rho\vert\,\vert^{-\zeta (A+1)}, \cdots, \rho\vert\,\vert^{-\zeta (A+T)}> \times \tau \simeq 
<\rho\vert\,\vert^{\zeta (A+T), \cdots, \rho\vert\,\vert^{\zeta (A+1)}}>\times \tau.
$$
Pour cela il suffit de savoir que $\rho\vert\,\vert^x \times \tau \simeq \rho\vert\,\vert^{-x}\times \tau $ pour tout $x \in [ (A+1),(A+T)]$. Or comme $Jac_{x}\tau=$ $0$, le terme de gauche a un unique sous-module irr\'eductible et il intervient avec multiplicit\'e 1 comme sous-quotient de l'induite $\rho\vert\,\vert^x \times \tau$. De m\^eme $Jac_{-x}\tau=0$ et il suffit donc de montrer que le sous-module irr\'eductible du terme de gauche est le sous-module irr\'eductible du terme de droite. Le plus simple ici est d'utiliser une description de $\tau$ comme sous-module de l'induite $S(\rho,A-1,B+1,\zeta)\times \pi(\psi',\epsilon')$ et donc de la proposition \ref{induction}  par r\'ecurrence; on peut donc admettre cela. Par hypoth\`ese  $\pi(\psi',\epsilon')$ s'obtient \`a partir d'un morphisme de restriction discr\`ete. Dans le cas de restriction discr\`ete \`a la diagonale, on a la propri\'et\'e plus g\'en\'erale:

Soit $x>0$ tel que $2x-1$ ne soit pas un bloc de Jordan de la restriction \`a la diagonale d'un morphisme $\tilde{\psi}'$. Alors $\rho\vert\,\vert^x \times \pi(\tilde{\psi}',\tilde{\epsilon}')$ est irr\'eductible pour tout choix de $\tilde{\epsilon}'$ (cf. \cite{paquetsdiscrets}, 6.4)

\

 Par les propri\'et\'es standard dans les groupes lin\'eaires, l'induite $$
<\rho\vert\,\vert^{\zeta B}, \cdots, \rho\vert\,\vert^{-\zeta A}> \times <\rho\vert\,\vert^{\zeta (A+T)}, \cdots, \rho\vert\,\vert^{\zeta (A+1)}>$$
est irr\'eductible (en effet pour tout $z\in [\zeta B, -\zeta A]$ et tout $z'\in [\zeta (A+T),\zeta (A+1)]$, $\vert z-z'\vert \geq \vert B-A-1\vert >1$).

On trouve ainsi une inclusion de $\pi(\psi,\epsilon,\underline{\ell},\underline{\eta})$ dans
$$
<\rho\vert\,\vert^{\zeta(B+T)}, \cdots, \rho\vert\,\vert^{-\zeta A}> \times S(\rho,A-1,B+1,T) \times  <\rho\vert\,\vert^{\zeta (A+T)}, \cdots, \rho\vert\,\vert^{\zeta (A+1)}> 
$$
$$
\times <\rho\vert\,\vert^{\zeta B}, \cdots, \rho\vert\,\vert^{-\zeta A}> \times \tau.
$$
On remplace encore  l'induite $<\rho\vert\,\vert^{\zeta(B+T)}, \cdots, \rho\vert\,\vert^{-\zeta A}> $ par son inclusion dans $$<\rho\vert\,\vert^{\zeta (B+T)}, \cdots, \rho\vert\,\vert^{\zeta (B+1)}> \times <\rho\vert\,\vert^{\zeta B}, \cdots, \rho\vert\,\vert^{-\zeta A}> $$ et on peut encore faire commuter comme ci-dessus pour obtenir une inclusion de $\pi(\psi,\epsilon,\underline{\ell},\underline{\eta})$ dans
$$
<\rho\vert\,\vert^{\zeta (B+T)}, \cdots, \rho\vert\,\vert^{\zeta (B+1)}>\times S(\rho,A-1,B+1,T) \times 
<\rho\vert\,\vert^{\zeta (A+T)}, \cdots, \rho\vert\,\vert^{\zeta (A+1)}> \eqno(5)
$$
$$
\times <\rho\vert\,\vert^{\zeta B}, \cdots, \rho\vert\,\vert^{-\zeta A}> \times <\rho\vert\,\vert^{\zeta B}, \cdots, \rho\vert\,\vert^{-\zeta A}> \times \tau.
$$
La  non nullit\'e de $Jac_{x\in {\cal C}(\zeta,A,B,T)}\pi(\psi,\epsilon,\underline{\ell},\underline{\eta})$ vient de la r\'eciprocit\'e de Frobenius et de (5).
D'o\`u la non nullit\'e annonc\'ee dans l'\'enonc\'e. 

\

On suppose maintenant que $\ell_{2}=0$ en gardant l'hypoth\`ese que $\ell_{1}\neq 0$. Montrons que le module de Jacquet cherch\'e est 0.  On revient au d\'ebut de la d\'emonstration pr\'ec\'edente. On consid\`ere le sous-tableau de ${\cal C}$ form\'e des $T-1$ premi\`eres lignes et des $A-B-1$ premi\`eres colonnes; ceci n'est autre que ${\cal C}(\zeta, A-1,B+1,T-1)$. D'autre part comme $\ell_{2}=0$, le bloc $(\rho,A+T,B+T,\zeta,\ell_{2},\eta_{2})$ peut \^etre remplac\'e par l'ensemble des blocs $\cup_{C\in [A+T,B+T]}(\rho,C,C,\zeta,\eta_{2}(-1)^{C-B-T})$ et  on peut donc remplacer $(\rho,A+T,B+T,\zeta,\ell_{2},\eta_{2})$ par les 3 blocs $(\rho,A+T-2,B+T,\zeta,t=0, \eta_{2}),$ $ (\rho,A+T,A+T,\zeta,\eta(-1)^{A-B}),$ $(\rho,A+T-1,A+T-1,\zeta,\eta_{2} (-1)^{A-B-1})$. Raisonnons par l'absurde. Le calcul d\'ej\`a fait montre que  la non nullit\'e cherch\'ee n\'ecessite aussi que
$$
Jac_{x\in {\cal C}(\zeta,A-1,B+1,T-1)} \pi(\psi',\epsilon', (\rho,A-1,B+1,\zeta,\ell_{1}-1,\eta_{1}), (\rho,A+T-2,B+T,\zeta,\ell_{2}=0, \eta), $$
$$(\rho,A+T,A+T,\zeta,\eta_{2}(-1)^{A-B}), (\rho,A+T-1,A+T-1,\zeta,\eta_{2} (-1)^{A-B-1}))\neq 0.$$
Par r\'ecurrence cela force $\ell_{1}-1=0$ et la condition $\eta_{2}=\eta_{1}(-1)^{A-B-2}$.

On d\'efinit le tableau:
$$
{\cal T}:=\begin{matrix}
&&&\zeta B &\cdots &-\zeta A\\
&\zeta (B+T)&\cdots &\zeta (B+1) &\cdots &-\zeta(A-1)\\
&\vdots &\vdots&\vdots&\vdots&\vdots\\
&\zeta (A+T-2) &\cdots &\zeta (A-1) &\cdots &-\zeta (B+1)\\
&\zeta (A+T-1) &\cdots &\zeta A\\
&\zeta (A+T) &\cdots &\zeta (A+1).
\end{matrix}
$$
Il lui correspond une unique repr\'esentation irr\'eductible du $GL$ convenable, bas\'ee sur $\rho$ que l'on note $\sigma_{\cal T}$. On veut montrer que $\pi(\psi,\epsilon,\underline{\ell},\underline{\eta})$ est alors un sous-module irr\'eductible de
$$\sigma_{\cal T}\times \pi(\underline{\psi}',(\rho,A-1,A-1,\zeta,\eta(-1)^{A-B-1}),(\rho,A,A,\zeta,\eta(-1)^{A-B})).
$$Si on admet cela, on peut conclure \`a la nullit\'e de $Jac_{x\in {\cal C}(\zeta,A,B,T)}\pi(\psi,\epsilon,\underline{\ell},\underline{\eta})$ car ceci est d\'ej\`a vrai pour l'induite ci-dessus: en effet, il suffit de remarquer que $Jac_{\zeta (B+T), \cdots, \zeta (B+1)}$ de cette induite est nul (on a pris la premi\`ere colonne de ${\cal C}(\zeta,A,B,T)$) et ceci r\'esulte des formules standards.

Il faut donc montrer cette inclusion. Comme on sait d\'ej\`a que $\ell_{1}=1$ et $\ell_{2}=0$ avec $\eta_{1}\eta_{2}=(-1)^{A-B}$, on peut \'ecrire $\pi(\psi,\epsilon,\underline{\ell},\underline{\eta})$ comme sous-module irr\'eductible de l'induite
$$
<\rho\vert\,\vert^{\zeta B}, \cdots, \rho\vert\,\vert^{-\zeta A}>\times \biggl(\begin{matrix}
&\zeta (B+T)&\cdots &\zeta (B+1) &\cdots &-\zeta(A-1)\\
&\vdots &\vdots&\vdots&\vdots&\vdots\\
&\zeta (A+T-2) &\cdots &\zeta (A-1) &\cdots &-\zeta (B+1)
\end{matrix}\biggr) $$
$$\times \pi(\psi',\epsilon',\underline{\ell}',\underline{\eta}',(\rho,A+T-1,A+T-1,\zeta,-\eta),(\rho,A+T,A+T,\zeta,\eta)).
$$
On peut remplacer la premi\`ere ligne ci-dessus par le sous-module irr\'eductible de l'induite \'ecrite car sinon on aurait une inclusion dans $$<\rho\vert\,\vert^{\zeta B}, \cdots, \rho\vert\,\vert^{-\zeta (A-1)}> \times
\biggl(\begin{matrix}
&\zeta (B+T)&\cdots &\zeta (B+1) &\cdots &-\zeta(A-1)\\
&\vdots &\vdots&\vdots&\vdots&\vdots\\
&\zeta (A+T-2) &\cdots &\zeta (A-1) &\cdots &-\zeta (B+1)
\end{matrix}\biggr)\times 
$$
$$
\rho\vert\,\vert^{-\zeta A}\times \pi(\psi',\epsilon',\underline{\ell}',\underline{\eta}',(\rho,A+T-1,A+T-1,\zeta,-\eta),(\rho,A+T,A+T,\zeta,\eta)).
$$
Mais l'induite de cette deuxi\`eme ligne est isomorphe \`a $$\rho\vert\,\vert^{\zeta A}\times \pi(\underline{\psi}',(\rho,A+T-1,A+T-1,\zeta,-\eta),(\rho,A+T,A+T,\zeta,\eta))$$ et on v\'erifie alors que l'on aurait $Jac_{\zeta A}\pi(\psi,\epsilon,\underline{\ell},\underline{\eta})\neq 0.$ Ceci est impossible d'apr\`es \ref{jacquetsimple}. C'\'etait le point cl\'e. La suite r\'esulte de l'inclusion:
$$\pi(\psi',\epsilon',\underline{\ell}',\underline{\eta}',(\rho,A+T-1,A+T-1,\zeta,-\eta),(\rho,A+T,A+T,\zeta,\eta))\hookrightarrow $$
$$<\rho\vert\,\vert^{\zeta (A+T-1)}, \cdots, \rho\vert\,\vert^{\zeta A}>\times <\rho\vert\,\vert^{\zeta (A+T)}, \cdots,\rho\vert\,\vert^{\zeta (A+1)}>\times$$
$$
\pi(\underline{\psi'},(\rho,A-1,A-1,\zeta,-\eta),(\rho,A,A,\zeta,\eta))$$et  des factorisations usuelles.

\

Consid\'erons maintenant le cas o\`u $\ell_{1}=0$. On regarde d'abord le cas o\`u  $\ell_{2}\neq 0$.  Ici, 
$$
\pi(\psi,\epsilon,\underline{\ell},\underline{\eta})\hookrightarrow
<\rho\vert\,\vert^{\zeta(B+T)}, \cdots, \rho\vert\,\vert^{-\zeta (A+T)}>\times 
$$
$$
\pi({\psi}',\epsilon',\underline{\ell}',\underline{\eta}',(\rho,A,B,\zeta,\ell_{1}=0,\eta_{1}),(\rho,A+T-1,B+T+1,\zeta,\ell_{2}-1,\eta_{2})).\eqno(6)
$$
On calcule le module de Jacquet $Jac_{x\in {\cal C}(\zeta,A,B,T)}$ de l'induite de droite et on trouve par des calculs analogues \`a ceux d\'ej\`a faits:
$$
<\rho\vert\,\vert^{\zeta B}, \cdots, \rho\vert\,\vert^{-\zeta A}>\times Jac_{x\in {\cal C}(\zeta, A-1,B+1,T)}\pi(\underline{\psi}',(\rho,A,B,\zeta,\ell_{1},\eta_{1}),(\rho,A-1+T,B+1+T,\ell_{2}-1,\eta_{2})).
$$
Comme $\ell_{1}=0$, on peut remplacer $(\rho,A,B,\zeta,\ell_{1},\eta_{1})$ par $\cup_{C\in [B,A]}(\rho,C,C,\zeta,\eta_{1}(-1)^{C-B})$ et donc encore par $(\rho,B,B,\zeta,\eta_{1}),(\rho,B+1,B+1,\zeta,-\eta_{1}),(\rho,A,B+2,\zeta,\ell=0,\eta_{1})$. On remarque que $(\rho,A-1+T,B+1+T,\ell_{2}-1,\eta_{2})=(\rho,A+(T-1),B+2+(T-1),t_{2}-1,\eta_{2})$. Or ${\cal C}(\zeta,A,B+2,T-1)$ est pr\'ecis\'ement l'ensemble des $T-1$ premi\`eres lignes de ${\cal C}(\zeta,A-1,B+1,T)$. Par r\'ecurrence, on sait donc que pour que $Jac_{x\in {\cal C}(\zeta,A,B,T)}$ du premier membre de (4) soit non nul, il faut n\'ecessairement que $\ell_{2}=1$ et $\eta_{2}=\eta_{1}(-1)^{A-B}$. Ayant cela, on montre l'existence d'une inclusion de $\pi(\psi,\epsilon,\underline{\ell},\underline{\eta})$ dans
$$
\biggl( \begin{matrix}&\zeta (B+T) &\cdots&\zeta (B+2)&\cdots &-\zeta (A+1) &\cdots &-\zeta (A+T)\\
&\zeta (B+T+1) &\cdots&\zeta (B+3)&\cdots &-\zeta A\\
&\vdots &\vdots &\vdots\\
&\zeta (A+T-1)&\cdots&\zeta (A+1)&\cdots &-\zeta (B+2)
\end{matrix}\biggr)
\times$$
$$ \pi(\underline{\psi}',(\rho,B,B,\zeta,\eta_{1}),(\rho,B+1,B+1,\zeta,-\eta_{1})),
$$o\`u la matrice signifie en fait la repr\'esentation associ\'ee par Zelevinky \`a cet ensemble de segments. Pour tout $x\in [-\zeta (A+2), -\zeta (A+T)]$ l'induite $\rho\vert\,\vert^x\times  \pi(\underline{\psi}',(\rho,B,B,\zeta,\eta_{1}),(\rho,B+1,B+1,\zeta,-\eta_{1}))$ est irr\'eductible. On peut remplacer l'inclusion ci-dessus par une inclusion dans l'induite:
$$\biggl(
\begin{matrix}
&\zeta (B+T) &\cdots&\zeta (B+2)&\cdots &-\zeta (A+1) \\
&\zeta (B+T+1) &\cdots&\zeta (B+3)&\cdots &-\zeta A\\
&\vdots &\vdots &\vdots\\
&\zeta (A+T-1)&\cdots&\zeta (A+1)&\cdots &-\zeta (B+2)\\
&\zeta (A+T)&\cdots & \zeta(A+2)
\end{matrix}\biggr)
\times$$
$$ \pi(\underline{\psi}',(\rho,B,B,\zeta,\eta_{1}),(\rho,B+1,B+1,\zeta,-\eta_{1})).
$$
La difficult\'e vient de ce que $Jac_{x\in {\cal C}(\zeta,A,B,T)}$ de la repr\'esentation induite \'ecrite n'est pas nulle. Il faut factoriser diff\'eremment; on sort les 2 derni\`eres colonnes de la matrice, c'est-\`a-dire que l'on \'ecrit $\pi(\psi,\epsilon,\underline{\ell},\underline{\eta})$ comme sous-module irr\'eductible de l'induite:
$$
\biggl(
\begin{matrix}
&\zeta (B+T) &\cdots&\zeta (B+2)&\cdots &-\zeta (A-1) \\
&\zeta (B+T+1) &\cdots&\zeta (B+3)&\cdots &-\zeta (A-2)\\
&\vdots &\vdots &\vdots\\
&\zeta (A+T-1)&\cdots&\zeta (A+1)&\cdots &-\zeta (B)\\
&\zeta (A+T)&\cdots & \zeta(A+2)
\end{matrix}\biggr)\eqno(7)$$
$$
\times
<\rho\vert\,\vert^{-\zeta A}, \cdots, \rho\vert\,\vert^{-\zeta (B+1)}>\times <\rho\vert\,\vert^{-\zeta (A+1)}, \cdots, \rho\vert\,\vert^{-\zeta (B+2)}>\eqno(8)$$
$$ \times
 \pi(\underline{\psi}',(\rho,B,B,\zeta,\eta_{1}),(\rho,B+1,B+1,\zeta,-\eta_{1})).\eqno(8)'
$$
On peut remplacer  (8)$\times $ (8)' par un sous-quotient irr\'eductible, $\tau''$. Mais en utilisant \ref{jacquetsimple}, on sait que $Jac_{x}\tau''=0$ pour tout $x\in ]\zeta B, \zeta (A+T)]$. On calcule d'abord $Jac_{x\in {\cal C}(\zeta,B+1,A+1,T-1)}$ de l'induite (7)$\times \tau''$ et on trouve l'induite
$$
\biggl(
\begin{matrix}
&\zeta (B+1) &\cdots &-\zeta (A-1)\\
&\vdots &\vdots &\vdots\\
&\zeta A &\cdots &-\zeta B
\end{matrix}\biggr)\times \tau''
$$
et quand on applique $Jac_{\zeta (B+1), \cdots, \zeta (A+1)}$ \`a ce r\'esultat, on trouve 0. C'est ce que l'on voulait.\

On est donc maintenant ramen\'e au cas $\ell_{1}=\ell_{2}=0$ et il faut montrer que la condition $\eta_{1}=\eta_{2}(-1)^{A-B}$ est la condition n\'ecessaire et suffisante pour que le module de Jacquet soit non nul. 
Ici $\pi(\psi,\epsilon,\underline{\ell},\underline{\eta})$ n'est autre que $
\pi({\psi}',\epsilon',\underline{\ell}',\underline{\eta}',\cup_{C\in [B,A]} (\rho,C,C,\zeta,\eta_{1}(-1)^{C-B}), \cup_{C'\in [B,A]}(\rho, C'+T,C'+T,\zeta,\eta_{2}(-1)^{(C'-B)})).$
La non nullit\'e du module de Jacquet cherch\'e n\'ecessite d'abord la non nullit\'e du module de Jacquet relativement \`a la premi\`ere colonne du tableau ${\cal C}(\zeta,A,B,T)$ c'est \`a dire $Jac_{\zeta (B+T), \cdots, \zeta (B+1)}$ et a fortiori la non nullit\'e du  $Jac_{\zeta (B+T), \cdots \zeta (A+1)}$ de cette repr\'esentation. Or par d\'efinition cette non nullit\'e est exactement \'equivalente  \`a ce que la valeur du caract\`ere sur le bloc $(\rho,A,A,\zeta)$ soit la m\^eme que sur le bloc $(\rho,B+T,B+T,\zeta)$; c'est-\`a-dire pr\'ecis\'ement $\eta_{1}(-1)^{(A-B)}=\eta_{2}$. Quand cette \'egalit\'e est v\'erifi\'ee, on a une inclusion:
$$\pi({\psi}',\epsilon',\underline{\ell}',\underline{\eta}',\cup_{C\in [B,A]} (\rho,C,C,\zeta,\eta_{1}(-1)^{C-B}), \cup_{C'\in [B,A]}(\rho, C'+T,C'+T,\zeta,\eta_{2}(-1)^{(C'-B)}))\hookrightarrow
$$
$$
<\rho\vert\,\vert^{\zeta (B+T)}, \cdots, \rho\vert\,\vert^{-\zeta A}> \times \cdots \times
<\rho\vert\,\vert^{\zeta (B+T+\ell)}, \cdots, \rho\vert\,\vert^{-\zeta (A-\ell)}> \times \cdots \times
<\rho\vert\,\vert^{\zeta (A+T)}, \cdots, \rho\vert\,\vert^{-\zeta B}> \eqno(9)$$
$$\times
\pi({\psi}',\epsilon',\underline{\ell}',\underline{\eta}').
$$Il est maintenant facile de conclure. Allons un peu plus loin.
On peut remplacer l'induite dans (9) par un sous-quotient irr\'eductible mais ce sous-quotient, $\sigma$, doit v\'erifier que si $Jac_{x}\sigma \neq 0$ alors $x=\zeta (B+T)$. Cela veut dire que l'on peut remplacer (9) par l'unique repr\'esentation bas\'ee sur $\rho$ associ\'ee au tableau:
$$
\begin{matrix}&\zeta (B+T) &\cdots &-\zeta A\\
&\vdots &\vdots &\vdots\\
&\zeta (A+T) &\cdots &-\zeta B.
\end{matrix}
$$ 
En calculant $Jac_{x\in {\cal C}(\zeta,A,B,T)}$, on trouve l'inclusion, pour $\eta=\pm$:
$$Jac_{x\in {\cal C}(\zeta,A,B,\zeta)}
\pi({\psi}',\epsilon',\underline{\ell}',\underline{\eta}',(\rho,A,B,\zeta,0,\eta),(\rho,A+T,B+T,\zeta,0,\eta(-1)^{(A-B)}))\hookrightarrow
S(\rho,A,B,\zeta)\times \pi(\underline{\psi}').
$$
\bf Remarque: \rm
On vient donc de construire 2 sous-modules irr\'eductibles non nuls de l'induite de droite. L'objet du paragraphe \ref{induction} est de montrer que cette induite est semi-simple constitut\'ee de tous les sous-modules irr\'eductibles $Jac_{x\in {\cal C}(\zeta,A,B,T)}(\pi(\underline{\psi}',(\rho,A,B,\zeta,\ell,\eta),(\rho,A+T,B+T,\zeta,\ell,\eta(-1)^{(A-B)})))$ quand $\ell$ et $\eta$ varient. 
%%%%%%%

%%%%%%
\section{Induction\label{induction}}

On fixe  $\psi',\epsilon', \underline{\ell}',\underline{\eta}'$ permettant de d\'efinir, $\pi':=\pi(\psi',\epsilon',\underline{\ell}',\underline{\eta}')$ une repr\'esentation irr\'eductible ou $0$ d'un groupe de m\^eme type que $G$; et on suppose ici que $\pi'\neq 0$. On fixe aussi un quadruplet $(\rho,A,B,\zeta)$ ayant la propri\'et\'e de parit\'e. On suppose que pour tout $(\rho,A',B',\zeta') \in Jord(\psi')$ on a soit $A'<B$ soit $B'>>A$.  On posera aussi $\underline{\psi}'$ pour l'ensemble des donn\'ees $\psi',\epsilon',\underline{\ell}',\underline{\eta}'$.

\

\bf Proposition: \sl L'induite $S(\rho,A,B,\zeta)\times \pi'$ est semi-simple de longueur $A-B+2$. Ses constituants sont l'ensemble des repr\'esentations $$Jac_{x\in {\cal C}(\zeta,A,B,T)}(\pi(\underline{\psi}', (\rho,A,B,\zeta,\ell,\eta),(\rho,A+T,B+T,\zeta,\ell,\eta(-1)^{(A-B)})))$$ quand $\ell$ varie dans $[0,[(A-B+1)/2]]$ et $\eta=\pm$.\rm

\

\noindent Si $A-B$ est impair, on peut alors avoir $\ell=(A-B+1)/2$ et dans ce cas, il n'y a pas de $\eta$. Et on v\'erifie facilement que nous avons bien d\'efini $A-B+2$ param\`etres: en effet si $A-B$ est pair, il y a $(A-B)/2+1$ valeurs de $\ell$ et pour chaque valeurs de $\ell$, 2 valeurs de $\eta$, ce qui donne bien $A-B+2$ param\`etres. Tandis que si $A-B$ est impair, il y a $(A-B+1)/2+1$ valeurs de $\ell$ et pour toute valeurs de $\ell$ sauf exactement 1, il y a 2 valeurs de $\eta$, pour l'autre $\eta$ n'intervient pas; on a donc $2(A-B+1)/2+1=A-B+2$ param\`etres.

\

\bf Remarque: \sl On rappelle que $(\rho,A,B,\zeta)$ est associ\'e \`a une repr\'esentation irr\'eductible de $W_{F}\times SL(2,{\mathbb C})\times SL(2,{\mathbb C})$ qui est un produit tensoriel de la repr\'esentation $\rho$ vue comme repr\'esentation de $W_{F}$ et de 2 repr\'esentations de $SL(2,{\mathbb C})$ l'une de dimension, disons $a$, l'autre de dimension disons $b$. Dans cette interpr\'etation $A-B+2=inf(a,b)+1$.
\rm

\

On pose $\sigma:=S(\rho,A,B,\zeta)\times \pi'$. Voici une description de la preuve.

1- On montrera ci-dessous (cf. \ref{sousquotient}) que les $A-B+2$ repr\'esentations d\'ecrites dans l'\'enonc\'e sont des sous-quotients tous disctincts de $\sigma$.

2- On  montrera que l'ensemble des sous-modules irr\'eductibles de $\sigma$ est un inclus dans l'ensemble des repr\'esentations d\'ecrites dans l'\'enonc\'e.

Donc si on sait que $\sigma$ est semi-simple, la proposition est d\'emontr\'ee. Or on saurait que $\sigma$ est semi-simple si l'on sait que $\pi'$ est unitaire; ce serait une hypoth\`ese tout \`a fait raisonable, car $\pi'$ est associ\'ee \`a un param\`etre de restriction discr\`ete \`a la diagonale, que l'on sait transf\'erer ces repr\'esentations mises en paquet \`a un groupe lin\'eaire convenable d\`es que l'on sait le faire pour les s\'eries discr\`etes; il faudrait encore une comparaison des formules de traces globales et on r\'ealiserait ainsi $\pi'$ comme composante locale d'une forme automorphe de carr\'e int\'egrable; ce sont les id\'ees d'Arthur, cf. par exemple \cite{arthurnouveau}. Toutefois, il manque encore les lemmes fondamentaux et les lemmes fondamentaux pond\'er\'es, nous allons donc nous passer de cette hypoth\`ese. On conna\^{\i}t quand m\^eme l'unitarit\'e si $\pi'$ est une s\'erie discr\`ete.

On donne donc une d\'emonstration qui ne suppose pas l'unitarit\'e de $\pi'$. On remplace l'unitarit\'e par l'utilisation d'une dualit\'e. En effet, on montrera que $\sigma$ est autoduale pour une dualit\'e que l'on pr\'ecisera (essentiellement la contragr\'ediente); cela entra\^{\i}ne que tout sous-module irr\'eductible de $\sigma$ est isomorphe \`a un quotient irr\'eductible de $\sigma$; si un tel sous-module n'est pas facteur direct, il intervient donc avec multiplicit\'e au moins 2. Ainsi il est \'equivalent de d\'emontrer que $\sigma$ est semi-simple et que toute repr\'esentation de l'\'enonc\'e intervient avec multiplicit\'e 1 comme sous-quotient irr\'eductible de $\sigma$. C'est cette 2e propri\'et\'e que l'on ram\`ene \`a son analogue quand $\psi'$ est le param\`etre d'un paquet de s\'eries discr\`etes; dans ce cas, on  a l'unitarit\'e, donc la semi-simplicit\'e et donc la propri\'et\'e sur les multiplicit\'es.

\subsection{Description de certains sous-quotients\label{sousquotient}}
On note ${\cal I}$ l'ensemble des repr\'esentations $$\tau_{\ell,\eta}:=Jac_{x\in {\cal C}(\zeta,A,B,T)}(\pi(\underline{\psi}', (\rho,A,B,\zeta,\ell,\eta),(\rho,A+T,B+T,\zeta,\ell,\eta(-1)^{(A-B)}))),$$ quand $\ell$ et $\eta$ varient. On garde les notations $\sigma,\pi',\underline{\psi}'$ introduites ci-dessus.
\nl
\bf Lemme: \sl Toute repr\'esentation de ${\cal I}$ intervient comme sous-quotient de $\sigma$.\rm
\nl
On calcule $Jac_{\zeta B, \cdots, -\zeta A}\sigma$ et on trouve (comme semi-simplifi\'e) la somme de l'induite, o\`u le tableau \'ecrit doit \^etre vu comme la repr\'esentation associ\'ee:
$$
\biggl(
\begin{matrix}
&\zeta (B+1) & \cdots &-\zeta (A-1)\\
&\vdots &\cdots &\vdots\\
&\zeta A &\cdots &-\zeta B
\end{matrix}
\biggr) \times \pi'
$$
$$
\oplus 
\biggl(
\begin{matrix}
&\zeta B & \cdots &-\zeta A\\
&\vdots &\cdots &\vdots\\
&\zeta (A-1) &\cdots &-\zeta (B+1)
\end{matrix}
\biggr) \times \pi'.
$$
Pour avoir ce r\'esultat, il faut couper l'intervalle $[\zeta B,-\zeta A]$ en trois sous-ensembles ${\cal E}_{i}$, pour $i=1,2,3$ tels que (en utilisant l'autodualit\'e de $S(\rho,A,B,\zeta)$)
$$
Jac_{x\in {\cal E}_{1}}S(\rho,A,B,\zeta)\neq 0; \qquad Jac_{x\in -^t{\cal E}_{2}}S(\rho,A,B,\zeta)\neq 0; \qquad Jac_{x\in {\cal E}_{3}}\pi'\neq 0.
$$
On sait gr\^ace \`a \ref{jacquetsimple} que $\zeta B$ n'est pas dans ${\cal E}_{3}$; d'apr\`es les formules standard si ${\cal E}_{i}\neq \emptyset$ que $i=1$ ou $2$ alors $\zeta B$ est le premier \'el\'ement de ${\cal E}_{i}$. Ainsi soit ${\cal E}_{1}$ soit ${\cal E}_{2}$ est vide et l'autre ensemble contient $\zeta B$. On sait aussi que $-\zeta A\notin {\cal E}_{3}$ (cf \cite{paquetsdiscrets}, 3.4) c'est donc un \'el\'ement du ${\cal E}_{i}$ pour $i=1$ ou $2$ qui est non vide. Mais alors cet ensemble contient tout le segment $[\zeta B, -\zeta A]$ par les formules standard; l'\'eventualit\'e ${\cal E}_{1}\neq \emptyset$ donne le premier terme et ${\cal E}_{2}\neq \emptyset$ donne le deuxi\`eme terme. La m\^eme d\'emonstration calcule
$$
Jac_{\zeta B, \cdots, -\zeta A}Jac_{\zeta B, \cdots, -\zeta A}\sigma \eqno(1)
$$
et on trouve exactement 2 copies de $\sigma':=S(\rho,A-1,B+1,\zeta)\times \pi'$. On conna\^{\i}t, par r\'ecurrence  sur $A-B$, la structure de cette repr\'esentation.

Puisque (1) est non nul, il existe un sous-quotient $\tau'$ de $\sigma$ tel que 
$$
Jac_{\zeta B, \cdots, -\zeta A}Jac_{\zeta B, \cdots, -\zeta A} \tau'\neq 0.
$$
La r\'eciprocit\'e de Frobenius entra\^{\i}ne alors l'existence d'une repr\'esentation irr\'eductible $\tau''$ et d'une inclusion de $\tau'$ dans l'induite
$$
 \rho\vert\,\vert^{\zeta B}\times \cdots \times \rho\vert\,\vert^{-\zeta A}\times  \rho\vert\,\vert^{\zeta B}\times \cdots \times \rho\vert\,\vert^{-\zeta A}\times \tau''.
 $$
 Et par exactitude du foncteur de Jacquet $\tau''$ est un sous-quotient irr\'eductible de $S(\rho,A-1,B+1,\zeta)\times \pi'$ ainsi $\tau''$ correspond \`a des donn\'ees $\ell'$ et $\eta$ quand on remplace $(A,B)$ par $(A-1,B+1)$. On v\'erifie que l'inclusion ci-dessus se factorise par le sous-module
 $$
 <\rho\vert\,\vert^{\zeta B}, \cdots, \rho\vert\,\vert^{-\zeta A}>\times  <\rho\vert\,\vert^{\zeta B}, \cdots, \rho\vert\,\vert^{-\zeta A}>\times \tau'',\eqno(2)
 $$
 simplement car pour tout $x\in ]\zeta B, -\zeta A]$, $Jac_{x, \cdots, -\zeta A
}\sigma=0$ et donc aussi pour $\tau'$.
On a vu que (2) a un unique sous-module irr\'eductible qui est l'\'el\'ement de ${\cal I}$ correspondant \`a $\ell'+1$ et $\eta$. Mais quand on calcule 
$$
Jac_{\zeta B, \cdots, -\zeta A}Jac_{\zeta B, \cdots, -\zeta A} \tau'
$$
on trouve exactement $\tau''$ avec multiplicit\'e 2. Pour \'epuiser tous les sous-quotients, $\tau$ de $\sigma$ v\'erifiant $Jac_{\zeta B, \cdots, -\zeta A}Jac_{\zeta B, \cdots, -\zeta A} \tau\neq 0$, il faut donc intervenir toutes les valeurs de $\ell'\in [0, [(A-B+1)/2]-1]$ et de $\eta'$ possibles. On obtient donc tous les \'el\'ements de ${\cal I}$ tels que $\ell \geq 1$. 

On a vu dans la remarque \`a la fin de \ref{nullitenonnullite} que les 2 repr\'esentations de ${\cal I}$ pour lesquelles $\ell=0$ sont des sous-modules irr\'eductibles de $\sigma$ et cela termine la preuve du lemme.

\subsection{Description de certains sous-modules\label{sousmodule}}
On garde la notation $\tau_{\ell,\eta}$ introduite dans \ref{sousquotient}.
On note $S_{-}$ l'unique sous-module irr\'eductible associ\'e \`a l'ensemble ordonn\'e (de haut en bas) de segments bas\'es sur $\rho$
$$
\begin{matrix}&\zeta B &\cdots &-\zeta A\\
&\vdots &\vdots &\vdots\\
& \zeta (A-1) &\cdots &-\zeta B
\end{matrix}.
$$
Et on a une inclusion de $S(\rho,A,B,\zeta)$ dans $S_{-}\times <\rho\vert\,\vert^{\zeta A}, \cdots, \rho\vert\,\vert^{-\zeta B}>$. On sait d\'ecomposer la repr\'esentation $<\rho\vert\,\vert^{\zeta A}, \cdots, \rho\vert\,\vert^{-\zeta B}>\times \pi'$; en effet  cette repr\'esentation est de longueur 3 exactement, contenant les 2 sous-modules irr\'eductibles, $\pi(\underline{\psi}', (\rho,A,A,\zeta,\epsilon),(\rho, B,B,\zeta,\epsilon))$ pour $\epsilon=\pm 1$ et l'unique sous-module irr\'eductible, $\pi_{u}$, de l'induite
$
<\rho\vert\,\vert^{\zeta B}, \cdots, \rho\vert\,\vert^{-\zeta A}>\times \pi'.
$. En effet,
un tel r\'esultat est connu si $A=B+1$, dans ce cas, $\pi_{u}=\pi(\underline{\psi}',(\rho,B+1,B,\zeta,+))$ et on a d\'ecrit la d\'ecomposition de cet \'el\'ement dans le groupe de Grothendieck dans \cite{paquetsdiscrets} et c'est exactement ce que nous avons anonc\'e. Si $A>B+1$, on d\'emontre que pour tout sous-quotient irr\'eductible de l'induite, $\tau$, $Jac_{\zeta A}\tau\neq 0$; je ne vois pas d'autre d\'emonstration que celle donn\'ee en \cite{paquetsdiscrets}, c'est-\`a-dire se ramener au cas o\`u $\pi'$ est cuspidal. Le cas o\`u $\pi'$ est cuspidal est facile car $\rho\vert\,\vert^{\zeta A}\times \pi'$ est irr\'eductible; ceci permet alors de remplacer $A$ par $A-1$ et de conclure par r\'ecurrence.

\noindent
On note $\pi_{d}$
la somme des 2 sous-modules irr\'eductibles d\'efinis ci-dessus et $\pi_{u}$ le quotient irr\'eductible d\'ej\`a introduit. On a donc une suite exacte:
$$
0\rightarrow S_{-}\times \pi_{d} \rightarrow S_{-}\times <\rho\vert\,\vert^{\zeta A}, \cdots, \rho\vert\,\vert^{-\zeta B}>\times \pi' \rightarrow S_{-}\times \pi_{u}\rightarrow 0.
$$
D'autre part $\sigma$ est une sous-repr\'esentation de la repr\'esentation du milieu; on note $\sigma_{d}$ l'intersection de $\sigma$ avec $S_{-}\times \pi_{d}$ et $\sigma_{u}$ l'image de $\sigma$ dans $S_{-}\times \pi_{u}$. On a donc la suite exacte:
$$
0\rightarrow \sigma_{d}\rightarrow \sigma \rightarrow \sigma_{u}\rightarrow 0.
$$
\bf Lemme: \sl Les sous-modules irr\'eductibles de $\sigma_{d}$ sont exactement les repr\'esentations $\tau_{\ell=0,\eta}$ pour $\eta=\pm 1$ et les sous-modules irr\'eductibles de $\sigma_{u}$ sont  de la forme $\tau_{\ell,\eta}$ pour $\ell>0$ chacune apparaissant au plus une fois.
\rm
\nl
L'assertion concernant les sous-modules de $\sigma_{u}$ est la plus simple. En effet un sous-module
 irr\'eductibles de $\sigma_{u}$ est aussi un sous-module irr\'eductible de $S_{-}\times \pi_{u}$ et donc de 
$$
S_{-}\times <\rho\vert\,\vert^{\zeta B}, \cdots, \rho\vert\,\vert^{-\zeta A}>\times \pi'.
$$
On utilise le fait que la repr\'esentation $S_{-}\times <\rho\vert\,\vert^{\zeta B}, \cdots, \rho\vert\,\vert^{-\zeta A}>$ est irr\'eductible (cf. \ref{undeuxiemeresultatdirreductibilite}) et on peut donc \'echanger les facteurs. Ainsi tout sous-module irr\'eductible de $S_{-}\times \pi_{u}$ est certainement un sous-module irr\'eductible de l'induite
$$
<\rho\vert\,\vert^{\zeta B}, \cdots, \rho\vert\,\vert^{-\zeta A}> \times <\rho\vert\,\vert^{\zeta B}, \cdots, \rho\vert\,\vert^{-\zeta A}> \times S(\rho,A-1,B+1,\zeta)\times \times \pi'.
$$
On a d\'ej\`a montr\'e que ces sous-modules irr\'eductibles sont exactement de la forme $\tau_{\ell,\eta}$ avec $\ell>0$. Cela d\'emontre, a fortiori, la 2e partie du lemme.
\nl
Montrons maintenant l'assertion concernant $\sigma_{d}$;  il suffit  de d\'emontrer que cette repr\'esenta\-tion a au plus 2 sous-modules irr\'eductibles, ils seront alors n\'ecessairement les repr\'esentations $\tau_{\ell=0,\eta}$ pour $\eta=\pm$ dont on a vu qu'elles n'\'etaient pas sous-modules de $\sigma_{u}$. On fixe $\epsilon=\pm$ et on note $\sigma_{d,\epsilon}$ l'intersection de $\sigma_{d}$ avec $S_{-}\times \pi(\underline{\psi'},(\rho,A,A,\zeta,\epsilon),(\rho, B,B,\zeta,\epsilon))$. 

Il faut d'abord consid\'erer le cas o\`u $A=B+1$. Ici on a simplement l'inclusion
$$
\sigma_{d,\epsilon}\hookrightarrow <\rho\vert\,\vert^{\zeta B}, \cdots, \rho\vert\,\vert^{-\zeta (B+1)}>\times \pi(\underline{\psi'},(\rho,B+1,B+1,\zeta,\epsilon),(\rho,B,B,\zeta,\epsilon)).
$$
On calcule par les formules standard, $Jac_{\zeta B, \cdots, -\zeta (B+1)}$ de l'induite de droite. On v\'erifie que pour tout $x$ tel que $\vert x\vert \leq (B+1)$, on a $$Jac_{x, \cdots, -\zeta (B+1)}\pi(\underline{\psi'},(\rho,B+1,B+1,\zeta,\epsilon),(\rho,B,B,\zeta,\epsilon))=0;$$
Pour cela on revient \`a la d\'efinition comme sous-module irr\'eductible de l'induite
$$
<\rho\vert\,\vert^{\zeta (B+1)}, \cdots, \rho\vert\,\vert^{-\zeta B}>\times \pi';
$$
Ceci entra\^{\i}ne qu'il faudrait une non nullit\'e analogue (\'eventuellement en changeant $x$ en $x'$ v\'erifiant encore $\vert x'\vert\leq (B+1)$) pour $\pi'$; ceci est exclu par \cite{paquetsdiscrets}, 3.4. Ceci entra\^{\i}ne que $$
Jac_{\zeta B, \cdots, -\zeta (B+1)}\biggl(<\rho\vert\,\vert^{\zeta B}, \cdots, \rho\vert\,\vert^{-\zeta (B+1)}>\times \pi(\underline{\psi'},(\rho,B+1,B+1,\zeta,\epsilon),(\rho,B,B,\zeta,\epsilon))\biggr)
$$
$$
=\pi(\underline{\psi'},(\rho,B+1,B+1,\zeta,\epsilon),(\rho,B,B,\zeta,\epsilon)).
$$
Cette induite a donc un unique sous-module irr\'eductible et il en est donc de m\^eme de $\sigma_{d,\epsilon}$. Cela termine la preuve dans ce cas.

On suppose donc que $A>B+1$.
On a l'inclusion
$
\sigma_{d,\epsilon}\hookrightarrow$
$$ <\rho\vert\,\vert^{\zeta B}, \cdots, \rho\vert\,\vert^{-\zeta A}>\times
S(\rho,A-1,B+1,\zeta)\times  \pi(\underline{\psi'},(\rho,A,A,\zeta,\epsilon),(\rho, B,B,\zeta,\epsilon)).\eqno(1)$$
Puisque $(A,B)$ est remplac\'e par $(A-1,B+1)$ en (1), on sait d\'ecomposer (1). La repr\'esentation (1) est une repr\'esentation semi-simple dont on note $\tau'_{\ell',\eta'}$ les diff\'erents composants de fa\c{c}on analogue aux $\tau_{\ell,\eta}$. On va montrer que l'inclusionde $\sigma_{d,\epsilon}$ est \`a valeurs dans
$
<\rho\vert\,\vert^{\zeta B}, \cdots, \rho\vert\,\vert^{-\zeta A}>\times \tau'_{\ell'=0,-\epsilon}.
$
Supposons qu'il n'en soit pas ainsi et fixons $\ell',\eta'$ tel qu'il existe un morphisme non nul de $\sigma_{d,\epsilon}$ dans l'induite:
$$
<\rho\vert\,\vert^{\zeta B}, \cdots, \rho\vert\,\vert^{-\zeta A}>\times \tau'_{\ell',\eta'}.
$$ 
Si $\ell'>0$, alors $Jac_{\zeta (B+1),\zeta (B+1)}\tau'_{\ell',\eta'}\neq 0$, c'est-\`a-dire qu'il existe une repr\'esentation irr\'eductible $\tau''$ et une inclusion de $\tau'_{\ell',\eta'}$ dans l'induite $\rho\vert\,\vert^{\zeta(B+1)}\times \rho\vert\,\vert^{\zeta (B+1)}\times \tau''$. Comme tout sous-quotient irr\'eductible, $\lambda$, de l'induite 
$
<\rho\vert\,\vert^{\zeta B}, \cdots, \rho\vert\,\vert^{-\zeta A}>\times \rho\vert\,\vert^{\zeta(B+1)}\times \rho\vert\,\vert^{\zeta (B+1)}
$
v\'erifie $Jac_{\zeta (B+1)}\lambda\neq 0$, on aura un sous-quotient de $\sigma_{d,\epsilon}$ dont le $Jac_{\zeta (B+1)}$ n'est pas nul; ceci est contradictoire avec le fait que $Jac_{\zeta (B+1)}\sigma=0$. Supposons maitenant que $\eta'\neq -\epsilon$; en partant des d\'efinitions, on v\'erifie que
$$
Jac_{\zeta A, \cdots, -\zeta (A-1)}\tau'_{\ell'=0,\epsilon}\neq 0.
$$
Comme $A>(B+1)$, l'induite $<\rho\vert\,\vert^{\zeta B}, \cdots, \rho\vert\,\vert^{-\zeta A}>\times\rho\vert\,\vert^{\zeta A}$ est irr\'eductible et on trouverait un morphisme non nul de $\sigma_{d,\epsilon}$ dans une induite de la forme $\rho\vert\,\vert^{\zeta A}\times \tau'''$ o\`u $\tau'''$ est quelconque. Par fonctorialit\'e cela force encore $Jac_{\zeta A}\sigma_{d,\epsilon}\neq 0$ ce qui est contradictoire avec le fait que $Jac_{\zeta A}\sigma=0$.
On vient donc de d\'emontrer l'inclusion
$
\sigma_{d,\epsilon}\hookrightarrow <\rho\vert\,\vert^{\zeta B}, \cdots, \rho\vert\,\vert^{-\zeta A}>\times \tau'_{0,-\epsilon}.
$
Ensuite on conclut comme dans le cas o\`u $A=B+1$ en montrant que $
Jac_{\zeta B, \cdots, -\zeta A}\biggl(<\rho\vert\,\vert^{\zeta B}, \cdots, \rho\vert\,\vert^{-\zeta A}>\times \tau'_{\ell',\eta'}\biggr)=\tau'_{0,-\epsilon}.
$
Cela termine la preuve du lemme.

%%%%%%%%

 \subsection{Autodualit\'e\label{autodualite}}
 Pour un groupe classique, la repr\'esentation duale d'une repr\'esentation irr\'eductible peut se r\'ealiser dans le m\^eme espace en  tordant \'eventuellement par un automorphisme venant du groupe des similitudes, cet \'el\'ement ne d\'epend que du groupe et non de la repr\'esentation. Soit $\pi_{0}$ une repr\'esentation irr\'eductible d'un groupe classique et $\delta$ une repr\'esentation autoduale d'un produit de groupes lin\'eaires. On consid\`ere l'induite $\delta\times \pi_{0}$ et sa duale; cette duale se r\'ealise dans le m\^eme espace en tordant \'eventuellement par un automorphisme comme ci-dessus qui ne d\'epend que du groupe de la repr\'esenta\-tion $\pi_{0}$. On a alors l'\'equivalence entre les 2 assertions suivantes:
 
 (1) $\delta\times \pi_{0}$ est semi-simple sans multiplicit\'e

 (2) tout sous-module irr\'eductible de $\delta\times \pi_{0}$ intervient avec multiplicit\'e 1 comme sous-quotient de $\delta\times \pi_{0}$.

 \

 \noindent En effet (1) entra\^{\i}ne \'evidemment (2) et la r\'eciproque se prouve ainsi. Notons $\gamma$ la similitude qui sert \`a entrelacer $\pi_{0}$ et sa duale et on voit cet \'el\'ement dans le ''gros'' groupe de l'induite et cet \'el\'ement normalise le parabolique qui sert \`a induire. Soit $\tau$ un sous-module irr\'eductible de $\delta\times \pi_{0}$; alors $^\gamma\tau$ est un sous-module irr\'eductible de $\delta\times \, ^\gamma\pi_{0}$ mais par dualit\'e $\sigma^*$ est un quotient de $\delta\times \pi_{0}^*$, c'est -\`a-dire que $^\gamma \tau$ est un quotient de $\delta\times \, ^\gamma\pi_{0}$. Ainsi $^\gamma \tau$ est isomorphe \`a un sous-module et \`a un quotient de cette induite et par dualit\'e $\tau$ est aussi isomorphe \`a un sous-module et \`a un quotient de $\delta\times \pi_{0}$. Par multiplicit\'e 1 comme sous-quotient, $\tau$ est donc facteur direct de $\delta\times \pi_{0}$. Ainsi la somme des sous-modules irr\'eductibles de $\delta\times \pi_{0}$ est facteur direct de $\delta\times \pi_{0}$ et co\"{\i}ncide donc avec toute l'induite. D'o\`u l'assertion.

 %%%%%%%%%%
 \subsection{Fin de la preuve de la semi-simplicit\'e\label{findelapreuve}}
 En tenant comte de \ref{autodualite}, il nous reste \`a d\'emontrer que les repr\'esentations $\tau_{\ell,\eta}$ d\'efinies en \ref{sousquotient}, interviennent dans $\sigma$ avec multiplicit\'e 1 en tant que sous-quotient irr\'eductible. On conna\^{\i}t le r\'esultat si $\pi'$ est unitaire donc en particulier si $\pi'$ est une s\'erie discr\`ete. On suppose d'abord que $\psi'$ est un morphisme \'el\'ementaire.
 
Ce cas se ram\`ene \`a celui des s\'eries discr\`etes en appliquant une involution d\'ecrite dans \cite{elementaire}. On applique cette involution non pas \`a $\pi(\psi')$ mais \`a l'induite $\sigma$. Montrons que cette involution commute \`a l'induction par $S(\rho,A,B,\zeta)$, ici il faut \'evidemment l'hypoth\`ese que pour tout $C\in [B,A]$, $(\rho,2C+1)\notin Jord(\psi'\circ \Delta)$ (sinon ce serait faux). Pr\'ecis\'ement, on a d\'efini pour tout entier $\alpha$ une involution $inv_{<\alpha}$ dans le groupe le groupe de Grothendieck; pour cela on note pour tout produit de groupes lin\'eaires $proj_{<\alpha}$ la projection sur l'ensemble des repr\'esentations dont le support cuspidal ne fait intervenir que des \'el\'ements $\rho\vert\,\vert^x$ avec $\vert x\vert <(\alpha-1)/2$ et  on \'etend cette application \`a tout L\'evi du groupe classique $G$ en la faisant porter sur le produit des groupes lin\'eaires intervenant dans le Levi. On a alors pos\'e, en imitant les formules de \cite{aubert} et \cite{SS}:
$$
inv_{<\alpha}\pi:= \sum_{P}(-1)^{(rg(G)-rg(P))} ind_{P}^G proj_{<\alpha}(res_{P}^G)\pi,
$$
o\`u la somme parcourt l'ensemble des paraboliques standard. On d\'efinit de la m\^eme fa\c{c}on $inv_{\leq \alpha}$ en rempla\c{c}ant les in\'egalit\'es strictes par des in\'egalit\'es larges. Par les d\'efinitions de \cite{elementaire}, pour toute application $\epsilon'$ de $Jord(\psi')$ dans $\pm 1$:
$$
\pi(\psi',\epsilon')=\prod_{\alpha\in Jord(\psi'\circ \Delta); \zeta(\alpha)=-}inv_{< \alpha}\circ inv_{\leq\alpha}\pi(\psi'\circ \Delta, \epsilon'),
$$
o\`u on a identifi\'e $Jord(\psi') $ et $Jord(\psi'\circ \Delta)$ pour d\'efinir le 2e $\epsilon'$, en oubliant la fonction $\zeta$. Or si $\alpha < 2B+1$, il r\'esulte pratiquement des d\'efinitions que $inv_{<\alpha}$ et $inv_{\leq \alpha}$ commute \`a l'induction par $S(\rho,A,B,\zeta)$ tout simplement parce que $Jac_{x}S(\rho,A,B,\zeta)\neq 0$ entra\^{\i}ne que $x=\zeta B$ et que $S(\rho,A,B,\zeta)$ est autoduale. Par contre,  si $\alpha>2A+1$ alors tout $\rho\vert\,\vert^x$ dans le support cuspidal de $S(\rho,A,B,\zeta)$ v\'erifie $\vert x\vert <(\alpha-1)/2$ et on a montr\'e en \cite{elementaire} (c'est en fait facile avec la d\'efinition) que 
$$
inv_{<\alpha}(S(\rho,A,B,\zeta)\times \tau)=inv(S(\rho,A,B,\zeta))\times inv_{<\alpha}\tau,
$$
pour toute repr\'esentation $\tau$, o\`u $inv$ est l'involution d'Iwahori-Matsumoto. Ici $inv(S(\rho,A,B,\zeta))=S(\rho,A,B,-\zeta)$. On a la m\^eme formule avec $inv_{\leq \alpha}$ et le produit $inv_{<\alpha}\circ inv_{\leq \alpha}$ commute donc avec l'induction par $S(\rho,A,B,\zeta)$.  D'o\`u le r\'esultat de commutation annonc\'e, puisque pour tout $\alpha\in Jord(\psi')$ on a soit $\alpha>2A+1$ soit $\alpha<2B+1$.
%%%%

Nous n'avons pas d\'emontr\'e que $inv_{\leq \alpha}$ envoie une repr\'esentation irr\'eductible sur une repr\'esentation irr\'eductible et c'est  d'ailleurs faux, dans notre situation, si on prenait $\alpha \in [B,A[$. Toutefois on a d\'emontr\'e cette conservation de l'irr\'eductibilit\'e pour les repr\'esentations  qui sont sous-quotients de $S(\rho,A,B,\zeta)\times \pi(\psi,\epsilon)$ d\`es que $\alpha \geq  A$ (cf. assertion de la preuve de \cite{elementaire}, 4.1) \`a un signe pr\`es. On peut donc conclure si pour tout $(\rho,A',B',\zeta')\in Jord(\psi)$ avec n\'ecessairement $A'=B'$, on a aussi $\zeta'=+$ d\`es que $B'<B$. On r\`egle donc d'abord par r\'ecurrence le cas o\`u ceci n'est pas r\'ealis\'e; on a \`a examiner 2 cas d'apr\`es les d\'efinitions de \cite{elementaire}.

1e cas, il existe  $(\rho,A',B',\zeta')$ avec $B'<B$ tel que soit $(\rho,A'-1,B'-1,\zeta')\notin Jord(\psi)$; alors $Jac_{\zeta' B'}\pi(\psi,\epsilon)=\pi(\psi',\epsilon)$ o\`u $Jord(\psi')$ se d\'eduit de $Jord(\psi)$ en rempla\c{c}ant simplement $(\rho,A',B',\zeta')$ par $(\rho,A'-1,B'-1,\zeta')$ (qui dispara\^{\i}t si $B'=1$)

2e cas, il existe $(\rho,A',B',\zeta')$ avec $B'<B$ tel que $(\rho,A'-1,B'-1,\zeta')\in Jord(\psi)$ et $$Jac_{\zeta' B', \cdots, -\zeta'(B'-1)}\pi(\psi,\epsilon)=\pi(\psi',\epsilon),$$ o\`u $Jord(\psi')$ se d\'eduit de $Jord(\psi)$ en enlevant les 2 blocs $(\rho,A',B',\zeta')$ et $(\rho,A'-1,B'-1,\zeta')$.

On v\'erifie que les $Jac$ consid\'er\'es, appliqu\'es \`a $S(\rho,A,B,\zeta)\times \pi(\psi,\epsilon)$ n'agissent que sur $\pi(\psi,\epsilon)$ et  n'annulent aucune des repr\'esentations $\tau_{\ell,\epsilon}$ qui nous int\'eressent. Ainsi si l'une de ces repr\'esentations intervenaient avec une multiplicit\'e plus grande que 1 dans l'induite avec $\psi$ cela entra\^{\i}nerait que l'induite avec $\psi'$ aurait de la multiplicit\'e elle aussi. On conclut alors facilement. On se ram\`ene ainsi au cas, o\`u il n'y a \`a consid\'erer que des $inv_{<\alpha}$ pour $\alpha >A$ et $inv_{\leq \alpha}$ pour $\alpha\geq A$ o\`u on sait que l'application consid\'er\'ee conserve l'irr\'eductibilit\'e des sous-quotients de notre induite.

%%%%
On ram\`ene maintenant le cas g\'en\'eral au cas des morphismes \'el\'ementaires. Si $\underline{\ell'}$ qui fait partie des donn\'ees d\'efinissant $\pi'$ est identiquement nul, $\pi'$ est une repr\'esentation dans un paquet associ\'e \`a un morphisme \'el\'ementaire et on a donc l'assertion. On suppose donc que $\underline{\ell}'$ n'est pas identiquement nul; on a donc un \'el\'ement $(\rho',A',B',\zeta')\in Jord(\psi')$ qui permet de construire une inclusion
$$
\pi'\hookrightarrow <\rho'\vert\,\vert^{\zeta' B'}, \cdots, \rho'\vert\,\vert^{-\zeta' A'}>\times \pi''
$$
o\`u $\pi''$ est associ\'e \`a des donn\'ees analogues \`a celles d\'efinissant $\pi'$; on change $(\rho',A',B',\zeta')$ en $(\rho',A'-1,B'+1,\zeta')$ et $\underline{\ell}'$ change uniquement sur cet \'el\'ement. Par hypoth\`ese, on sait que soit $B>A'$ soit $B'>A$. Pour \'eviter des probl\`emes de notations, on suppose que $\rho'=\rho$ (ce qui est de toute fa\c{c}on le cas le plus difficile, a priori). On v\'erifie sur la d\'efinition que 
$
Jac_{\zeta' B', \cdots, -\zeta' A'}\tau_{\ell,\eta}\neq 0
$
contient l'analogue,$\tau'_{\ell',\eta}$, de $\tau_{\ell,\eta}$ quand on remplace $\pi'$ par $\pi''$ (en fait il y a \'egalit\'e mais nous n'avons pas besoin d'un r\'esultat aussi pr\'ecis). De plus
$$
Jac_{\zeta' B', \cdots, -\zeta' A'}\sigma= S(\rho,A,B,\zeta) \times \pi''.
$$Par exactitude du foncteur de Jacquet, $\tau_{\ell,\eta}'$ intervient avec multiplicit\'e au moins 2 comme sous-quotient de $S(\rho,A,B,\zeta)\times \pi''$; ce qui est exclu. Cela termine la preuve.

  %%%%%%%%
 \subsection{Induction suite \label{inductionsuite}}
 On fixe des donn\'ees $\underline{\psi}':=(\psi',\epsilon',\underline{\ell}',\underline{\eta}')$ permettant de construire une repr\'esentation irr\'eductible $\pi(\underline{\psi}')$. On fixe aussi $(\rho,A,B,\zeta)$ ayant la condition de parit\'e et un entier $a$. Le but de ce paragraphe est de d\'ecrire la repr\'esentation induite:
 $$
 \sigma_{a}:=S(\rho,A,B,\zeta) \times \cdots \times S(\rho,A,B,\zeta)\times \pi(\underline{\psi}'),
 $$o\`u il y a $a$ copies de $S(\rho,A,B,\zeta)$.

 On note $\psi_{a}$ le morphisme qui se d\'eduit de $\psi'$ en ajoutant $2a$ copies de $(\rho,A,B,\zeta)$ \`a l'ensemble $Jord(\psi')$ et on fixe un morphisme $\tilde{\psi}_{a}$ dominant $\psi_{a}$. En particulier il existe des entiers $T_{i}$ pour $i\in [1,2a]$ tels que chaque $(\rho,A+T_{i},B+T_{i},\zeta) \in Jord(\tilde{\psi}_{a})$ corresponde \`a une des copies de $(\rho,A,B,\zeta)\in Jord(\psi)$ ajout\'ee. Pour fixer les notations, on suppose que $T_{1}<T_{2}<\cdots <T_{2a}$. Toutefois, il continue d'y avoir une ambiguit\'e si $(\rho,A,B,\zeta)\in Jord(\psi')$ pour la lever, on suppose que dans ce cas l'\'el\'ement de $Jord(\tilde{\psi}_{a})$ qui correspond \`a ce bloc est de la forme $(\rho,A+T_{0},B+T_{0},\zeta)$ avec $T_{0}<T_{1}$. Etant donn\'e le th\'eor\`eme ci-dessus, il n'est pas utile de consid\'erer le cas o\`u $(\rho,A,B,\zeta)$ appara\^{\i}t dans $Jord(\psi')$ avec une multiplicit\'e plus grande que 1.

 Pour $\ell\in [0,[(A-B+1)/2]]$ et pour $\eta=\pm$ qui n'est d\'efini que si $\ell<(A-B+1)/2$, on prolonge $\underline{\ell}',\underline{\eta}'$ \`a $Jord(\tilde{\psi}_{a})$ en  $\underline{\ell}$, $\underline{\eta}$, en posant $\underline{\ell}(\rho,A+T_{i},B+T_{i},\zeta)=\ell$ pour tout $i\in [1,2a]$ et $\underline{\eta}(\rho,A+T_{i},B+T_{i},\zeta)=\eta (-1)^{(i-1)(A-B)}$. Si $(\rho,A,B,\zeta)\in Jord(\psi')$ on suppose que $\ell=\underline{\ell}'(\rho,A,B,\zeta)$ et que $\eta=\underline{\eta}(\rho,A,B,\zeta)$. 
 
 On pose alors $\pi(\psi_{a},\epsilon,\underline{t},\underline{\eta}):=Jac_{{\cal C}}\pi(\tilde{\psi}_{a},\epsilon,\underline{t},\underline{\eta})$, o\`u ${\cal C}$ est un ensemble convenable comme en \ref{definitiongenerale}.
 
 \
 
 \bf Th\'eor\`eme: \sl La repr\'esentation $\sigma_{a}$ est semi-simple de longueur inf\'erieure ou \'egale \`a $A-B+2$; elle est la somme directe des repr\'esentations $\pi(\psi_{a},\epsilon,\underline{t},\underline{\eta})$ d\'efinies ci-dessus. Elle est en particulier irr\'eductible si $(\rho,A,B,\zeta)\in Jord(\psi')$.\rm

\

Pour faire la d\'emonstration, on  travaille dans le groupe de Grothendieck, on d\'ecrit donc le semi-simplifi\'e de $\sigma_{a}$. En particulier on va montrer qu'il est sans multiplicit\'e car constitu\'e exactement des repr\'esentations d\'ecrites par l'\'enonc\'e (cf. \ref{irreductibilite}) et cela suffit pour savoir que $\sigma_{a}$ est semi-simple (cf \ref{autodualite}).

\

On note $\tilde{\psi}'$ le morphisme qui se d\'eduit de $\tilde{\psi}_{a}$ en enlevant les blocs $(\rho,A+T_{i},B+T_{i},\zeta)$ pour $i\in [1,2a]$; il domine $\psi'$ et $\underline{\ell}',\underline{\eta}'$ sont bien d\'efinis sur $Jord(\tilde{\psi}')$. On reprend la notation $\underline{\tilde{\psi}'}$ pour l'ensemble des donn\'ees $\tilde{\psi}',\epsilon',\underline{\ell}',\underline{\eta}'$. 
Avant l'\'enonc\'e, on a d\'efini certains prolongements de $\underline{\ell}'$ et $\underline{\eta}'$ \`a $Jord(\tilde{\psi}_{a})$ parce que ce sont les seuls qui vont intervenir. Mais pour la preuve on a besoin de consid\'erer tous les prolongements possibles, $\underline{\tilde{\ell}}$ et $\underline{\tilde{\eta}}$. Gr\^ace \`a \ref{induction} dont les hypoth\`eses sont remplies:
$$
\times_{j\in [1,a]}S(\rho,A+T_{2j-1},B+T_{2j-1},\zeta) \times \pi(\underline{\tilde{\psi}'})=$$
$$\oplus_{\underline{\tilde{\ell}},\underline{\tilde{\eta}}}Jac_{x\in \cup_{j\in [1,a]} {\cal C}(\zeta, A+T_{2j-1},B+T_{2j-1},T_{2j}-T_{2j-1})}\pi(\tilde{\psi},\underline{\tilde{\ell}},\underline{\tilde{\eta}}).\eqno(1)
$$
A droite n'intervient d\'ej\`a que les prolongements qui v\'erifient  $$\forall j\in [1,a]\, \underline{\tilde{\ell}}(\rho,A+T_{2j},B+T_{2j},\zeta)=
\underline{\tilde{\ell}}(\rho,A+T_{2j-1},B+T_{2j-1},\zeta)$$ et $\underline{\eta}(\rho,A+T_{2j},B+T_{2j},\zeta)=
\underline{\eta}(\rho,A+T_{2j-1},B+T_{2j-1},\zeta)(-1)^{(A-B)}$ (cf. \ref{decomposition}). On consid\`ere les \'el\'ements de $Jord(\psi')$ qui sont inf\'erieurs ou \'egaux \`a $(\rho,A,B,\zeta)$ et ${\cal C}'_{\leq}$ les \'el\'ements qui font ``passer'' des blocs de $\tilde{\psi}'$ \`a ces blocs. On applique $Jac_{x\in {\cal C}'_{\leq }}$ aux 2 membres de (1). Cette op\'eration commute avec l'induction par les $S(\rho,A+T_{2j-1},B+T_{2j-1},\zeta)$ car ${\cal C}'_{\leq }$ ne contient aucun terme de la forme $\zeta(B+T_{2j-1})$; cette op\'eration commute aussi au $Jac$ du membre de droite parce que les \'el\'ements des ${\cal C}(\zeta, A+T_{2j-1},B+T_{2j-1},T_{2j}-T_{2j-1})$ sont tr\`es diff\'erents des \'el\'ements de ${\cal C}'_{\leq}$. D'o\`u
$$
\times_{j\in [1,a]}S(\rho,A+T_{2j-1},B+T_{2j-1},\zeta) \times Jac_{x\in {\cal C}'_{\leq}} \pi(\underline{\tilde{\psi}'})=$$
$$\oplus_{\underline{\tilde{\ell}},\underline{\tilde{\eta}}}Jac_{x\in \cup_{j\in [1,a]} {\cal C}(\zeta, A+T_{2j-1},B+T_{2j-1},T_{2j}-T_{2j-1})}Jac_{x\in {\cal C}'_{\leq}}\pi(\tilde{\psi},\underline{\tilde{\ell}},\underline{\tilde{\eta}}).\eqno(2)
$$
On applique maintenant progressivement de $j=1$ \`a $j=a$, $Jac_{x\in {\cal C}(\zeta,A,B,T_{2j-1})\cup {\cal C}(\zeta,A,B,T_{2j-1})}$ aux 2 membres de (2). On v\'erifie que le membre de gauche devient
$$
S(\rho,A,B,\zeta)\times \cdots \times S(\rho,A,B,\zeta)\times Jac_{x\in {\cal C}'_{\leq}} \pi(\underline{\tilde{\psi}'})$$
avec la multiplicit\'e $2^a$.
Il est plus difficile de d\'ecrire le terme de droite; le r\'esultat est que l'on obtient avec multiplicit\'e, $2^a$, les repr\'esentations
$
Jac_{x\in \cup_{i\in [1,2a]}{\cal C}(\rho,A,B,T_{i})}Jac_{x\in {\cal C}'_{\leq }}\pi(\tilde{\psi},\underline{\tilde{\ell}},\underline{\tilde{\eta}}).$

\noindent
Admettons cela pour le moment et concluons. 
Cette assertion donne une \'egalit\'e:
$$
S(\rho,A,B,\zeta)\times \cdots \times S(\rho,A,B,\zeta)\times Jac_{x\in {\cal C}'_{\leq}} \pi(\underline{\tilde{\psi}'})=\oplus_{\underline{\tilde{\ell}},\underline{\tilde{\eta}}}
Jac_{x\in \cup_{i\in [1,2a]}{\cal C}(\rho,A,B,T_{i})}Jac_{x\in {\cal C}'_{\leq }}\pi(\tilde{\psi},\underline{\tilde{\ell}},\underline{\tilde{\eta}}).\eqno(3)$$
ou encore pour tout $\underline{\tilde{\ell}},\underline{\tilde{\eta}}$ tel que $Jac_{x\in \cup_{i\in [1,2a]}{\cal C}(\rho,A,B,T_{i})}Jac_{x\in {\cal C}'_{\leq }}\pi(\tilde{\psi},\underline{\tilde{\ell}},\underline{\tilde{\eta}})$ soit non nul, une inclusion:
$$
Jac_{x\in {\cal C}'_{\leq }}\pi(\tilde{\psi},\underline{\tilde{\ell}},\underline{\tilde{\eta}})\hookrightarrow\times_{j\in [1,2a]}S(\zeta,A,B,T_{j})\times Jac_{x\in {\cal C}'_{\leq}} \pi(\underline{\tilde{\psi}'})
$$Par irr\'eductibilit\'e de la repr\'esentation $\times_{j\in [1,2a]}S(\zeta,A,B,T_{j})$, on peut permuter les facteurs. Fixons $i\in [1,a[$, on peut alors mettre en premi\`ere place (\`a gauche) $S(\zeta,A,B,T_{2j+1})$ et en utilisant l'inclusion de $S(\zeta, A,B,T_{2j+1})$ dans $S(\zeta,A+T_{2j},B+T_{2j},T_{2j+1}-T_{2j})\times S(\zeta,A,B,T_{2j})$, on obtient
$$
Jac_{x\in {\cal C}(\zeta,A+T_{2j},B+T_{2j},T_{2j+1}-T_{2j})}Jac_{x\in {\cal C}'_{\leq }}\pi(\tilde{\psi},\underline{\tilde{\ell}},\underline{\tilde{\eta}})\neq 0.
$$
De \ref{definition}, on en d\'eduit que l'on a aussi $$\underline{\tilde{\ell}}(\rho,A+T_{2j},B+T_{2j},\zeta)= \underline{\tilde{\ell}}(\rho,A+T_{2j+1},B+T_{2j+1},\zeta); $$
$$ \underline{\tilde{\eta}}(\rho,A+T_{2j},B+T_{2j},\zeta)= \underline{\tilde{\eta}}(\rho,A+T_{2j+1},B+T_{2j+1},\zeta)(-1)^{A-B}.
$$
Si $(\rho,A,B,\zeta)\in Jord(\psi')$ on fait la m\^eme chose pour $j=0$ (cf. le d\'ebut de la d\'emonstration) pour obtenir les conditions de l'\'enonc\'e.
On note encore ${\cal C}''$ les \'el\'ements qui font ''passer'' des blocs de $\tilde{\psi}_{a}$ strictement sup\'erieurs \`a $(\rho,A+T_{2a},B+T_{2a},\zeta)$ \`a ceux de $\psi'$ strictement sup\'erieurs \`a $(\rho,A,B,\zeta)$. On applique $Jac_{x\in {\cal C}''}$ aux 2 membres de (3); \`a gauche, cette op\'eration commute \`a l'induction et on obtient donc que les sous-quotients irr\'eductibles de $\sigma_{a}$ sont exactement les repr\'esentations d\'ecrites dans l'\'enonc\'e et avec multiplicit\'e 1; on sait qu'elles sont non isomorphes (ou nulles) gr\^ace \`a \ref{proprietesgenerales} et cela termine la preuve modulo l'assertion admise.

\nl
D\'emontrons cette assertion; soit $\underline{\tilde{\ell}},\underline{\tilde{\eta}}$ tel que $$Jac_{x\in \cup_{j\in [1,a]}{\cal C}(\zeta,A,B,T_{2j-1})\cup {\cal C}(\zeta,A,B,T_{2j-1})}Jac_{x\in \cup_{j\in [1,a]} {\cal C}(\zeta, A+T_{2j-1},B+T_{2j-1},T_{2j}-T_{2j-1})}Jac_{x\in {\cal C}'_{\leq }}\pi(\tilde{\psi},\underline{\tilde{\ell}},\underline{\tilde{\eta}})\neq 0.$$
Alors il existe une repr\'esentation irr\'eductible $\pi$ et une inclusion de
$Jac_{x\in {\cal C}'_{\leq }}\pi(\tilde{\psi},\underline{\tilde{\ell}},\underline{\tilde{\eta}})$ dans
$$\times_{j\in [1,a]}S(\zeta,A+T_{2j-1},B+T_{2j-1},T_{2j}-T_{2j-1}) \times
\biggl(\times_{j\in [1,a]} S(\zeta,A,B,T_{2j-1})\times S(\zeta,A,B,T_{2j-1})\biggr) \times  \pi.
$$
On sait en plus que $Jac_{z}$ de la repr\'esentation de gauche vaut $0$ pour tout $z\in ]\zeta B,\zeta(A+T_{2a})]$ sauf pour les valeurs $\zeta(B+T_{j})$ avec $j\in [1,2a]$. On en d\'eduit que l'inclusion se factorise par 
$$
\times_{i\in [1,2a]}S(\zeta,A,B,T_{i})\times \pi.\eqno(4)
$$
Cela force $\pi=Jac_{x\in \cup_{i\in [1,a]}{\cal C}(\zeta,A,B,T_{i})}Jac_{x\in {\cal C}'_{\leq }}\pi(\tilde{\psi},\underline{\tilde{\ell}},\underline{\tilde{\eta}})$. On peut facilement calculer $$Jac_{x\in \cup_{j\in [1,a]}{\cal C}(\zeta,A,B,T_{2j-1})\cup {\cal C}(\zeta,A,B,T_{2j-1})}Jac_{x\in \cup_{j\in [1,a]} {\cal C}(\zeta, A+T_{2j-1},B+T_{2j-1},T_{2j}-T_{2j-1})}$$
$$\times_{i\in [1,a]}S(\zeta,A,B,T_{i})\times \pi$$et cela vaut $2^a$ fois $\pi$. On a donc d\'emontr\'e que l'on n'a pas plus que ce qui est annonc\'e mais il faut d\'emontrer que l'on a bien toutes ces repr\'esentations. 

On reprend la d\'emonstration en sens inverse, c'est-\`a-dire que l'on suppose que $$Jac_{x\in \cup_{i\in [1,a]}{\cal C}(\zeta,A,B,T_{i})}Jac_{x\in {\cal C}'_{\leq }}\pi(\tilde{\psi},\underline{\tilde{\ell}},\underline{\tilde{\eta}}) \neq 0.$$ On sait cf \ref{irreductibilite} que cette repr\'esentation est irr\'eductible, on la note $\pi$ est on a une inclusion comme dans (4). Puis on utilise encore 
$$
\times_{i\in [1,2a]}S(\zeta,A,B,T_{i})\simeq \times_{j\in [1,a]}S(\zeta, A,B,T_{2j})\times_{j\in [1,a]}S(\zeta,A,B,T_{2j-1})$$et pour $j\in [1,a]$
$$
S(\zeta, A,B,T_{2j})\times_{k\in ]j,a]}S(\zeta,A,B,T_{2k})\hookrightarrow$$
$$ S(\zeta,A+T_{2j-1},B+T_{2j-1},T_{2j}-T_{2j-1})\times S(\zeta,A,B,T_{2j-1})\times_{k\in ]j,a]}S(\zeta,A,B,T_{2k})$$
$$
\simeq S(\zeta,A+T_{2j-1},B+T_{2j-1},T_{2j}-T_{2j-1}) \times_{k\in ]j,a]}S(\zeta,A,B,T_{2k})\times S(\zeta,A,B,T_{2j-1}).
$$ 
Ainsi de proche en proche, on obtient une inclusion
$
\times_{i\in [1,2a]}S(\zeta,A,B,T_{i})\hookrightarrow$
$$ \times_{j\in [1,a]}S(\zeta,A+T_{2j-1},B+T_{2j-1},T_{2j}-T_{2j-1}) \times_{j\in [1,a]}S(\zeta,A,B,T_{2j-1})\times_{j\in [1,a]}S(\zeta,A,B,T_{2j-1}).
$$
Et par irr\'eductibilit\'e, on peut permuter comme on veut les 2a derniers facteurs et on peut donc prolonger l'inclusion (4) en
$$ \times_{j\in [1,a]}S(\zeta,A+T_{2j-1},B+T_{2j-1},T_{2j}-T_{2j-1}) \times_{j\in [1,a]}\biggl(S(\zeta,A,B,T_{2j-1})\times S(\zeta,A,B,T_{2j-1})\biggr)\times \pi.
$$
D'o\`u le fait que $\pi$ intervient dans $$Jac_{x\in \cup_{j\in [1,a]}{\cal C}(\zeta,A,B,T_{2j-1})\cup {\cal C}(\zeta,A,B,T_{2j-1})}Jac_{x\in \cup_{j\in [1,a]} {\cal C}(\zeta, A+T_{2j-1},B+T_{2j-1},T_{2j}-T_{2j-1})}Jac_{x\in {\cal C}'_{\leq }}\pi(\tilde{\psi},\underline{\tilde{\ell}},\underline{\tilde{\eta}})$$avec multiplicit\'e au moins 1. Mais comme elle ne peut intervenir que pour cette valeur de $\underline{\tilde{\ell}}$ et $\underline{\tilde{\eta}}$ elle intervient n\'ecessairement avec multiplicit\'e $2^a$ d'apr\`es la partie directe. Cela termine la preuve.
%%%%%%%
\section{Param\'etrisation et induction\label{parametrisationetinduction}}
Les d\'efinitions de \ref{definition} sont maintenant pleinement justifi\'ees. On peut relier cela \`a l'induction de la fa\c{c}on suivante. Soit $\psi$ tel que tout \'el\'ement de $Jord(\psi)$ ait la condition de parit\'e et notons $\psi_{0}$ le morphisme qui se d\'eduit de $\psi$ en posant $$Jord(\psi_{0}):=\{(\rho,A,B,\zeta)\in Jord(\psi); (-1)^{mult_{Jord(\psi)}(\rho,A,B,\zeta)}=-1\}.$$C'est un ensemble sans multiplicit\'e.  On fixe $\epsilon_{0}$ un caract\`ere du centralisateur de $\psi_{0}$, c'est-\`a-dire un morphisme de $Jord(\psi_{0})$ dans $\pm 1$. 

\

\bf Corollaire: \sl $\oplus_{\epsilon; \epsilon_{\vert Jord(\psi_{0})}=\epsilon_{0}}\pi(\psi,\epsilon)=$
$$\times_{(\rho,A,B,\zeta)\in Jord(\psi)} \underbrace{S(\rho,A,B,\zeta)\times \cdots \times S(\rho,A,B,\zeta)}_{1/2(mult_{Jord(\psi)}(\rho,A,B,\zeta)-mult_{Jord(\psi_{0})}(\rho,A,B,\zeta))} \times \pi(\psi_{0},\epsilon_{0}).$$\rm

\

C'est un corollaire imm\'ediat de \ref{inductionsuite}

%%%%%%%
\section {Induction, le cas g\'en\'eral\label{inductionfin}}
Ici on fixe $\psi$ tel que tous les \'el\'ements de $Jord(\psi)$ ait la condition de parit\'e. Pour tout choix d'application $\underline{\ell}$ et $\underline{\eta}$ compatible avec un caract\`ere, $\epsilon$ du centralisateur de $\psi$ (cf. \ref{definition}), on d\'efinit $\pi(\psi,\epsilon,\underline{\ell},\underline{\eta})$ et pour simplifier on pose $\underline{\psi}$ cet ensemble de donn\'ees $\psi,\epsilon,\underline{\ell},\underline{\eta}$.

Et on fixe une collection, ${\cal E}$, de quadruplets $(\rho,A,B,\zeta)$ dont aucun n'a la condition de parit\'e.

\

\bf Th\'eor\`eme: \sl L'induite $\times_{(\rho,A,B,\zeta)\in {\cal E}}S(\rho,A,B,\zeta) \times \pi(\underline{\psi})$ est irr\'eductible. \rm

\

Pour simplifier les notations de la d\'emonstration, on suppose que pour tout $(\rho,A,B,\zeta)\in {\cal E}$, $\rho\simeq \rho^*$, ce qui est le cas le plus difficile.
On fait une premi\`ere r\'eduction: on ordonne totalement les \'el\'ements de ${\cal E}$ de telle sorte que si $(\rho,A,B,\zeta),(\rho,A',B',\zeta')\in {\cal E}$ avec $B>B'$ alors le premier terme est strictement sup\'erieur au deuxi\`eme. Il y a de nombreux choix, d'autant plus que ${\cal E}$ peut avoir de la multiplicit\'e. Pour tout $(\rho,A,B,\zeta)\in {\cal E}$ on fixe $T_{\rho,A,B,\zeta}\in {\mathbb N}$ de telle sorte que si $(\rho,A,B,\zeta)>(\rho',A',B,',\zeta')$ alors $B+T_{\rho,A,B,\zeta}>A'+T_{\rho',A',B',\zeta'}$. On suppose que l'on sait d\'emontrer que 
$$\sigma_{T}:=
\times _{(\rho,A,B,\zeta)\in {\cal E}}S(\rho,A+T_{\rho,A,B,\zeta},B+T_{\rho,A,B,\zeta},\zeta) \times \pi(\underline{\psi})$$est irr\'eductible et on va en d\'eduire le th\'eor\`eme.

Fixons $(\rho,A,B,\zeta)\in {\cal E}$ et posons:
$$
\sigma_{\geq (\rho,A,B,\zeta)}:= \times_{(\rho',A',B',\zeta')\in {\cal E}; (\rho',A',B',\zeta')\geq (\rho,A,B,\zeta)}S(\rho',A'+T',B'+T',\zeta') $$
$$\times _{(\rho',A',B',\zeta')\in {\cal E}; (\rho',A',B',\zeta')<(\rho,A,B,\zeta)}S(\rho',A',B',\zeta')  \times \pi(\underline{\psi});$$
$$
\sigma_{> (\rho,A,B,\zeta)}:= \times_{(\rho',A',B',\zeta')\in {\cal E}; (\rho',A',B',\zeta')>(\rho,A,B,\zeta)}S(\rho',A'+T',B'+T',\zeta')$$
$$ \times _{(\rho',A',B',\zeta')\in {\cal E}; (\rho',A',B',\zeta')\leq (\rho,A,B,\zeta)}S(\rho',A',B',\zeta')  \times \pi(\underline{\psi}).
$$
On admet que $\sigma_{\geq (\rho,A,B,\zeta)}$ est irr\'eductible et on d\'emontre que $\sigma_{>(\rho,A,B,\zeta)}$ est aussi irr\'eductible.

En effet, on calcule ais\'ement:
$$
Jac_{x\in {\cal C}(\zeta,A,B,T_{\rho,A,B,\zeta})\cup {\cal C}(\zeta,A,B,T_{\rho,A,B,\zeta})}\sigma_{\geq (\rho,A,B,\zeta)}.
$$
Montrons plus g\'en\'eralement l'assertion suivante: soit $\tau$ une repr\'esentation (non n\'ecessairement irr\'eductible) et soit $(\rho,A,B,\zeta)$  $T\in {\mathbb N}$ comme ci-dessus. On suppose que pour tout $x\in {\cal C}(\zeta,A,B,T)$, $Jac_{x}\tau=0$ alors, dans le groupe de Grothendieck:
$$
Jac_{x\in {\cal C}(\zeta,A,B,T)\cup {\cal C}(\zeta,A,B,T)}\biggl(S(\rho,A+T,B+T,\zeta)\times \tau\biggr)=S(\rho,A,B,\zeta)\times \tau \oplus S(\rho,A,B,\zeta)\times \tau .
$$
En effet, on calcule d'abord $Jac_{x\in {\cal C}(\zeta,A,B,T)}\biggl(S(\rho,A+T,B+T,\zeta)\times \tau\biggr)$. On commence par calculer $Jac_{x\in [\zeta (B+T), \zeta (A+T)]}$. Pour cela on d\'ecompose $[\zeta (B+T), \zeta (A+T)]$ en trois sous-ensembles ${\cal I}_{j}$ pour $j\in \{1,2,3\}$ v\'erifiant:
$$
Jac_{x\in {\cal I}_{1}}Jac^d_{x\in -^t{\cal I}_{2}}S(\rho,A+T,B+T,\zeta)\neq 0; \qquad Jac_{x\in {\cal I}_{3}}\tau\neq 0.
$$
Par l'hypoth\`ese que l'on a faite ${\cal I}_{3}$ est n\'ecessairement vide puisque qu'il est inclus dans ${\cal C}(\zeta,A,B,T)$. Si ${\cal I}_{1}$ est non vide, il contient n\'ecessairement $\zeta(B+T)$ et si ${\cal I}_{2}$ est non vide, $-^t{\cal I}_{2}$ contient n\'ecessairement $-\zeta (B+T)$ c'est-\`a-dire que ${\cal I}_{2}$ contient $\zeta (B+T)$. Mais comme l'intervalle ne contient qu'une fois $\zeta (B+T)$,  ${\cal I}_{j}
$ est vide pour une valeur de $j=1,2$ et vaut $[\zeta (B+T),\zeta (A+T)]$ pour l'autre valeur. Dans le cas o\`u c'est ${\cal I}_{1}$ qui est non vide, on enl\`eve la premi\`ere colonne du tableau d\'efinissant $S(\rho,A+T,B+T,\zeta)$ et dans le deuxi\`eme cas on enl\`eve la derni\`ere colonne. Donc $$
Jac_{x\in [\zeta (B+T),\zeta (A+T)}\biggl(S(\rho,A+T,B+T,\zeta)\times \tau\biggr)=$$
$$ \biggl(\begin{matrix}
&\zeta (B+T-1) &\cdots &-\zeta (A+T)\\
&\vdots &\vdots &\vdots\\
&\zeta (A+T-1) &\cdots &-\zeta (B+T)
\end{matrix}\biggr) \times \tau
\oplus
\biggl(\begin{matrix}
&\zeta (B+T) &\cdots &-\zeta (A+T-1)\\
&\vdots &\vdots &\vdots\\
&\zeta (A+T) &\cdots &-\zeta (B+T-1)
\end{matrix}\biggr) \times \tau.
$$Dans le groupe de Grothendieck les 2 induites sont les m\^emes. Ensuite on calcule le $Jac_{x\in {\cal C}(\zeta,A,B,T-1)}$ de cette induite. Comme $\zeta (B+T)$ n'est plus dans l'ensemble ${\cal C}(\zeta,A,B,T-1)$ il n'y a plus qu'une possibilit\'e enlever les $T-1$ premi\`eres colones quand l'induite est \'ecrite sous la premi\`ere forme (les $T-1$ derni\`eres colonnes sous la deuxi\`eme forme). Finalement on trouve:
$$
Jac_{x\in {\cal C}(\zeta,A,B,T)}\biggl(S(\rho,A+T,B+T,\zeta)\times \tau\biggr)=$$
$$
 \biggl(\begin{matrix}
&\zeta B&\cdots &-\zeta (A+T)\\
&\vdots &\vdots &\vdots\\
&\zeta A &\cdots &-\zeta (B+T)
\end{matrix}\biggr) \times \tau
\oplus
\biggl(\begin{matrix}
&\zeta (B+T) &\cdots &-\zeta A\\
&\vdots &\vdots &\vdots\\
&\zeta (A+T) &\cdots &-\zeta B
\end{matrix}\biggr) \times \tau.
$$Dans le groupe de Grothendieck les 2 induites sont les m\^emes. Il reste encore \`a calculer $Jac_{x\in {\cal C}(\zeta,A,B,T)}$, cela se fait de la m\^eme fa\c{c}on et finalement on trouve la repr\'esentation $S(\rho,A,B,T)\times \tau$ avec multiplicit\'e 2.

On applique ce calcul,  en prenant pour $$\tau:=\times_{(\rho',A',B',\zeta')\in {\cal E}; (\rho',A',B',\zeta')>(\rho,A,B,\zeta)}S(\rho',A'+T',B'+T',\zeta') $$
$$\times _{(\rho',A',B',\zeta')\in {\cal E}; (\rho',A',B',\zeta')<(\rho,A,B,\zeta)}S(\rho',A',B',\zeta')  \times \pi(\underline{\psi}).$$
Les hypoth\`eses sont satisfaites, car $Jac_{x}\pi(\underline{\psi})\neq 0$ n\'ecessite que $x$ soit entier si les \'el\'ements de ${\cal C}(\zeta,A,B,T)$ sont demi-entiers non entiers et vice et versa et $$Jac_{x}\biggl(\times_{(\rho',A',B',\zeta')\in {\cal E}; (\rho',A',B',\zeta')>(\rho,A,B,\zeta)}S(\rho',A'+T',B'+T',\zeta')$$
$$ \times _{(\rho',A',B',\zeta')\in {\cal E}; (\rho',A',B',\zeta')<(\rho,A,B,\zeta)}S(\rho',A',B',\zeta')\biggr)\neq 0$$n\'ecessite que $x=\zeta' (B'+T')$ avec $B'+T'>A+T$ ou $x=\zeta' B'$ avec $B'\leq B<(B+1)$. Or les \'el\'ements de ${\cal C}(\zeta,A,B,T)$ sont de valeur absolue comprise entre $(B+1)$ et $(A+T)$. Il faut la m\^eme propri\'et\'e avec $Jac^d_{-x}$  mais c'est automatique, puisque la repr\'esentation est autoduale. Pour vraiment se ramener \`a notre calcul, il faut encore utiliser que les repr\'esentations $S(\rho_{1},C_{1},C_{1}',\zeta_{1})\times S(\rho_{2},C_{2},C_{2}',\zeta_{2})$ sont irr\'eductibles et donc permutables (cf \ref{unresultatdirreductibilite} pour ce r\'esultat).

Finalement on trouve que $Jac_{x\in {\cal C}(\zeta,A,B,T)\cup {\cal C}(\zeta,A,B,T)}\sigma_{\geq (\rho,A,B,T) }$ vaut  la repr\'esentation $\sigma_{>(\rho,A,B,T)}$ avec multiplicit\'e 2. Il en existe donc un sous-quotient irr\'eductible, not\'e $\sigma'$, de $\sigma_{>(\rho,A,B,T)}$ tel que par r\'eciprocit\'e de Frobenius:
$$
\sigma_{\geq (\rho,A,B,\zeta)}\hookrightarrow
\times_{x\in{\cal C}(\zeta,A,B,T)\cup {\cal C}(\zeta,A,B,T)}\rho\vert\,\vert^x \times \sigma'.
$$
Avec la propri\'et\'e que $Jac_{x}\sigma_{\geq (\rho,A,B,\zeta)}=0$ pour tout $x\in {\cal C}(\rho,A,B,\zeta)$ sauf \'eventuellement $x=\zeta (B+T)$, on voit que l'inclusion se factorise par:
$$
\sigma_{\geq (\rho,A,B,\zeta)}\hookrightarrow S(\zeta,A,B,T)\times S(\zeta,A,B,T)\times \sigma'.
$$
On v\'erifie sans probl\`eme que $Jac_{x}\sigma_{>(\rho,A,B,\zeta)}=0$ pour tout $x\in {\cal C}(\zeta,A,B,T)$ (on l'a essentiellement fait ci-dessus) et on v\'erifie essentiellement comme ci-dessus que
$$
Jac_{x\in {\cal C}(\zeta,A,B,T)\cup {\cal C}(\zeta,A,B,T)}S(\zeta,A,B,T)\times S(\zeta,A,B,T)\times \sigma'$$vaut 2 fois $\sigma'$. Or ce module de Jacquet doit contenir celui de $\sigma_{\geq (\rho,A,B,\zeta)}$ et donc $\sigma'$ doit co\"{\i}ncider avec $\sigma_{>(\rho,A,B,T)}$. D'o\`u l'irr\'eductibilit\'e annonc\'ee.

\

Dans la suite de la d\'emonstration, on suppose donc que pour tout $(\rho,A,B,\zeta),(\rho',A',B',\zeta') \in {\cal E}$ on a soit $B>A'$ soit $B'>A$; on suppose m\^eme que ces nombres sont tr\`es disjoints. On se ram\`ene aussi au cas o\`u $\psi$ est de restriction discr\`ete \`a la diagonale comme dans la preuve de \ref{inductionsuite}.

On va d'abord d\'emontrer le th\'eor\`eme dans le cas o\`u $\psi$ est trivial sur la 2e copie de $SL(2,{\mathbb C})$. En particulier, il n'y a ici que $\epsilon$ et $\pi(\psi,\epsilon)$ est une s\'erie discr\`ete. On a aussi suppos\'e que la notion de mauvaise parit\'e pour $\rho$ autodudale fix\'e, traduit le fait que les induites:
$$
<\rho\vert\,\vert^{C}, \cdots, \rho\vert\,\vert^{-C}>\times \pi(\psi'
,\epsilon)\eqno(1)$$sont irr\'eductibles soit pour tout $C$ entier (les entiers ont alors la mauvaise parit\'e, il faut penser \`a la parit\'e de $2C+1$) soit  pour tout $C$ demi-entier non entier, ce sont alors les demi-entiers non entiers qui ont la mauvaise parit\'e. 
Les r\'esultats d'Harish-Chandra (cf. \cite{JIMJ}) g\'en\'eralisent alors cette assertion quand on remplace le segment $[C,-C]$ par un ensemble fini de segments d\'ecroissants. On veut la m\^eme assertion pour un ensemble de segments, ${\cal E}$,  de la forme $[\zeta C,-\zeta C]$ o\`u $\zeta=+$ ou $-$; je n'ai pas trouv\'e d'argument \'el\'ementaire pour prouver cela. La d\'emonstration que je propose consiste \`a utiliser l'application $inv_{<\alpha}\circ inv_{\leq \alpha}$ de \cite{elementaire} 4.1 pour se ramener au cas des segments tous d\'ecroissants; pour r\'eussir, il faut que l'on ne puisse pas avoir dans l'ensemble des segments \`a la fois un segment $[C,-C]$ et un segment $[-C,C]$; on a  fait ce qu'il fallait ci-dessus pour se ramener au cas o\`u cette \'eventualit\'e ne se produit pas. Ensuite on ne peut pas invoquer directement l'assertion de \cite{elementaire} 4.1 pour savoir que cette application conserve l'irr\'eductibilit\'e (au signe pr\`es) des sous-quotients irr\'eductibles de l'induite
$$
\times _{[\zeta C,-\zeta C]\in {\cal E}}<\rho\vert\,\vert^{\zeta C}, \cdots, \rho\vert\,\vert^{-\zeta C}> \times \pi(\psi,\epsilon).
$$
Mais fort heureusement la d\'emonstration s'applique telle quelle pour $inv_{\leq \alpha}$ quand on remarque que dans l'\'enonc\'e de l'assertion de loc. cite, on peut parfaitement remplacer $\pi(\psi,\epsilon)$ par $$\times_{[\zeta C,-\zeta C]\in {\cal E}; C>\alpha}<\rho\vert\,\vert^{\zeta C}, \cdots, \rho\vert\,\vert^{-\zeta C}> \times \pi(\psi,\epsilon).$$Le parabolique qui intervient est exactement celui de la preuve donn\'ee. Cela permet de se ramener au cas temp\'er\'e.

\nl On suppose donc pour le moment que $\pi(\psi,\epsilon)$ est une s\'erie discr\`ete.
Fixons $(\rho,A,B,\zeta)\in {\cal E}$. On montre d'abord que l'induite:
$$\tau:=
<\rho\vert\,\vert^{\zeta A}, \cdots, \rho\vert\,\vert^{-\zeta B}>\times \pi(\psi,\epsilon)$$est irr\'eductible. Pour cela on remarque d'abord qu'elle a un unique sous-module irr\'eductible, $\tau_{s}$: cela r\'esulte des 2 faits suivants. D'une part pour tout $x\in [\zeta A, \zeta(B+1)]$, $Jac_{x}\pi(\psi,\epsilon)=0$, d'o\`u
$$
Jac_{\zeta A, \cdots, \zeta (B+1)}\tau=<\rho\vert\,\vert^{\zeta B}, \cdots, \rho\vert\,\vert^{-\zeta B}>\times \pi(\psi,\epsilon).$$
Et d'autre part, l'induite $<\rho\vert\,\vert^{\zeta B}, \cdots, \rho\vert\,\vert^{-\zeta B}>\times \pi(\psi,\epsilon)$ est irr\'eductible. 

\noindent
L'induite $<\rho\vert\,\vert^{\zeta B}, \cdots, \rho\vert\,\vert^{-\zeta A}>\times \pi(\psi,\epsilon)$ a aussi un unique sous-module irr\'eductible car
$$
Jac_{\zeta B, \cdots, -\zeta A} <\rho\vert\,\vert^{\zeta B}, \cdots, \rho\vert\,\vert^{-\zeta A}>\times \pi(\psi,\epsilon)=\pi(\psi,\epsilon)$$
par r\'eciprocit\'e de Frobenius, on voit que $\tau$ a aussi un unique quotient irr\'eductible $\tau_{q}$. Il suffit donc de d\'emontrer que $\tau_{s}\simeq \tau_{q}$ puisque chacun intervient avec multiplicit\'e 1 comme sous-quotient irr\'eductible de $\tau$ d'apr\`es les calculs de modules de Jacquet faits ci-dessus.

On consid\`ere $\tilde{\tau}:=<\rho\vert\,\vert^{\zeta A}, \cdots, \rho\vert\,\vert^{-\zeta A}> \times \pi(\psi,\epsilon)$. C'est  une repr\'esentation irr\'eductible par hypoth\`ese. On calcule le module de Jacquet, $Jac_{\zeta A, \cdots, \zeta (B+1)}\tilde{\tau}$. On a une filtration:
$$
0 \rightarrow <\rho\vert\,\vert^{\zeta A}, \cdots, \rho\vert\,\vert^{\zeta(B+1)}>\otimes \biggl(<\rho\vert\,\vert^{\zeta A}, \cdots, \rho\vert\,\vert^{-\zeta B}>\times \pi(\psi,\epsilon)\biggr) \rightarrow$$
$$
Jac_{\zeta A, \cdots, \zeta (B+1)}\tilde{\tau} \rightarrow <\rho\vert\,\vert^{\zeta A}, \cdots, \rho\vert\,\vert^{\zeta(B+1)}>\otimes \biggl(<\rho\vert\,\vert^{\zeta B}, \cdots, \rho\vert\,\vert^{-\zeta A}>\times \pi(\psi,\epsilon) \biggr)\rightarrow 0.\eqno(2)
$$
Ici $ <\rho\vert\,\vert^{\zeta A}, \cdots, \rho\vert\,\vert^{\zeta(B+1)}>\otimes \tau_{s}$ est quotient de la repr\'esentation la plus \`a droite; d'o\`u par r\'eciprocit\'e de Frobenius un homorphisme non nul:
$$
\tilde{\tau} \rightarrow <\rho\vert\,\vert^{\zeta A}, \cdots, \rho\vert\,\vert^{\zeta(B+1)}>\times \tau_{s}.\eqno(3)
$$Par irr\'eductibilit\'e ce morphisme est une inclusion. 
On recalcule dans le groupe de Grothendieck $Jac_{\zeta A, \cdots, \zeta (B+1)}$ du 2e membre de (3); on obtient la somme de $\tau_{s}$ avec $<\rho\vert\,\vert^{\zeta A}, \cdots, \rho\vert\,\vert^{\zeta (B+1)} >\times Jac_{\zeta A, \cdots, \zeta (B+1)}\tau_{s}$. Or $Jac_{\zeta A, \cdots, \zeta (B+1)}\tau_{s}$ est la repr\'esentation irr\'eductible $<\rho\vert\,\vert^{\zeta B}, \cdots, \rho\vert\,\vert^{-\zeta B}>\times \pi(\psi,\epsilon)$. Ainsi, dans le groupe de Grothendieck, $Jac_{\zeta A,\cdots, \zeta (B+1)}\tilde{\tau}$ est au plus la somme de $\tau_{s}$ et de $$<\rho\vert\,\vert^{\zeta A}, \cdots, \rho\vert\,\vert^{\zeta (B+1)}>\times <\rho\vert\,\vert^{\zeta B}, \cdots, \rho\vert\,\vert^{-\zeta B}>\times \pi(\psi,\epsilon).$$
On suppose que $\tau_{s}\not\simeq \tau_{q}$ et on calcule la multiplicit\'e de $\tau_{q}$ dans $Jac_{\zeta A, \cdots, \zeta (B+1)}\tilde{\tau}$ avec l'inclusion (3). D'apr\`es ce que l'on vient d'\'ecrire; cette multiplicit\'e est inf\'erieure ou \'egal \`a  1. Par contre dans la suite exacte \'ecrite en (2), elle est 2; d'o\`u la contradiction qui donne $\tau_{s}\simeq \tau_{q}$ et l'irr\'eductibilit\'e.

 On revient \`a ${\cal E}$ et on \'ecrit ${\cal E}':={\cal E}-\{(\rho,A,B,\zeta)\}$. On rappelle, car on l'utilisera librement que pour tout $(\rho,A',B',\zeta'),(\rho,A'',B'',\zeta'')\in {\cal E}$ l'induite $S(\rho,A',B',\zeta')\times S(\rho,A'',B'',\zeta'')$ est irr\'eductible. On suppose que $A>B$, sinon on change d'\'el\'ement  et si pour tout \'el\'ement dans ${\cal E}$, $A=B$, on a d\'ej\`a montr\'e le r\'esultat. 
 
 Soit $\tau$ un sous-module irr\'eductible de $\times_{(\rho',A',B',\zeta')\in{\cal E}}S(\rho',A',B',\zeta') \times \pi(\psi,\epsilon)$. On a la suite de morphismes:
 $$
 \tau \hookrightarrow \times_{(\rho',A',B',\zeta')\in {\cal E}'}S(\rho',A',B',\zeta') \times S(\rho,A,B,\zeta)\times \pi(\psi,\epsilon) \hookrightarrow
 $$
 $$
 \times_{(\rho',A',B',\zeta')\in {\cal E}'}S(\rho',A',B',\zeta')\times \biggl(
 \begin{matrix} &\zeta B &\cdots &-\zeta A\\
 &\vdots &\vdots &\vdots\\
 &\zeta (A-1) &\cdots &-\zeta (B+1)
 \end{matrix}\biggr) \times <\rho\vert\,\vert^{\zeta A}, \cdots, \rho\vert\,\vert^{-\zeta B}> \times \pi(\psi,\epsilon)$$d'apr\`es l'irr\'eductibilit\'e qui vient d'\^etre d\'emontr\'ee
 $$
 \simeq  \times_{(\rho',A',B',\zeta')\in {\cal E}'}S(\rho',A',B',\zeta')\times \biggl(
 \begin{matrix} &\zeta B &\cdots &-\zeta A\\
 &\vdots &\vdots &\vdots\\
 &\zeta (A-1) &\cdots &-\zeta (B+1)
 \end{matrix}\biggr) \times <\rho\vert\,\vert^{\zeta B}, \cdots, \rho\vert\,\vert^{-\zeta A}>\times \pi(\psi,\epsilon)$$
 d'apr\`es \ref{undeuxiemeresultatdirreductibilite}
 $$
 \simeq
  \times_{(\rho',A',B',\zeta')\in {\cal E}'}S(\rho',A',B',\zeta')\times <\rho\vert\,\vert^{\zeta B}, \cdots, \rho\vert\,\vert^{-\zeta A}>\times \biggl(
 \begin{matrix} &\zeta B &\cdots &-\zeta A\\
 &\vdots &\vdots &\vdots\\
 &\zeta (A-1) &\cdots &-\zeta (B+1)
 \end{matrix}\biggr) \times\pi(\psi,\epsilon)$$
 $$
 \hookrightarrow \times_{(\rho',A',B',\zeta')\in {\cal E}'}S(\rho',A',B',\zeta')\times <\rho\vert\,\vert^{\zeta B}, \cdots, \rho\vert\,\vert^{-\zeta A}>\times <\rho\vert\,\vert^{\zeta B}, \cdots, \rho\vert\,\vert^{-\zeta A}>$$
 $$
 \times \biggl(
 \begin{matrix} &\zeta (B+1) &\cdots &-\zeta (A-1)\\
 &\vdots &\vdots &\vdots\\
 &\zeta (A-1) &\cdots &-\zeta (B+1)
 \end{matrix}\biggr) \times\pi(\psi,\epsilon)
 $$
 $$
 =
 \times_{(\rho',A',B',\zeta')\in {\cal E}'}S(\rho',A',B',\zeta')\times <\rho\vert\,\vert^{\zeta B}, \cdots, \rho\vert\,\vert^{-\zeta A}>\times <\rho\vert\,\vert^{\zeta B}, \cdots, \rho\vert\,\vert^{-\zeta A}>$$
 $$
 \times S(\rho,A-1,B+1,\zeta)\times \pi(\psi,\epsilon).
 $$
 On peut encore faire commuter $<\rho\vert\,\vert^{\zeta B}, \cdots, \rho\vert\,\vert^{-\zeta A}>$ au dessus de chaque $S(\rho',A',B',\zeta')$ car les segments qui d\'efinissent $S(\rho',A',B',\zeta')$ sont soit inclus dans le segment $[\zeta B,-\zeta A]$ (dans le cas o\`u $A'<B$) soit contiennent ce segment (cas o\`u $A<B'$). Donc finalement on a:
 $$
 \sigma\hookrightarrow <\rho\vert\,\vert^{\zeta B}, \cdots, \rho\vert\,\vert^{-\zeta A}>\times <\rho\vert\,\vert^{\zeta B}, \cdots, \rho\vert\,\vert^{-\zeta A}>
 \times_{(\rho',A',B',\zeta')\in {\cal E}'}S(\rho',A',B',\zeta')\times$$
 $$  S(\rho,A-1,B+1,\zeta)\times \pi(\psi,\epsilon).\eqno(1)
$$
Par r\'ecurrence sur $\sum_{(\rho',A',B',\zeta')\in {\cal E}}(A'-B')$ on admet que la repr\'esentation 
$$ \times_{(\rho',A',B',\zeta')\in {\cal E}'}S(\rho',A',B',\zeta')\times  S(\rho,A-1,B+1,\zeta)\times \pi(\psi,\epsilon)$$ est irr\'eductible. On montre que l'induite de droite dans (1) a un unique sous-module irr\'eductible et qu'il intervient avec multiplicit\'e 1 comme sous-quotient de cette induite: pour cela on calcule $Jac_{x\in [\zeta B, -\zeta A] \cup [\zeta B, -\zeta A]}$ de cette induite. Il faut le faire en 2 temps, d'abord calculer $Jac_{x\in [\zeta B, -\zeta A]}$; on d\'ecompose ce segment en 5 sous-ensembles ${\cal I}_{j}$ pour $j\in [1,5]$ tels que
$$
Jac_{x\in {\cal I}_{1}}Jac^d_{x\in -^t{\cal I}_{2}}<\rho\vert\,\vert^{\zeta B}, \cdots, \rho\vert\,\vert^{-\zeta A}>\neq 0
$$
$$
Jac_{x\in {\cal I}_{3}}Jac^d_{x\in -^t{\cal I}_{4}}<\rho\vert\,\vert^{\zeta B}, \cdots, \rho\vert\,\vert^{-\zeta A}>\neq 0
$$
$$
Jac_{x\in {\cal I}_{5}} \times_{(\rho',A',B',\zeta')\in {\cal E}'}S(\rho',A',B',\zeta')\times S(\rho,A-1,B+1,\zeta)\times \pi(\psi,\epsilon)\neq 0.
$$
On cherche quel ensemble peut contenir $-\zeta A$; certainement ni ${\cal I}_{2}$ ni ${\cal I}_{4}$ tout simplement car $\zeta A$ n'est pas dans le support cuspidal de la repr\'esentation $<\rho\vert\,\vert^{\zeta B}, \cdots, \rho\vert\,\vert^{-\zeta A}>$. Il ne peut pas non plus \^etre dans ${\cal I}_{5}$ car il faudrait qu'il existe $(\rho,A',B',\zeta')\in {\cal E}'$ avec $\zeta' B' \in [\zeta B,-\zeta A]$ et $-\zeta A \in [\zeta'A',-\zeta'A']$. Or la premi\`ere condition entra\^{\i}ne $B'\leq A$ et nos hypoth\`eses forcent alors $A'<B$; la 2e condition devient alors impossible \`a r\'ealiser.

Soit maintenant $j=1$ ou $3$ tel que ${\cal I}_{j}$ contienne $-\zeta A$, n\'ecessairement pour cette valeur de $j$, ${\cal I}_{j}=[\zeta B, -\zeta A]$ et le module de Jacquet cherch\'e est, pour ce choix, l'induite:
$$
<\rho\vert\,\vert^{\zeta B}, \cdots, \rho\vert\,\vert^{-\zeta A}>
 \times_{(\rho',A',B',\zeta')\in {\cal E}'}S(\rho',A',B',\zeta')\times  S(\rho,A-1,B+1,\zeta)\times \pi(\psi,\epsilon).
 $$
Quand on applique encore $Jac_{x\in [\zeta B, -\zeta A]}$ \`a ce r\'esultat, l'argument que l'on vient de donner montre qu'on obtient $\tau':= \times_{(\rho',A',B',\zeta')\in {\cal E}'}S(\rho',A',B',\zeta')\times  S(\rho,A-1,B+1,\zeta)\times \pi(\psi,\epsilon).$ Le module de Jacquet de l'induite de (1) cherch\'e est donc $\tau'$ avec multiplicit\'e 2.

Comme l'induite dans le groupe lin\'eaire convenable $<\rho\vert\,\vert^{\zeta B}, \cdots, -\zeta A>\times <\rho\vert\,\vert^{\zeta B}, \cdots, \rho\vert\,\vert^{-\zeta A}>$ est irr\'eductible, il est automatique que $Jac_{x\in {[\zeta B, -\zeta A]}\cup [\zeta B,-\zeta A]}\tau$ contienne au moins $\tau'$ avec multiplicit\'e 2. Ainsi $\tau$ est l'unique sous-module irr\'eductible de l'induite de droite de (1) et par exactitude du foncteur de Jacquet, y intervient au plus avec multiplicit\'e 1 comme sous-quotient. D'o\`u nos assertions et le fait que $\tau$ est compl\`etement d\'etermin\'ee par $\tau'$. Ainsi $\sigma$ a un unique sous-module irr\'eductible et comme $\sigma$ aussi est un sous-module de l'induite de droite de (1), $\tau$ y intervient avec multiplicit\'e 1 comme sous-quotient irr\'eductible. Cela prouve l'irr\'eductibilit\'e de $\sigma$, argument que nous avons d\'ej\`a utilis\'e: $\sigma$ \'etant autodual tout sous-module irr\'eductible de $\sigma$ qui intervient avec multiplicit\'e 1 comme sous-quotient irr\'eductible est n\'ecessairement facteur direct. Comme $\sigma$ n'a qu'un seul sous-module irr\'eductible et que ce sous-module intervient avec multiplicit\'e un comme sous-quotient irr\'eductible de $\sigma$, $\sigma$ est irr\'eductible.

\

Ensuite on passe du cas des s\'eries discr\`etes au cas g\'en\'eral  comme expliqu\'e dans la preuve de \ref{findelapreuve}.
\section{Conclusion\label{conclusion}}
Soit $\psi$ un morphisme de $W_{F}\times SL(2,{\mathbb C})\times SL(2,{\mathbb C})$ dans $^LG$ soumis aux propri\'et\'es de \ref{definitiongenerale}. Soit $\epsilon$ un caract\`ere du centralisateur de $\psi$ comme dans \ref{definitiongenerale}. On a construit dans les paragraphes pr\'ec\'edents, une repr\'esentation semi-simple sans multiplicit\'e, $\pi(\psi,\epsilon)$ qui a le mauvais gout de pouvoir \^etre nulle mais dont on sait, gr\^ace \`a \cite{selecta} par 4, que, modulo les propri\'et\'es de transfert de fonction conjectur\'ees mais pas encore connues, le caract\`ere
$$
\sum_{\epsilon}\epsilon(\psi(1,1,\bigl(\begin{matrix}-1&0\\0&-1&\end{matrix}\bigr)) \, tr\, \pi(\psi,\epsilon)
$$
se transf\`ere en le caract\`ere de la repr\'esentation du GL(n,F) convenable d\'efinie par $\psi$.

Il me semble plus naturelle de d\'ecomposer $\psi$ en repr\'esentations irr\'eductibles de la forme $\rho\otimes [a]\otimes [b]$, o\`u $\rho$ est une repr\'esentation irr\'eductible de $W_{F}$ et $a,b$ sont des entiers donnant chacun la dimension de la repr\'esentation de $SL(2,{\mathbb C})$ \'ecrite entre crochets. Avec ces notations, $Jord(\psi)$ est l'ensemble des triplets $(\rho,a,b)$; c'est un ensemble avec multiplicit\'e. On dit que $(\rho,a,b)$ a la bonne parit\'e si la repr\'esentation $\rho\otimes [a]\otimes [b]$ est \`a valeurs dans un groupe de m\^eme type que $^LG$ et on note $Jord(\psi)_{p}$ le sous-ensemble de $Jord(\psi)$ form\'e des triplets ayant bonne parit\'e; c'est encore un ensemble avec multiplicit\'e. Et $\epsilon$ s'identifie \`a une fonction de $Jord(\psi)$ (vu sans multiplicit\'e) dans $\{\pm 1\}$ qui vaut 1 sur tous les triplets n'ayant pas la bonne parit\'e. Pour d\'ecomposer $\pi(\psi,\epsilon)$, on a utilis\'e des fonctions annexes, $\underline{\ell}$ de $Jord(\psi)$ dans ${\mathbb N}$ v\'erifiant
$$
\forall (\rho,a,b)\in Jord(\psi)_{p}, \quad \underline{\ell}(\rho,a,b)\in [0,[inf(a,b)/2]]
$$
et une fonction $\underline{\eta}$ d\'efinie elle aussi sur $Jord(\psi)_{p}$ est \`a valeurs dans $\{\pm 1\}$ v\'erifiant
$$
\forall(\rho,a,b)\in Jord(\psi) \quad \underline{\eta}(\rho,a,b)^{inf(a,b)}(-1)^{[inf(a,b)/2]+\underline{\ell}(\rho,a,b)}=\epsilon(\rho,a,b).
$$
A de telles donn\'ees, on a associ\'e une sous-repr\'esentation irr\'eductible ou nulle de $\pi(\psi,\epsilon)$ not\'ee $\pi(\psi,\epsilon,\underline{\ell},\underline{\eta})$. Il n'y a qu'un seul cas, o\`u je peux d\'emontrer que ces repr\'esentations sont toutes non nulles; on  note $\psi_{0}$ le morphisme analogue \`a $\psi$ mais tel que $Jord(\psi_{0})$ soit sans multiplicit\'e et co\"{\i}ncide avec $Jord(\psi)_{p}$ aux multiplicit\'e pr\`es. L'hypoth\`ese est que $\psi_{0}$ soit un morphisme de restriction discr\`ete \`a la diaogonale. Dans ce cas, on a:
\nl
\bf Th\'eor\`eme: \sl on suppose que $\psi_{0}$ d\'efinie ci-dessus est de restriction discr\`ete \`a la diagonale, alors les repr\'esentations $\pi(\psi,\epsilon,\underline{\ell},\underline{\eta})$ sont toutes non nulles pour les donn\'ees ci-dessus et d\'ecomposent la repr\'esentation $\pi(\psi,\epsilon)$. \rm
\nl
Etant donn\'e les propri\'et\'es d'irr\'eductibilit\'e d\'ej\`a d\'emontr\'ees pour les induites et \ref{parametrisationetinduction}, il faut prouver le th\'eor\`eme dans le cas o\`u $\psi=\psi_{p}$ et les multiplicit\'es dans $Jord(\psi_{p})$ sont inf\'erieures ou \'egale \`a 2. On note ${\cal E}$ le sous-ensemble de $Jord(\psi)$ form\'e des \'el\'ements dont la multiplicit\'e est exactement 2 et on consid\`ere ${\cal E}$ comme ensemble sans multiplicit\'e. On note aussi $\psi'$ le morphisme qui se d\'eduit de $\psi$ en enlevant les \'el\'ements de ${\cal E}$ de $Jord(\psi)$. On reprend les notations habituelles pour $Jord(\psi)$; on fixe des fonctions $\underline{\ell}',\underline{\eta}'$ telles que $\pi(\psi',\epsilon',\underline{\ell}',\underline{\eta}')$ soit l'une des composantes irr\'eductibles de $\pi(\psi',\epsilon')$. On doit d\'emontrer que l'induite
$$
\times _{(\rho,A,B,\zeta)\in {\cal E}}S(\rho,A,B,\zeta)\times \pi(\psi',\epsilon',\underline{\ell}',\underline{\eta}')
$$
est de longueur exactement $\sum_{(\rho,A,B,\zeta)\in {\cal E}}(A-B+2)$. On fait la d\'emonstration par induction sur le nombre d'\'el\'ement de ${\cal E}$; on note $(\rho,A,B,\zeta)$ l'\'el\'ement de ${\cal E}$ tel que $B$ soit maximum pour cette propri\'et\'e d'appartenance. Avec l'hypoth\`ese de r\'ecurrence appliqu\'e \`a ${\cal E}-\{(\rho,A,B,\zeta)\}$ et \ref{induction}, on a cette assertion \`a condition que pour tout $(\rho,A',B',\zeta')\in Jord(\psi)$ v\'erifiant $B'>B$ alors $B'>>A$. Il faut donc passer de $B'>>A$ \`a $B'>A$. Cela va \^etre cons\'equence du r\'esultat plus g\'en\'eral suivant

\nl
\bf Assertion 1: \sl soit $\psi,\epsilon,\underline{\ell},\underline{\eta}$ permettant de d\'efinir une repr\'esentation $\pi(\psi,\epsilon,\underline{\ell},\underline{\eta})$ dont on suppose qu'elle est non nulle. Soit $(\rho,A,B,\zeta)\in Jord(\psi)$ y intervenant avec multiplicit\'e 1 et tel que pour tout $(\rho,A',B',\zeta')\in Jord(\psi)$ un \'el\'ement diff\'erent, soit $B'>A+1$ soit $A'<B-1$. Alors
$$
Jac_{\zeta B, \cdots, \zeta A}\pi(\psi,\epsilon,\underline{\ell},\underline{\eta})\neq 0.
$$
\rm
Il faut remarquer que les hypoth\`eses de cette assertion entra\^{\i}ne que $(\rho,A,B,\zeta)$ intervient avec multiplicit\'e 1 dans $Jord(\psi)$.

On fait une r\'ecurrence sur $\sum_{(\rho',A',B',\zeta')\in Jord(\psi)}A'-B'$. Le cas des morphismes \'el\'ementaires est connus (cf. \cite{elementaire}) et on n'a d'ailleurs besoin que de l'hypoth\`ese $B'>A$ et non $B'>A+1$. Avant de passer au cas g\'en\'eral, on g\'en\'eralise un autre r\'esultat du cas discret
\nl
\bf Assertion 2: \sl soit $\psi$ ayant les propri\'et\'es du th\'eor\`eme et soit $\underline{\ell}$ comme dans cet \'enonc\'e. On suppose qu'il existe $(\rho,A,B,\zeta)\in Jord(\psi)$ tel que $\underline{\ell}(\rho,A,B,\zeta)>0$. On note $m$ la multiplicit\'e de $(\rho,A,B,\zeta)$ dans $Jord(\psi)$. Alors $Jac_{\zeta B, \cdots, -\zeta A}$ appliqu\'e $m$ fois \`a $\pi(\psi,\epsilon,\underline{\ell},\underline{\eta})$ est une repr\'esentation irr\'eductible non nulle. On obtient exactement $\pi(\psi',\epsilon,\underline{\ell}',\underline{\eta})$, o\`u $\psi'$ se d\'eduit de $\psi$ en rempla\c{c}ant chaque copie de $(\rho,A,B,\zeta)$ par $(\rho,A-1,B+1,\zeta)$ et $\underline{\ell}'$ ne diff\`ere de $\underline{\ell}$ que sur $(\rho,A-1,B+1,\zeta)$ o\`u cette fonction vaut $\underline{\ell}(\rho,A,B,\zeta)-1$.\rm
\nl
On fixe une des copies de $(\rho,A,B,\zeta)$ dans $Jord(\psi)$; si pour tout  \'el\'ement, $(\rho,A',B',\zeta')$ de $Jord(\psi)$ strictement plus grand que cette copie fix\'ee pour l'ordre mis sur $Jord(\psi)$, on a $B'>>A$, il suffit d'appliquer la d\'efinition ( si  $(\rho,A,B,\zeta)$ intervient avec multiplicit\'e, on renvoie \`a \ref{nullitenonnullite} et si cette multiplicit\'e est plus grande que 2 \`a l'irr\'eductibilit\'e \ref{parametrisationetinduction} et \ref{inductionfin}). Pour la situation g\'en\'erale, une part d'une situation comme ci-dessus et on calcule des $Jac_{x\in {\cal C}(\zeta',A',B',T')}$; mais pour tout \'el\'ement $x$ dans un tel ensemble ${\cal C}(\zeta',A',B',T')$, on a s\^urement $$\vert x\vert \geq B'+1\geq A+2.$$
Il y a donc une commutation entre cette op\'eration et $Jac_{\zeta B, \cdots, -\zeta A}$. D'o\`u l'assertion 2.
\nl
Revenons \`a la preuve de l'assertion 1. Gr\^ace \`a l'assertion 2, on se ram\`ene au cas o\`u $\underline{\ell}(\rho,A',B',\zeta')=0$ pour tout $(\rho,A',B',\zeta')\in Jord(\psi)$ diff\'erent de $(\rho,A,B,\zeta)$. La r\'ecurrence porte donc uniquement sur $A-B$ et si $A=B$, on a d\'ej\`a vue cette assertion.
On suppose donc que $A>B$. Et on suppose d'abord que  $\underline{\ell}(\rho,A,B,\zeta)=0$; on est directement ramen\'e au cas d'un morphisme \'el\'ementaire. On suppose maintenant que $\underline{\ell}(\rho,A,B,\zeta)>0$

 %%%%%%
 \section{Appendice}
 \subsection{Un r\'esultat d'irr\'eductibilit\'e\label{unresultatdirreductibilite}}
 
 On rappelle et red\'emontre (pour la commoditi\'e du lecteur) un r\'esultat de \cite{mw}. Pour cela on a besoin d'adapter les notations. Soient $(\rho,A,B,\zeta)$ comme dans les paragraphes pr\'ec\'edents. On fixe $T$ un entier et on consid\`ere le tableau:
 $$
 \begin{matrix}
 &\zeta(B+T) &\cdots & \zeta(B+1)\\
 &\vdots &\vdots &\vdots\\
 &\zeta (A+T) &	\cdots &\zeta (A+1)
 \end{matrix}
 $$
 On a associ\'e une repr\'esentation d'un groupe lin\'eaire convenable $S(\zeta,A,B,T)$, \`a ce tableau. 
 
 {\bf Proposition}: \sl
 Soit $\{T_{1},\cdots, T_{\ell}\}$ une collection d'entiers; l'induite $\times _{i\in [1,\ell]}S(\zeta,A,B,T_{i})$ est irr\'educ\-tible comme repr\'esentation d'un groupe lin\'eaire convenable. 
 
 \
 
 \rm
Soit $\sigma$ un sous-quotient irr\'eductible de cette induite c'est une repr\'esentation d'un groupe lin\'eaire $GL(N)$ pour $N$ convenable. Pour ${\cal C}$ un ensemble totalement ordonn\'e on garde la notation d\'ej\`a utilis\'ee $Jac^d_{x\in {\cal C}}\sigma$ pour repr\'esenter la projection du module de Jacquet de $\sigma$ le long du parabolique standard de Levi $GL(N-\vert {\cal C}\vert d_{\rho})$ $\times GL(d_{\rho}) \times \cdots \times GL(d_{\rho})$ o\`u il y a $\vert {\cal C}\vert$-copies $GL(d_{\rho})$, sur le support cuspidal des copies $GL(d_{\rho})$ qui est pr\'ecis\'ement $\otimes_{x\in {\cal C}}\rho\vert\,\vert^x$ avec l'ordre de ${\cal C}$. On prend l'expression $\zeta$ plus petit pour dire le plus petit si $\zeta=+$ et le plus grand si $\zeta=-$. Le $\zeta$ plus petit \'el\'ement du support cuspidal de $\sigma$ est $\zeta (B+1)$ et il intervient avec multiplicit\'e $\ell$ dans le support cuspidal de l'induite de l'\'enonc\'e et donc aussi de $\sigma$. Il existe donc des \'el\'ements $Y_{i}$ pour $i=1,\cdots, \ell$ tels que $[\zeta(B+1),Y_{i}]$ soit un segment et en notant ${\cal C}$ l'ensemble totalement ordonn\'e
 $
 \cup_{i\in [\ell, 1]}[\zeta(B+1),Y_{i}]$ on ait $$Jac^d_{x\in {\cal C}}\sigma\neq 0.\eqno(1)$$On fixe $Y_{1}$ $\zeta$-minimal avec cette propri\'et\'e puis $Y_{2}$ $\zeta$-minimal, etc... On va montrer que $Y_{i}=\zeta(A+1)$ pour tout $i=1, \cdots, \ell$. Pour $x$ fix\'e, on a $Jac^d_{x}S(\zeta,A,B,T)=0$ quelque soit $T$ sauf pour $x=\zeta(A+1)$. La non nullit\'e de (1) entra\^{\i}ne \`a fortiori que $Jac^d_{Y_{1}}\sigma\neq 0$ et comme ce module de Jacquet est un sous-quotient du module de Jacquet de toute l'induite, il existe $i\in [1,\ell]$ tel que $Jac^d_{Y_{1}}S(\zeta,A,B,T)\neq 0$. D'o\`u d\'ej\`a $Y_{1}=\zeta(A+1)$. On calcule maintenant $Y_{2}$: on v\'erifie par r\'eciprocit\'e de Frobenius qu'il existe une repr\'esentation $\sigma'$ et une inclusion
 $$
 \sigma\hookrightarrow \sigma'\times \rho\vert\,\vert^{\zeta(B+1)} \times \cdots \times \rho\vert\,\vert^{Y_{2}}\times \rho\vert\,\vert^{\zeta (B+1)}\times \cdots \times \rho\vert\,\vert^{\zeta(A+1)}.
 $$
 On peut remplacer $\rho\vert\,\vert^{\zeta (B+1)}\times \cdots \times \rho\vert\,\vert^{\zeta(A+1)}$ par $<\rho\vert\,\vert^{\zeta(B+1)}, \cdots, \rho\vert\,\vert^{\zeta(A+1)}>$ car sinon il existerait $x\in [\zeta(B+1),\zeta(A+1)[$ tel que $Jac^d_{x}\sigma \neq 0$. Par minimalit\'e de $Y_{2}$ on peut aussi remplacer $ \rho\vert\,\vert^{\zeta(B+1)} \times \cdots \times \rho\vert\,\vert^{Y_{2}}$ par $<\rho\vert\,\vert^{\zeta(B+1)}, \cdots, \rho\vert\,\vert^{Y_{2}}>$. La repr\'esentation induite $$<\rho\vert\,\vert^{\zeta(B+1)}, \cdots, \rho\vert\,\vert^{Y_{2}}>\times
 <\rho\vert\,\vert^{\zeta(B+1)}, \cdots, \rho\vert\,\vert^{\zeta(A+1)}>$$est irr\'eductible car les segments ont m\^eme propri\'et\'e de croissance et m\^eme premier terme. Cette irr\'eductibilit\'e permet d'\'echanger les 2 facteurs et prouve que $Jac^d_{Y_{2}}\sigma \neq 0$; d'o\`u comme ci-dessus $Y_{2}=\zeta(A+1)$ et l'argument s'applique de proche en proche pour prouver que $Y_{i}=\zeta(A+1)$ pour tout $i\in [1,\ell]$. 
 
 On note ici ${\cal C}$ l'union de $\ell$-copies de $[\zeta(B+1),\zeta(A+1)]$ et $c$ le cardinal de cet ensemble ; ici l'ensemble n'est plus ordonn\'e.

 On calcule le module de Jacquet de $\sigma$ relativement au parabolique standard de Levi $GL(N-cd_{\rho})$ $ \times$ $ GL(cd_{\rho})$ et on prend encore la projection de ce module de Jacquet sur l'ensemble des rep\'esentations pourlesquelles $GL(cd_{\rho})$ agit par une repr\'esentation dont le support cuspidal est form\'e des repr\'esenta\-tions $\rho\vert\,\vert^x$ o\`u $x$ parcourt ${\cal C}$. On note $Jac_{[{\cal C]}}\sigma$ ce module de Jacquet. Par exactitude du foncteur de Jacquet, ceci est un sous-quotient de son analogue quand on remplace $\sigma$ par toute l'induite $\tau:=\times_{i\in [1,\ell]}S(\zeta,A,B,T_{i})$. On pose  $\tau':=\times_{i\in [1,\ell]}S(\zeta,A+1,B+1,T_{i}-1)$, c'est-\`a-dire que  l'on enl\`eve les derni\`eres lignes de chaque tableau d\'efinissant $\tau$. On v\'erifie par les formules standard que $$Jac_{[{\cal C}]}
\tau=\tau'\otimes \bigl (<\rho\vert\,\vert^{\zeta(B+1)}, \cdots, \rho\vert\,\vert^{\zeta(A+1)}>\times \cdots \times <\rho\vert\,\vert^{\zeta(B+1)}, \cdots, \rho\vert\,\vert^{\zeta(A+1)}>\bigr),$$
 o\`u il y a $\ell$ copies cach\'ees dans les $\cdots$. Un argument par r\'ecurrence dit que $\tau'$ est irr\'eductible ce qui entra\^{\i}ne l'irr\'eductibilit\'e du module de $Jac_{[{\cal C}]}\tau$; en particulier $Jac_{[{\cal C}]}\tau=Jac_{[{\cal C}]}\sigma$. On peut aussi calculer $Jac_{[{\cal C]}}\tau''$ o\`u  $$\tau'':=\tau'\times <\rho\vert\,\vert^{\zeta(B+1)}, \cdots, \rho\vert\,\vert^{\zeta(A+1)}>\times \cdots \times <\rho\vert\,\vert^{\zeta(B+1)}, \cdots, \rho\vert\,\vert^{\zeta(A+1)}>.$$ C'est tr\`es facile car si $x$ est tel que $Jac^d_{x}\tau'\neq 0$ alors $x=\zeta(A+2)$; or $\zeta(A+2)\notin {\cal C}$. Ainsi, on a aussi $Jac_{[{\cal C}]}\tau''=Jac_{[{\cal C}]}\sigma$. On en d\'eduit que $\tau''$ 
 a un unique sous-module irr\'eductible qui intervient avec multiplicit\'e 1 comme sous-quotient irr\'eductible et qui est n\'ecessairement $\sigma$. Cela d\'etermine uniquement $\sigma$ et prouve que $\sigma$ intervient a fortiori avec multiplicit\'e 1 dans $\tau$. Donc $\tau=\sigma$. Cela termine la preuve de l'irr\'eductibilit\'e. 

 %On peut certainement g\'en\'eraliser
 
 \subsection{Un deuxi\`eme r\'esultat d'irr\'eductibilit\'e\label{undeuxiemeresultatdirreductibilite}}
On reprend la notation $S_{-}$ de \ref{induction} pour la repr\'esentation bas\'ee sur $\rho$ et associ\'ee au tableau 
$$
\begin{matrix}
&\zeta B &\cdots &-\zeta A\\
&\vdots &\vdots &\vdots \\
&\zeta (A-1) &\cdots &-\zeta (B+1)
\end{matrix}
$$
\bf Lemme: \sl La repr\'esentation $S_{-}\times <\rho\vert\,\vert^{\zeta B}, \cdots, \rho\vert\,\vert^{-\zeta A}>$ est irr\'eductible.\rm

\

\noindent
Soit $\tau$ un sous-quotient irr\'eductible de cette induite; le support cuspidal de $\tau$ contient avec multiplicit\'e $2$ la repr\'esentation $\rho\vert\,\vert^{-\zeta A}$; comme dans la preuve ci-dessus, il existe $x_{1},x_{2}$ tels que pour $i=1,2$, $x_{i}+\zeta A\in \mathbb Z$ et il existe une repr\'esentation irr\'eductible $\tau'$ avec une inclusion:
$$
\tau \hookrightarrow <\rho\vert\,\vert^{x_{1}}, \cdots, \rho\vert\,\vert^{-\zeta A}> \times <\rho\vert\,\vert^{x_{2}}, \cdots, \rho\vert\,\vert^{-\zeta A}> \times \tau'.
$$
Et on ne peut encore avoir que $x_{1}=x_{2}=\zeta B$ car n\'ecessairement $Jac_{x_{i}}\tau\neq 0$ pour $i=1,2$ et par exactitude du foncteur de Jacquet, aussi $Jac_{x_{i}}(S_{-}\times <\rho\vert\,\vert^{\zeta B}, \cdots, \rho\vert\,\vert^{-\zeta A}>)\neq 0$. Dans ce cas, il n'y a que la possibilit\'e $\tau'=S(\rho,A-1,B+1,\zeta)$ par un calcul facile et $\tau$ est uniquement d\'etermin\'e. Cela montre l'unicit\'e de $\tau$ et l'irr\'eductibilit\'e.

 %%%%%%%%%%%%%%%%% 

 %%%%%%%%%%%%%%%%


\begin{thebibliography}{1}
\bibitem{arthur} \sc Arthur J.: \sl Unipotent automorphic
representations: conjectures \rm  in Orbites unipotentes et
repr\'esentations II, Ast\'erisque 171-172, 1989, pp. 13-72

\bibitem{arthurnouveau}\sc Arthur J.: \sl An introduction to the trace formula \rm pr\'epublication 

\bibitem{aubert} \sc Aubert A.-M. : \sl Dualit\'e dans le
groupe de Grothendieck de la cat\'egorie des
repr\'esentations lisses de longueur finie d'un groupe
r\'eductif p-adique, \rm TAMS, 347, 1995, pp. 2179-2189 avec
l'erratum publi\'e dans TAMS, 348, 1996, pp. 4687-4690



\bibitem{bernsteinzelevinsky} \sc Bernstein I. N.,
Zelevinsky A. V.: \sl Induced Representations of Reductive
p-adic groups 1 \rm Ann de l'ENS, 10, 1977, pp. 147-185

\bibitem{debacker} \sc Debacker S., Reeder M.: \sl Depth-zero supercuspidal L-packets and their stability \rm pr\'epublication 2004

\bibitem{gangurevich} \sc Gan W.T., Gurevich N.: \sl Non-tempered Arthur Packets of $G_{2}$\rm in 
Automorphic Representations, 
L-Functions and Applications: Progress and Prospects, W de G, OSU math. res. inst. pub. 11, 2005,  pp.129-156


\bibitem{harris} \sc Harris, M.; Taylor, R.: \sl The geometry and
cohomology of some simple Shimura varieties,\rm  Annals of Math
Studies, 151, Princeton Univ. Press, 2001

\bibitem{henniart} \sc Henniart, G.: \sl Une preuve simple des
conjectures de Langlands pour $GL_n$ sur un corps p-adique,\rm
Invent. Math., 139, 2000, pp. 439-455

\bibitem{kazhdan} \sc Kazhdan D., Varshavsky Y.: \sl Endoscopic decomposition of 
characters of certain cuspidal representations, \rm ERA,  2004

\bibitem{laumonngo} \sc Laumon G., Ngo B. -C.: \sl Le lemme fondamental pour les 
groupes unitaires, \rm 
arXiv AG/ 0404454




\bibitem{elementaire} \sc M{oe}glin C.: \sl
Sur certains paquets d'Arthur et involution
d'Aubert-Schneider-Stuhler g\'e\-n\'e\-ralis\'ee, \rm representation theory, volume 10, 2006, pp. 86-129

\bibitem{discretorthogonaux} \sc M{oe}glin C.: \sl Classification  des s\'eries discr\`etes pour 
certains groupes classiques $p$-adiques, \rm pr\'epublication, 2006

\bibitem{discretunitaire}\sc M{oe}glin C.: \sl Classification et Changement de base pour les  
s\'eries discr\`etes des groupes unitaires 
p-adiques, \rm pr\'epublication, 2006

\bibitem{paquetsdiscrets}\sc M{\oe}glin C.: \sl Paquets d'Arthur discrets pour les groupes classiques, \rm 
pr\'epublication

\bibitem{europe} \sc M{\oe}glin C.: \sl Classification des
s\'eries discr\`etes: param\`etres de Langlands et
exhaustivit\'e, \rm JEMS, 4, 143-200, 2002

\bibitem{ams} \sc M{\oe}glin C., Tadic M.: \sl Construction
of discrete series for classical p-adic groups, \rm journal
de l'AMS, volume 15, 2002, pp 715-786

\bibitem{inventiones} \sc M{\oe}glin C., Waldspurger J.-L.:
\sl 
Paquets stables de repr\'esentations
temp\'er\'ees et de r\'eduction unipotente pour $SO(2n+1)$,
\rm Inventiones, 152, ,  2003, 461-623

\bibitem{selecta} \sc M{\oe}glin C., Waldspurger J.-L.:
\sl Sur le transfert des traces d'un groupe classique p-adique \`a un groupe lin\'eaire 
tordu
\rm pr\'epublication 2006, \`a para\^{\i}tre dans Selecta Math.

\bibitem{mw}  \sc M{\oe}glin C., Waldspurger J.-L.:
\sl Spectre r\'esiduel de GL(n), \rm Ann de l'ENS, 22, 1989, pp. 605-674

\bibitem{ngo} \sc Ngo B. -C.: \sl Fibration de Hitchin et endoscopie, \rm Inventiones, 164, 2006, pp. 399-453


\bibitem{silberger}\sc Silberger A. : \sl Special representations of reductive p-adic groups are 
not integrable, \rm Ann. of Math., 111, 1980, pp. 571-587.

\bibitem{SS}\sc Schneider M., Stuhler U. : \sl Representation
theory and sheaves on the Bruhat-Tits building \rm Publ. Math.
IHES 85, 1997, pp. 97-191

\bibitem{JIMJ}\sc Waldspurger J.-L.: \sl La formule de
Plancherel d'apr\`es Harish-Chandra, \rm   JIMJ, vol2, fasc 2,
2003


\bibitem{waldspurger} \sc Waldspurger J.-L.: \sl Le groupe GL$_{N}$ tordu sur un corps 
$p$-adique, 
1e et 2e partie, \rm pr\'epublication 2005, \`a para\^{\i}tre au Duke math. 



\bibitem{waldspurgerJIMJ} \sc Waldspurger J. -L.: \sl Endoscopie et Changement de 
caract\'eristique, \rm JIMJ, 5.3, 2006, pp. 423-525

\bibitem{zelevinsky} \sc 
Zelevinsky A. V.: \sl Induced Representations of Reductive
p-adic groups II, \rm Ann de l'ENS, 13, 1980, pp. 165-210
\end{thebibliography}
\end{document}